\documentclass[compress,preprint]{elsarticle}

\usepackage[abs]{overpic}
\usepackage{graphicx}
\usepackage{stmaryrd}
\usepackage{amsfonts, color, amssymb}
\usepackage{amsmath, bm}
\usepackage{amsthm}
\usepackage{epsfig}
\usepackage{psfrag}
\usepackage{paralist}
\usepackage{algorithm}
\usepackage{algorithmicx}
\usepackage{algpseudocode}
\usepackage{appendix}
\usepackage{caption}
\usepackage{hyperref}
\usepackage{changes}
\usepackage{ulem}
\usepackage{subcaption}

\usepackage[misc]{ifsym}
\usepackage{placeins}

\usepackage[a4paper, total={6in,8in}]{geometry}

\allowdisplaybreaks[4]
\newtheorem{theorem}{Theorem}[section] %
\newtheorem{definition}{Definition}[section] %

\newtheorem{example}{Example}[section]

\newtheorem{lemma}{Lemma}[section]

\newcommand{\bu}{{\bf u}}

\def\bn{{\bf n}}
\def\3bar{{|\hspace{-.02in}|\hspace{-.02in}|}}
\allowdisplaybreaks[4]

\newcommand{\trb}[1]{|\!|\!|#1|\!|\!|}

\allowdisplaybreaks 

\numberwithin{equation}{section}

\journal{Computers $\&$ Mathematics with Applications}

\begin{document}

\begin{frontmatter}

\title{A weak Galerkin finite element method for the coupled Stokes--Biot problem}


\author[1]{Rui Li}
\ead{lirui25@mails.jlu.edu.cn}


\author[2]{\texorpdfstring{Hui Peng\corref{cor1}}{Hui Peng}}
\ead{penghui@njust.edu.cn}
\cortext[cor1]{Corresponding author.}

\author[1]{Xiuli Wang}
\ead{xiuli19@jlu.edu.cn}

\address[1]{School of Mathematics, Jilin University, Changchun, 130012, Jilin, P. R. China}
\address[2]{School of Mathematics and Statistics, Nanjing University of Science and Technology, Nanjing, 210094, Jiangsu, P. R. China}

\begin{abstract}	
In this paper, we consider a stable numerical method for the Stokes--Biot model. We adopt the standard formulation of the Stokes equations coupled with the mixed formulation of the Biot system involving the additional total pressure variable. We construct the semi-discrete scheme on polygonal meshes by employing weak Galerkin finite element methods for spatial discretization. It is convenient to incorporate the interface conditions into the semi-discrete scheme via the edge-based degrees of freedom. Furthermore, we give a fully discrete scheme based on the backward Euler method. We establish the well-posedness of these two schemes and provide a rigorous error analysis. Finally, we present some numerical examples to verify the convergence and stability of our algorithm under a wide range of some physical parameters. 
\end{abstract}


\end{frontmatter}


\section{Introduction}
\label{section:introduction}
The interaction between a free fluid and poroelastic medium frequently occurs in scientific and engineering applications, such as blood flow in arteries, surface water-groundwater exchange, and wave-induced seabed responses \cite{MR1936102-2001-ANM-background2,MR1885704-2001-Siam-background1,MR3475830-2016-IJAESCAE-background1-1}. 
In order to simulate the complex behavior, the Stokes--Biot model is proposed and widely investigated. However, the analytical solution is hard to obtain, many scholars turn their attention to numerical solutions.

During the past few years, 
various numerical methods have been {\color{black}developed} for this coupled model, including mixed finite element methods \cite{MR4487559-2022-NM-mixed,MR4362198-2022-Esaim-mixed}, discontinuous Galerkin methods \cite{MR4598417-2023-CMA-total,MR4124408-2020-CMA-primal,MR4766788-2024-JCAM-total}, virtual element methods \cite{MR4941775-2025-Esaim-three-VEM}, 
 diffuse interface methods \cite{MR5024699-2026-IMAJNA-primal} and so on. Within these existing studies, much focus lies on two key numerical issues.  
 On the one hand, the appropriate treatment of interface conditions is essential, and several approaches have been developed for this purpose.
 For instance, a Lagrange multiplier method was proposed in \cite{MR3851065-2018-NM-mixed-Lagrange} 
 to weakly enforce {\color{black}the continuity of the normal velocity across the interface},  
 while Nitsche's method was adopted in \cite{MR3347244-2015-CMAME-nitsche} to impose the interface conditions within a consistent and stable variational framework. 
 On the other hand, numerical instabilities require attention, namely Poisson locking and spurious pressure oscillations in the poroelastic region.
 Researchers have made several attempts to mitigate these numerical instabilities. Stabilization techniques were considered in \cite{MR2455881-2008-stabilization,MR3427711-2016-stabilization}, where artificial diffusion terms were added to the mass conservation equation to obtain stable numerical solutions.  
 In addition, reformulating the governing equations by introducing additional unknowns is another common strategy. Typical examples include mixed formulations that take the Darcy flux as an independent unknown \cite{MR3820123-2018-three-darcy-flux3,MR2327964-2007-three-darcy-flux2,MR2327966-2007-three-darcy-flux1} and the three-field formulation with total pressure \cite{MR3590654-2017-total,MR3552204-2016-Siam-total-pressure}.  
 Previous works have demonstrated that these formulations improve the robustness of numerical approximations for single-field Biot problems, and hence possess the potential to be applied to the coupled Stokes–Biot system.

The aim of this paper is to develop a stable weak Galerkin (WG) finite element method for the coupled system consisting of standard Stokes equations and the mixed displacement-pressure-total pressure formulation of Biot system. The WG method was proposed in \cite{MR2994424-WG-first}. The key point of this method lies in substituting  classical differential operators with weakly‑defined ones for discontinuous approximating functions in the numerical scheme. Since it is straightforward to construct  
{\color{black}arbitrary‑order elements} and 
easy to implement on polygonal meshes, {\color{black}WG methods have been widely applied to} various problems,  
{\color{black}such as} elasticity equations \cite{MR3816361-elasticity1,MR4489858-elasticity2}, Brinkman equations \cite{MR3227259-brinkman1,MR3501216-brinkman3}, biharmonic equations \cite{MR4215881-biharmonic3,MR3190348-biharmonic2}, Maxwell equations \cite{MR3394450-maxwell1,MR4473646-maxwell2},  eigenvalue problems \cite{MR3968996-eigenvalue2,MR3919912-eigenvalue1} and so on. In particular, WG methods have also been separately developed for
Stokes equations \cite{MR4346856-2022-SFWG-AAMM,MR3452926-WG-stokes-2016,MR3475512-2016-WG-JCAM} and Biot's consolidation models \cite{MR3717602-2018-JCAM-WG-biot, MR4588484-2023-ANM, MR3775100-2018-CMA-WG-biot,MR4799262-2025-WG-Biot}, providing stable and accurate discretizations for the individual single-physics problems.

On the basis of these existing studies, we extend the WG method to spatial discretization of the coupled Stokes--Biot system. 
One benefit of the three-field formulation with total pressure is that we only need to consider the compatibility between the discrete spaces of displacement and total pressure \cite{MR3590654-2017-total}.  Since WG methods employ totally discontinuous approximating functions, it is not difficult to achieve the requirement of the compatibility. 
The same element pair is also available for {\color{black}the velocity variable and pressure variable in Stokes equations}.
In addition, the employment of  WG elements makes it convenient to treat the interface conditions via the degrees of freedom defined on edges of the mesh. Moreover, WG methods support polygonal mesh partitions, so meshes of the two subdomains do not need to match along the interface.  As for temporal discretization, we use the backward Euler method. 
The semi-discrete scheme and fully discrete scheme are constructed. We prove the existence and uniqueness of the solutions of these two discrete schemes. Furthermore, we give error equations and derive error estimates.
{\color{black}Both our theoretical and numerical results show robustness against Poisson locking and spurious‑pressure oscillations.}
 
{\color{black}An outline of the paper goes as follows}. In Section~\ref{model-weak-form}, we describe the Stokes--Biot model and its weak formulation. In Section~\ref{Weak Galerkin finite element discretization}, we introduce the WG method, and give the semi-discrete and fully discrete schemes. In Section~\ref{section:exist-unique}, we establish the well-posedness of the numerical schemes. The optimal order error estimates are derived in Section~\ref{section:error-estimate}. Some numerical examples are presented in Section~\ref{Section:example} to demonstrate the theoretical results. Finally, in Section \ref{sec:7}, {\color{black}we offer a conclusion}.

\section{Model problem and weak form}\label{model-weak-form}
\begin{figure}
	\centering
	\setlength{\unitlength}{1bp}%
	\begin{picture}(125.31, 186.33)(0,0)
		\put(0,0){\includegraphics{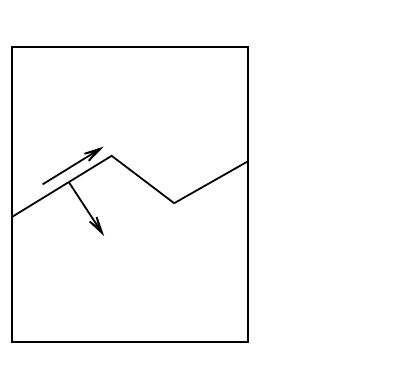}}
		\put(56.22,169.55){\fontsize{14.23}{17.07}\selectfont $\Gamma_f$}
		\put(56.12,8.73){\fontsize{14.23}{17.07}\selectfont $\Gamma_m$}
		\put(84.11,144.92){\fontsize{14.23}{17.07}\selectfont $\Omega_f$}
		\put(19.92,44.55){\fontsize{14.23}{17.07}\selectfont $\Omega_m$}
		\put(29.80,70.32){\fontsize{14.23}{17.07}\selectfont $\mathbf{n}_{\Gamma}$}
		\put(21.87,115.91){\fontsize{14.23}{17.07}\selectfont $\bm{\tau}_{\Gamma}$}
		\put(88.01,100.90){\fontsize{14.23}{17.07}\selectfont $\Gamma$}
	\end{picture}%
	\caption{The model domain of the coupled Stokes--Biot problem.}\label{domain}
\end{figure}

 Let $\Omega \subset \mathbb{R}^2$ be a bounded domain with Lipschitz boundary. The domain $\Omega$ is divided into two parts: a free-flow region $\Omega_f$ and a poroelastic medium region $\Omega_m$. The interface between two subdomains is denoted by $\Gamma = \partial\Omega_f \cap \partial\Omega_m$. The unit outward normal vector and unit tangent vector on the interface are denoted by $\mathbf{n}_{\Gamma}$ and $\bm{\tau}_{\Gamma}$, respectively. As shown in Figure~\ref{domain}, the outer boundaries are defined as $\Gamma_f = \partial \Omega_f \setminus \Gamma$ and $\Gamma_m = \partial \Omega_m\setminus\Gamma$. 

In $\Omega_f$, the Stokes equations are employed to describe the dynamics of the free flow:
\begin{align}
	\label{problem-eq1-s}
	-\nabla\cdot\bm{\sigma}_f(\mathbf{u}_f,p_f)&=\mathbf{f}_f,\quad\text{in}~\Omega_f\times\left(0,T\right] ,
	\\\label{problem-eq2-s}
	\nabla\cdot \mathbf{u}_f&=0,\quad\ \text{in}~\Omega_f\times\left(0,T\right] ,
\end{align}
\textcolor{black}{
where $\mathbf{u}_f$ and $p_f$ represent the fluid velocity and pressure, respectively, $\mathbf{f}_f$ is the prescribed body force and $T>0$ is the final time. The stress tensor is given by $\bm{\sigma}_f(\mathbf{u}_f,p_f)
=2\mu_f\mathbf{D}(\mathbf{u}_f)-p_f\mathbf I,$  $\mu_f$ denotes the fluid viscosity, $\mathbf{D}(\mathbf{u}_f)
=\frac12(\nabla\mathbf{u}_f+(\nabla\mathbf{u}_f)^T)$
denotes the deformation tensor, and $\mathbf I$ is the identity matrix.}

The deformation and fluid transport in $\Omega_m$ are described by the following three-field Biot system:
\begin{align}
	\label{problem-eq1-phi-b}
	-\nabla\cdot\bm{\sigma}_m(\mathbf{u}_m,\phi_m)&=\mathbf{f}_m,\quad\text{in}~\Omega_m\times\left(0,T\right] ,
	\\\label{problem-eq2-phi-b}
	\Bigl(c_0+\dfrac{\alpha^2}{\lambda}\Bigr)\partial_t  p_m-\dfrac{\alpha}{\lambda}\partial_t \phi_m-\nabla\cdot(\kappa\nabla p_m)&=g,\quad\ \text{in}~\Omega_m\times\left(0,T\right] ,
	\\\label{problem-eq3-phi-b}
	\dfrac{1}{\lambda}\phi_m-\dfrac{\alpha}{\lambda}p_m+\nabla\cdot\mathbf{u}_m&=0,\quad\ \text{in}~\Omega_m\times\left(0,T\right],
\end{align}
where $\mathbf{u}_m$ represents the displacement of the porous skeleton, $p_m$ is fluid pressure, $\phi_m$ is the total pressure. The stress is defined by $\bm{\sigma}_m(\mathbf{u}_m,\phi_m)=2\mu_m\mathbf{D}(\mathbf{u}_m)-\phi_m\mathbf{I}$. The functions $\mathbf{f}_m$ and $g$ denote the prescribed body force and source term. The parameters $\alpha$, $c_0\geq0$ represent the Biot--Willis coefficient and the constrained specific storage coefficient, respectively.  $\kappa$ denotes the permeability of the porous medium. The coefficients $\mu_m$ and $\lambda$ are Lam\'e parameters, which are determined  by Young's modulus $E$ and Poisson's ratio $0<\gamma<\frac{1}{2}$ through $\mu_m=\frac{E}{2(1+\gamma)},\ \lambda=\frac{E\gamma}{(1+\gamma)(1-2\gamma)}.$

On the interface, four interface conditions \cite{MR2573342-interface-condition1,MR4362198-2022-Esaim-mixed} are imposed:
\begin{align}
	\label{problem-interface1}
	\mathbf{u}_f\cdot\mathbf{n}_{\Gamma}&=\left(\partial_t \mathbf{u}_m-\kappa\nabla p_m\right)\cdot\mathbf{n}_{\Gamma},\quad\ \text{on}~\Gamma\times\left(0,T\right],
	\\\label{problem-interface2}
	-\bm{\sigma}_f(\mathbf{u}_f,p_f)\mathbf{n}_{\Gamma}\cdot\mathbf{n}_{\Gamma}&=p_m,\quad\ \text{on}~\Gamma\times\left(0,T\right],
	\\\label{problem-interface3}
	\bm{\sigma}_f(\mathbf{u}_f,p_f)\mathbf{n}_{\Gamma}&=\bm{\sigma}_m(\mathbf{u}_m,\phi_m)\mathbf{n}_{\Gamma},\quad\ \text{on}~\Gamma\times\left(0,T\right],
	\\\label{problem-interface4}
	-\bm{\sigma}_f(\mathbf{u}_f,p_f)\mathbf{n}_{\Gamma}\cdot\bm{\tau}_{\Gamma}&=\beta_{\Gamma}\left(\mathbf{u}_f-\partial_t \mathbf{u}_m\right)\cdot\bm{\tau}_{\Gamma},\quad\ \text{on}~\Gamma\times\left(0,T\right],
\end{align}
\textcolor{black}{
The coefficient $\beta_{\Gamma}$ is given by
$\beta_{\Gamma}=\mu_f\alpha_{\rm BJS}\kappa^{-1/2},$
and $\alpha_{\rm BJS}$ is a coefficient determined by experiment.}

On the outer boundary, we consider the following boundary conditions:
\begin{align}
	&\label{problem-boundary1}
	\mathbf{u}_f=\mathbf{0},\quad\ \text{on}~\Gamma_f^D\times\left(0,T\right],& &\bm{\sigma}_f(\mathbf{u}_f,p_f)\mathbf{n}=0,\quad\ \text{on}~\Gamma_f^N\times\left(0,T\right],
	\\\label{problem-boundary2}
	&\mathbf{u}_m=\mathbf{0},\quad\ \text{on}~\Gamma_m^{D}\times\left(0,T\right],& &\bm{\sigma}_m(\mathbf{u}_m,\phi_m)\mathbf{n}=0,\quad\ \text{on}~\Gamma_m^{N}\times\left(0,T\right],
	\\\label{problem-boundary3}
	&p_m=0,\quad\ \text{on}~\Gamma_m^{D_1}\times\left(0,T\right],& &(-\kappa\nabla p_m)\cdot\mathbf{n}=0,\quad\ \text{on}~\Gamma_m^{N_1}\times\left(0,T\right],
\end{align}
where $\mathbf{n}$ is the unit outward normal vector, $|\Gamma_f^D|>0$ and $|\Gamma_f^N|>0$ are parts of the boundary of Stokes subdomain and satisfy $\Gamma_f=\Gamma_f^D\cup\Gamma_f^N$,  $\Gamma_f^D\cap\Gamma_f^N=\emptyset$. $|\Gamma_m^{D}|>0$ and $|\Gamma_m^{N}|>0$ consist of a partition of $\Gamma_m$:  $\Gamma_m=\Gamma_m^{D}\cup\Gamma_m^{N}$ and $\Gamma_m^{D}\cap\Gamma_m^{N}=\emptyset$, the other partition of $\Gamma_m$ is  $\Gamma_m=\Gamma_m^{D_1}\cup\Gamma_m^{N_1}$ and $\Gamma_m^{D_1}\cap\Gamma_m^{N_1}=\emptyset$.

The initial conditions are given as
\begin{align}
	\label{problem-initial}
	\mathbf{u}_m(\cdot,0)=\mathbf{u}_m^0,\qquad\ p_m(\cdot,0)=p_m^0,\qquad\ \phi_m(\cdot,0)=\phi_m^0,
\end{align}
where $\phi_m^0=-\lambda\nabla\cdot\mathbf{u}_m^0+\alpha p_m^0$.
\par
For the Stokes and Biot subdomains, we define the following Sobolev spaces:
\begin{align*}
	V^i&=\{\mathbf{v}_i\in[H^1(\Omega_i)]^2,\,\mathbf{v}_i|_{\Gamma_i^D}=\mathbf{0}\},\ W^i=L^2(\Omega_i),\quad i=f,m.\\
	Z&=\{q_m\in H^1(\Omega_m),\,q_m|_{\Gamma_m^{D_1}}=0\}.
\end{align*}
Then, we present the weak formulation of the coupled problem \eqref{problem-eq1-s}--\eqref{problem-initial}: Seek $\mathbf{u}_f\in V^f,\, p_f\in W^f,\,\mathbf{u}_m\in V^m,\,p_m\in Z$, and $\phi_m\in W^m$ such that
\begin{align*}
	&2\mu_f(\mathbf{D}(\mathbf{u}_f),\mathbf{D}(\mathbf{v}_f))_{\Omega_f}-(p_f,\nabla\cdot\mathbf{v}_f)_{\Omega_f}+\langle\beta_{\Gamma}\left(\mathbf{u}_f-\partial_t \mathbf{u}_m\right)\cdot\bm{\tau}_{\Gamma},\mathbf{v}_f\cdot\bm{\tau}_{\Gamma}\rangle_{\Gamma}\nonumber\\
	&\qquad+\langle p_m\mathbf{n}_{\Gamma},\mathbf{v}_f\rangle_{\Gamma}=(\mathbf{f}_f,\mathbf{v}_f)_{\Omega_f},\\
	&(\nabla\cdot\mathbf{u}_f,q_f)_{\Omega_f}=0,\\
	&2\mu_m(\mathbf{D}(\mathbf{u}_m),\mathbf{D}(\mathbf{v}_m))_{\Omega_m}-(\phi_m,\nabla\cdot\mathbf{v}_m)_{\Omega_m}-\langle\beta_{\Gamma}\left(\mathbf{u}_f-\partial_t \mathbf{u}_m\right)\cdot\bm{\tau}_{\Gamma},\mathbf{v}_m\cdot\bm{\tau}_{\Gamma}\rangle_{\Gamma}\nonumber\\
	&\qquad-\langle p_m\mathbf{n}_{\Gamma},\mathbf{v}_m\rangle_{\Gamma}=(\mathbf{f}_m,\mathbf{v}_m)_{\Omega_m},\\
	&\Bigl(c_0+\frac{\alpha^2}{\lambda}\Bigr)(\partial_t  p_m,q_m)_{\Omega_m}-\frac{\alpha}{\lambda}(\partial_t \phi_m,q_m)_{\Omega_m}+(\kappa\nabla p_m,\nabla q_m)_{\Omega_m}+\langle\partial_t \mathbf{u}_m\cdot\mathbf{n}_{\Gamma},q_m\rangle_{\Gamma}\nonumber\\
	&\qquad-\langle\mathbf{u}_f\cdot\mathbf{n}_{\Gamma},q_m\rangle_{\Gamma}=(g,q_m)_{\Omega_m},\\
	&\frac{1}{\lambda}(\phi_m,\psi_m)_{\Omega_m}-\frac{\alpha}{\lambda}(p_m,\psi_m)_{\Omega_m}+(\nabla\cdot\mathbf{u}_m,\psi_m)_{\Omega_m}=0,
\end{align*}
for all $\mathbf{v}_f\in V^f,\,q_f\in W^f,\,\mathbf{v}_m\in V^m,\,q_m\in Z$ and $\psi_m\in W^m$.

\section{Weak Galerkin finite element discretization}\label{Weak Galerkin finite element discretization}
In this section, \textcolor{black}{
we begin by defining the WG finite element spaces and subsequently give the semi-discrete and fully discrete formulations for the coupled Stokes--Biot system.}

Let $\mathcal{T}_h^f$ and $\mathcal{T}_h^m$ be the shape regular polygonal partitions  \cite{MR3223326-shape-regular-inequality} of $\Omega_f$ and $\Omega_m$, respectively. 
Let $\mathcal{T}_h =\mathcal{T}_h^f\cup\mathcal{T}_h^m$.
For each element $T\in\mathcal{T}_h$, \textcolor{black}{ 
let $|T|$ and $h_T$ denote the area and diameter of $T$, respectively. The mesh size of $\mathcal{T}_h$ is given by $h=\max\limits_{T\in\mathcal{T}_h}h_T$.} 
\textcolor{black}{
The set of all edges in $\mathcal{T}_h$ is denoted by $\mathcal{E}_h$, and the subsets of edges belonging to $\mathcal{T}_h^f$ and $\mathcal{T}_h^m$ are denoted by $\mathcal{E}_h^f$ and $\mathcal{E}_h^m$, respectively. The set of edges located on $\Gamma$ is denoted by $\mathcal{E}_h^I$. We note that $\mathcal{T}_h^f$ and $\mathcal{T}_h^m$ are not necessary to align on the interface $\Gamma$. For any integer $k\geq1$, $P_k(T)$ represents the polynomial space consisting of all polynomials with degree no more than $k$ on $T$.}
\par
\textcolor{black}{
First, we define the WG spaces associated with the velocity and displacement as}
\begin{align*}
	V_h^i&=\{\mathbf{v}_{i,h}=\{\mathbf{v}_{i,0},\mathbf{v}_{i,b}\}:\mathbf{v}_{i,0}|_T\in[P_k(T)]^2,\,\forall\,T\in\mathcal{T}^i_h;\,\mathbf{v}_{i,b}|_e\in[P_k(e)]^2,\,\forall\,e\in\mathcal{E}^i_h\},\quad i=f,m.
\end{align*}
The subspaces satisfying the homogeneous Dirichlet boundary condition on $\Gamma_i^D, i=f,m$ are given by
\begin{align*}
	V_h^{i,D}&=\{\mathbf{v}_{i,h}\in V_h^i:\mathbf{v}_{i,b}|_{\Gamma^D_i}=\mathbf{0}\}.
\end{align*} 
Define the discrete spaces for the fluid pressure in $\Omega_f$ and the total pressure in $\Omega_m$.
\begin{align*}
	W_h^i&=\{q_{i,h}\in L^2(\Omega_i):q_{i,h}|_T\in P_{k-1}(T),\,\forall\,T\in\mathcal{T}_h^i\},\quad i=f,m.
\end{align*}
\textcolor{black}{
For the pore pressure variable in $\Omega_m$, we introduce the following scalar-valued WG space:}
\begin{align*}
	Z_h&=\{q_{m,h}=\{q_{m,0},q_{m,b}\}:q_{m,0}|_T\in P_k(T),\,\forall\,T\in\mathcal{T}^m_h;\, q_{m,b}|_e\in P_k(e),\,\forall\,e\in\mathcal{E}^m_h\}.
\end{align*}
The corresponding subspace satisfying the homogeneous Dirichlet boundary condition on $\Gamma_m^{D_1}$ is defined by
\begin{align*}
	Z_h^D&=\{q_{m,h}\in Z_h:q_{m,b}|_{\Gamma_m^{D_1}}=0\}.
\end{align*}

\textcolor{black}{ 
For weak functions with discontinuous interior and boundary components, we introduce} the discrete weak differential operators locally on each element. 
\begin{definition}
	For each $T\in \mathcal{T}_h^i$ and each $\mathbf{v}_{i,h}\in V^i_h,\,i=f,m$, the discrete weak gradient $\nabla_w\mathbf{v}_{i,h}$ is \textcolor{black}{
    determined by} the unique \textcolor{black}{tensor-valued} polynomial in $[P_{k-1}(T)]^{2\times 2}$ such that
	\begin{equation}
		(\nabla_w \mathbf{v}_{i,h},\bm{\sigma})_T=-(\mathbf{v}_{i,0},\nabla\cdot\bm{\sigma})_T+\left\langle\mathbf{v}_{i,b},\bm{\sigma}\bn\right\rangle_{\partial T},\quad \forall\, \bm{\sigma} \in \left[P_{k-1}(T)\right]^{2 \times 2}.\label{define-gradient1}
	\end{equation}
\end{definition}
\textcolor{black}{
Based on the definition of the discrete weak gradient, the weak strain tensor is given by} \[\mathbf{D}_w(\mathbf{v}_{i,h})=\frac{1}{2}(\nabla_w\mathbf{v}_{i,h}+(\nabla_w\mathbf{v}_{i,h})^T).\]
\begin{definition}
	\textcolor{black}{
    For each $T\in \mathcal{T}_h^i$ and each $\mathbf{v}_{i,h}\in V^i_h,\,i=f,m$,
    its discrete weak divergence $\nabla_w\cdot\mathbf v_{i,h}$ is the unique polynomial in $P_{k-1}(T)$ such that}
	\begin{equation}
		\left(\nabla_w \cdot \mathbf{v}_{i,h}, \varphi\right)_T=-\left(\mathbf{v}_{i,0}, \nabla \varphi\right)_T+\left\langle \mathbf{v}_{i,b},\varphi\mathbf{n} \right\rangle_{\partial T}, \quad \forall\, \varphi \in P_{k-1}(T).\label{define-divergence}
	\end{equation}
\end{definition}
\begin{definition}
	For each $T\in\mathcal{T}_h^m$ and each $q_{m,h}=\{q_{m,0},q_{m,b}\}\in Z_h$, the \textcolor{black}{discrete} weak gradient $\widehat{\nabla}_w q_{m,h}\in [P_{k-1}(T)]^2$ satisfies
	\begin{equation}
		\big(\widehat{\nabla}_w q_{m,h},\bm{\psi}\big)_T=-(q_{m,0},\nabla\cdot\bm{\psi})_T+\left\langle q_{m,b}, {\bm{\psi}\cdot\bn}\right\rangle_{\partial T},\quad \forall\, \bm{\psi}\in \left[P_{k-1}(T)\right]^2.\label{define-gradient2}
	\end{equation}
\end{definition}
Next, we introduce some projection operators. For any $T\in\mathcal{T}_h$, denote by $Q_0,\,\widehat{Q}_h,\,\mathbf{Q}_0,\,\widehat{\mathbf{Q}}_h$ and $\mathbb{Q}_h$ the $L^2$ projection operators onto $P_k(T),\,P_{k-1}(T),\,[P_k(T)]^2,\,[P_{k-1}(T)]^{2}$ and $[P_{k-1}(T)]^{2\times 2}$, respectively. For each $e\in\mathcal{E}_h$, denote by $Q_b$ and $\mathbf{Q}_b$ the $L^2$ projection operators onto $P_k(e)$ and $[P_k(e)]^{2}$, respectively. We further define the composite operators $Q_h=\{Q_0,Q_b\}$ and $\mathbf{Q}_h=\{\mathbf{Q}_0,\mathbf{Q}_b\}$.
\par
Now, we proceed to define the bilinear forms. For $i=f,m$, let
\begin{align*}
	a_i(\mathbf{z}_{i,h},\mathbf{v}_{i,h})&:= 2\mu_i\sum_{T\in\mathcal{T}_h^i}\big(\mathbf{D}_w(\mathbf{z}_{i,h}),\mathbf{D}_w(\mathbf{v}_{i,h})\big)_T+s_i(\mathbf{z}_{i,h},\mathbf{v}_{i,h}),
	\\
	b_i(\mathbf{v}_{i,h},q_{i,h})&:= \sum_{T\in\mathcal{T}_h^i}\big(\nabla_w\cdot\mathbf{v}_{i,h},q_{i,h}\big)_T,
	\\
	a_p(y_{m,h},q_{m,h})&:= \kappa\sum_{T\in\mathcal{T}_h^m}\big(\widehat{\nabla}_w y_{m,h},\widehat{\nabla}_w q_{m,h}\big)_T+s_p(y_{m,h},q_{m,h}).
\end{align*}
Here, $a_i(\cdot,\cdot)$ is defined on $V_h^i\times V_h^i$, $b_i(\cdot,\cdot)$ on $V_h^i\times W_h^i$, and $a_p(\cdot,\cdot)$ on $Z_h\times Z_h$. The stabilization terms are given by
\begin{align*}
	s_i(\mathbf{z}_{i,h},\mathbf{v}_{i,h})&:=2\mu_i \sum_{T\in\mathcal{T}_h^i}h_T^{-1}\left\langle\mathbf{z}_{i,0}-\mathbf{z}_{i,b},\mathbf{v}_{i,0}-\mathbf{v}_{i,b}\right\rangle_{\partial T},
	\\
	s_p(y_{m,h},q_{m,h})&:=\kappa \sum_{T\in\mathcal{T}_h^m}h_T^{-1}\left\langle y_{m,0}-y_{m,b},q_{m,0}-q_{m,b}\right\rangle_{\partial T}.
\end{align*}
For $\mathbf z_{f,h},\mathbf v_{f,h}\in V_h^f,\ \mathbf z_{m,h},\mathbf v_{m,h}\in V_h^m$ and $ y_{m,h},q_{m,h}\in Z_h$, the interface terms are defined by
\begin{align*}
	a_{\Gamma}(\mathbf{z}_{f,h},\partial_t \mathbf{z}_{m,h};\mathbf{v}_{f,h},\mathbf{v}_{m,h})&:=\sum_{e\in\mathcal{E}_h^I}\beta_{\Gamma}\left\langle(\mathbf{z}_{f,b}-\partial_t \mathbf{z}_{m,b})\cdot\bm{\tau}_{\Gamma},(\mathbf{v}_{f,b}-\mathbf{v}_{m,b})\cdot\bm{\tau}_{\Gamma}\right\rangle_e,\\
	b_{\Gamma}(\mathbf{z}_{f,h},y_{m,h};\mathbf{v}_{f,h},q_{m,h})&:=\sum_{e\in\mathcal{E}_h^I}\left\langle\mathbf{v}_{f,b}\cdot\mathbf{n}_{\Gamma},y_{m,b}\right\rangle_e-\sum_{e\in\mathcal{E}_h^I}\left\langle\mathbf{z}_{f,b}\cdot\mathbf{n}_{\Gamma},q_{m,b}\right\rangle_e,\\
	c_{\Gamma}(\partial_t \mathbf{z}_{m,h},y_{m,h};\mathbf{v}_{m,h},q_{m,h})&:=\sum_{e\in\mathcal{E}_h^I}\left\langle\partial_t \mathbf{z}_{m,b}\cdot\mathbf{n}_{\Gamma},q_{m,b}\right\rangle_e-\sum_{e\in\mathcal{E}_h^I}\left\langle\mathbf{v}_{m,b}\cdot\mathbf{n}_{\Gamma},y_{m,b}\right\rangle_e.
\end{align*}

With the above preparations, the semi-discrete scheme of the coupled model \eqref{problem-eq1-s}--\eqref{problem-initial} is \textcolor{black}{described} as:
\textcolor{black}{
For any $t\in[0,T]$, seek
\[
\mathbf{u}_{f,h}(t)\in V_h^{f,D},\quad
p_{f,h}(t)\in W_h^f,\quad
\mathbf{u}_{m,h}(t)\in V_h^{m,D},\quad
p_{m,h}(t)\in Z_h^D,\quad
\phi_{m,h}(t)\in W_h^m,
\]
with the initial data
\[
\mathbf{u}_{m,h}(0)=\mathbf{Q}_h\mathbf{u}_m^0,\qquad
p_{m,h}(0)=Q_hp_m^0,\qquad
\phi_{m,h}(0)=\widehat{Q}_h\phi_m^0,
\]
such that the following equations hold:}
\begin{align}
	&a_f(\mathbf{u}_{f,h},\mathbf{v}_{f,h})-b_f(\mathbf{v}_{f,h},p_{f,h})+a_{\Gamma}(\mathbf{u}_{f,h},\partial_t \mathbf{u}_{m,h};\mathbf{v}_{f,h},\mathbf{0})+b_{\Gamma}(\mathbf{0},p_{m,h};\mathbf{v}_{f,h},0)=(\mathbf{f}_f,\mathbf{v}_{f,0})_{\Omega_f},\label{semi-scheme1}\\
	&b_f(\mathbf{u}_{f,h},q_{f,h})=0,\label{semi-scheme2}\\
	&a_m(\mathbf{u}_{m,h},\mathbf{v}_{m,h})-b_m(\mathbf{v}_{m,h},\phi_{m,h})+a_{\Gamma}(\mathbf{u}_{f,h},\partial_t \mathbf{u}_{m,h};\mathbf{0},\mathbf{v}_{m,h})\nonumber\\
	&\qquad+c_{\Gamma}(\mathbf{0},p_{m,h};\mathbf{v}_{m,h},0)=(\mathbf{f}_m,\mathbf{v}_{m,0})_{\Omega_m},\label{semi-scheme3}\\
	&\Bigl(c_0+\frac{\alpha^2}{\lambda}\Bigr)(\partial_t p_{m,0},q_{m,0})_{\Omega_m}-\frac{\alpha}{\lambda}(\partial_t \phi_{m,h},q_{m,0})_{\Omega_m}+a_p(p_{m,h},q_{m,h})+b_{\Gamma}(\mathbf{u}_{f,h},0;\mathbf{0},q_{m,h})\nonumber\\
	&\qquad+c_{\Gamma}(\partial_t \mathbf{u}_{m,h},0;\mathbf{0},q_{m,h})=(g,q_{m,0})_{\Omega_m},\label{semi-scheme4}\\
	&\frac{1}{\lambda}(\phi_{m,h},\psi_{m,h})_{\Omega_m}-\frac{\alpha}{\lambda}(p_{m,0},\psi_{m,h})_{\Omega_m}+b_m(\mathbf{u}_{m,h},\psi_{m,h})=0,\label{semi-scheme5}
\end{align}
for all $\mathbf{v}_{f,h}\in V_h^{f,D},\,q_{f,h}\in W_h^f,\,\mathbf{v}_{m,h}\in V_h^{m,D},\,q_{m,h}\in Z_h^D$ and $\psi_{m,h}\in W_h^m$.

\textcolor{black}{
As for the fully discrete formulation, we first introduce a uniform partition of $[0,T]$:
\[
0=t^0<t^1<\cdots<t^N=T,
\]
where the time step size is denoted by $\ell=\frac{T}{N}$ and $t^n=n\ell$ for $n=0,1,\ldots,N$.} For a time-dependent function $\omega$, we \textcolor{black}{
denote} $\omega^n=\omega(t^n)$ and define the backward \textcolor{black}{Euler} difference operator by
\begin{align*}
\bar{\partial}_t\omega^n=\frac{\omega^n-\omega^{n-1}}{\ell}.
\end{align*}
The fully discrete scheme of the coupled problem 
\eqref{problem-eq1-s}--\eqref{problem-initial} is presented as follows:
For each $n=1,2,\ldots,N$, seek
\[
\mathbf{u}_{f,h}^n\in V_h^{f,D},\quad
p_{f,h}^n\in W_h^f,\quad
\mathbf{u}_{m,h}^n\in V_h^{m,D},\quad
p_{m,h}^n\in Z_h^D,\quad
\phi_{m,h}^n\in W_h^m,
\]
with the initial approximations
\[
\mathbf{u}_{m,h}^0=\mathbf{Q}_h\mathbf{u}_m^0,\qquad
p_{m,h}^0=Q_hp_m^0,\qquad
\phi_{m,h}^0=\widehat{Q}_h\phi_m^0,
\]
such that 
\begin{align}
	&a_f(\mathbf{u}_{f,h}^n,\mathbf{v}_{f,h})-b_f(\mathbf{v}_{f,h},p_{f,h}^n)+a_{\Gamma}(\mathbf{u}_{f,h}^n,\bar{\partial}_t \mathbf{u}_{m,h}^n;\mathbf{v}_{f,h},\mathbf{0})+b_{\Gamma}(\mathbf{0},p_{m,h}^n;\mathbf{v}_{f,h},0)=(\mathbf{f}_f^n,\mathbf{v}_{f,0})_{\Omega_f},\label{full-scheme1}\\
	&b_f(\mathbf{u}_{f,h}^n,q_{f,h})=0,\label{full-scheme2}\\
	&a_m(\mathbf{u}_{m,h}^n,\mathbf{v}_{m,h})-b_m(\mathbf{v}_{m,h},\phi_{m,h}^n)+a_{\Gamma}(\mathbf{u}_{f,h}^n,\bar{\partial}_t \mathbf{u}_{m,h}^n;\mathbf{0},\mathbf{v}_{m,h})\nonumber\\
	&\qquad+c_{\Gamma}(\mathbf{0},p_{m,h}^n;\mathbf{v}_{m,h},0)=(\mathbf{f}_m^n,\mathbf{v}_{m,0})_{\Omega_m},\label{full-scheme3}\\
	&\Bigl(c_0+\frac{\alpha^2}{\lambda}\Bigr)(\bar{\partial}_t p_{m,0}^n,q_{m,0})_{\Omega_m}-\frac{\alpha}{\lambda}(\bar{\partial}_t \phi_{m,h}^n,q_{m,0})_{\Omega_m}+a_p(p_{m,h}^n,q_{m,h})+b_{\Gamma}(\mathbf{u}_{f,h}^n,0;\mathbf{0},q_{m,h})\nonumber\\
	&\qquad+c_{\Gamma}(\bar{\partial}_t \mathbf{u}_{m,h}^n,0;\mathbf{0},q_{m,h})=(g^n,q_{m,0})_{\Omega_m},\label{full-scheme4}\\
	&\frac{1}{\lambda}(\phi_{m,h}^n,\psi_{m,h})_{\Omega_m}-\frac{\alpha}{\lambda}(p_{m,0}^n,\psi_{m,h})_{\Omega_m}+b_m(\mathbf{u}_{m,h}^n,\psi_{m,h})=0,\label{full-scheme5}
\end{align}
for all $\mathbf{v}_{f,h}\in V_h^{f,D},\,q_{f,h}\in W_h^f,\,\mathbf{v}_{m,h}\in V_h^{m,D},\,q_{m,h}\in Z_h^D$ and $\psi_{m,h}\in W_h^m$. 

\section{Existence, uniqueness and stability}\label{section:exist-unique}
The aim of this section 
is to establish the well-posedness of the semi-discrete scheme \eqref{semi-scheme1}--\eqref{semi-scheme5} and the fully discrete scheme \eqref{full-scheme1}--\eqref{full-scheme5}. In the following content, unless otherwise specified, the general positive constant $C$ is independent of the mesh size h.

We first introduce two discrete semi-norms. For $i=f,m$ and $\mathbf{v}_{i,h}\in V_h^i$, define 
\begin{align}
    \trb{\mathbf{v}_{i,h}}_i^2:=\sum_{T\in\mathcal{T}_h^i}\bigl(\|\mathbf{D}_w(\mathbf{v}_{i,h})\|_T^2+h_T^{-1}\|\mathbf{v}_{i,0}-\mathbf{v}_{i,b}\|_{\partial T}^2\bigr).\label{norm-vector-1}
\end{align}
For $q_{m,h}\in Z_h$, define
\begin{align}
	\trb{q_{m,h}}_p^2:=\sum_{T\in\mathcal{T}_h^m}\bigl(\|\widehat{\nabla}_wq_{m,h}\|_T^2+h_T^{-1}\|q_{m,0}-q_{m,b}\|_{\partial T}^2\bigr).\label{norm-scale}
\end{align}
\begin{lemma}
	For $i=f,m$, $\trb{\cdot}_i$ provides a norm on the space $V_h^{i,D}$.
\end{lemma}
\begin{proof}
    Let $\mathbf{v}_{i,h}\in V_h^{i,D}$ and assume that $\trb{\mathbf{v}_{i,h}}_i=0$, we get
    \begin{align*}
    	&\sum_{T\in\mathcal{T}_h^i}(\|\mathbf{D}_w(\mathbf{v}_{i,h})\|_T^2+h_T^{-1}\|\mathbf{v}_{i,0}-\mathbf{v}_{i,b}\|_{\partial T}^2)=0,
    \end{align*}
    which implies that
    \begin{align}
    	\mathbf{D}_w(\mathbf{v}_{i,h})=\mathbf{0}~\text{in}~T\in\mathcal{T}_h^i,\ \mathbf{v}_{i,0}=\mathbf{v}_{i,b}~\text{on}~e\in\mathcal{E}_h^i.\label{norm-equality}
    \end{align}
    For any $\bm{\sigma}\in [P_{k-1}(T)]^{2\times 2}$, it follows from the definition of the weak gradient \eqref{define-gradient1} that
    \begin{align}
    	(\nabla_w\mathbf{v}_{i,h},\bm{\sigma})_T=-(\mathbf{v}_{i,0},\nabla\cdot\bm{\sigma})_T+\langle\mathbf{v}_{i,b},\bm{\sigma}\mathbf{n}\rangle_{\partial T}=(\nabla\mathbf{v}_{i,0},\bm{\sigma})_T.\label{norm-vector-proof1}
    \end{align}
    Similarly, we have
    \begin{align}
    	(\nabla_w^T\mathbf{v}_{i,h},\bm{\sigma})_T=(\nabla_w\mathbf{v}_{i,h},\bm{\sigma}^T)_T=(\nabla\mathbf{v}_{i,0},\bm{\sigma}^T)_T=(\nabla^T\mathbf{v}_{i,0},\bm{\sigma})_T.\label{norm-vector-proof2}
    \end{align}
    Adding \eqref{norm-vector-proof1} and \eqref{norm-vector-proof2} to get
    \begin{align*}
    	(\mathbf{D}_w(\mathbf{v}_{i,h}),\bm{\sigma})_T=(\mathbf{D}(\mathbf{v}_{i,0}),\bm{\sigma})_T.
    \end{align*}
    Taking $\bm{\sigma}=\mathbf{D}(\mathbf{v}_{i,0})$ leads to $\mathbf{D}(\mathbf{v}_{i,0})=0$ in $T\in\mathcal{T}_h^i$. Thus, $\mathbf{v}_{i,0}\in \text{RM}$ in $T$, where
    \[\mathrm{RM}=\mathrm{span}\left\{
    \begin{pmatrix}1\\0\end{pmatrix},
    \begin{pmatrix}0\\1\end{pmatrix},
    \begin{pmatrix}-y\\x\end{pmatrix}\right\},\]
    represents \textcolor{black}{
    the rigid body motion space}. Combining with \eqref{norm-equality}, we know that $\mathbf{v}_{i,0}$ is continuous in $\Omega_i$. 
    Moreover, since $\mathbf{v}_{i,b}|_{\Gamma_i^D}=\mathbf{0}$, we have $\mathbf{v}_{i,0}|_{\Gamma_i^D}=\mathbf{0}$. Applying Korn's inequality \cite{MR3508859-korn-inequality}, we obtain
    \begin{align*}
    	\|\mathbf{v}_{i,0}\|_1\leq C(\|\mathbf{D}(\mathbf{v}_{i,0})\|_{\Omega_i}+\|\mathbf{v}_{i,0}\|_{\Gamma_i^D}),
    \end{align*}
    which implies $\mathbf{v}_{i,0}\equiv \mathbf{0}$ in $\Omega_i$, and $\mathbf{v}_{i,b}\equiv \mathbf{0}$ on $e\in \mathcal{E}_h^i$. The proof is complete.
\end{proof}
	Similarly, $\trb{\cdot}_p$ defines a norm on the space $Z_h^D$ \cite{MR3452926-WG-stokes-2016}.

\begin{lemma}[\cite{MR3325251-2015-WG-ellipse,MR3286455-2015-v0-vh}]
	\textcolor{black}{
    There exist constants $C_1,C_2>0$ such that}
	\begin{align}
		\sum_{T\in\mathcal{T}_h^m}\|\nabla q_{m,0}\|_T^2&\leq C_1\trb{q_{m,h}}_p^2,\quad\forall\,q_{m,h}\in Z_h,\label{nabla-weak-nabla1}\\
		\sum_{T\in\mathcal{T}_h^m}\|q_{m,0}\|_T^2&\leq C_2\trb{q_{m,h}}_p^2,\quad\forall\,q_{m,h}\in Z_h^D.\label{v0-vh1}
	\end{align}
\end{lemma}
\begin{lemma}
	For each $q_{m,h}\in Z_h^D$, \textcolor{black}{there holds}
	\begin{align}
		\sum_{e\in\mathcal{E}_h^I}\|q_{m,b}\|_e^2\leq C\trb{q_{m,h}}_p^2.\label{vb-vh1}
	\end{align}
\end{lemma}
\begin{proof}
	Using the triangle inequality, we have
	\begin{align}
		\sum_{e\in\mathcal{E}_h^I}\|q_{m,b}\|_e^2\leq \sum_{e\in\mathcal{E}_h^I}\|q_{m,b}-q_{m,0}\|_e^2+\sum_{e\in\mathcal{E}_h^I}\|q_{m,0}\|_e^2.\label{vb-vh-proof1}
	\end{align}
    It follows from \eqref{norm-scale} that 
    \begin{align*}
    	\sum_{e\in\mathcal{E}_h^I}\|q_{m,b}-q_{m,0}\|_e^2\leq h\trb{q_{m,h}}_p^2\leq \trb{q_{m,h}}_p^2.
    \end{align*}
    According to the broken trace theorem \cite{MR2557047-2009-broken-trace}, we have
    \begin{align*}
    	\sum_{e\in\mathcal{E}_h^I}\|q_{m,0}\|_e^2\leq C\left(\|q_{m,0}\|^2_{L^1(\Omega_m)}+\|\nabla q_{m,0}\|^2_{\Omega_m}+\sum_{T\in\mathcal{T}_h^m}\int_{\partial T}h_T^{-1}|\llbracket q_{m,0} \rrbracket|^2\mathrm{d}s\right),
    \end{align*}
    where $\llbracket q_{m,0} \rrbracket$ denotes the jump across element interfaces. The estimate \eqref{v0-vh1} gives
    \begin{align*}
    	\|q_{m,0}\|_{L^1(\Omega_m)}=\int_{\Omega_m}|q_{m,0}|\,\mathrm{d}x\leq\left(\int_{\Omega_m}|q_{m,0}|^2\,\mathrm{d}x\right)^{\frac{1}{2}}\left(\int_{\Omega_m}1^2\,\mathrm{d}x\right)^{\frac{1}{2}}=|\Omega_m|^{\frac{1}{2}}\|q_{m,0}\|_{\Omega_m}\leq C\trb{q_{m,h}}_p.
    \end{align*}
    From \eqref{nabla-weak-nabla1} and \eqref{norm-scale}, we obtain
    \begin{align*}
    	&\|\nabla q_{m,0}\|^2_{\Omega_m}+\sum_{T\in\mathcal{T}_h^m}\int_{\partial T}h_T^{-1}|\llbracket q_{m,0} \rrbracket|^2\mathrm{d}s\\
    	\leq\,& C\trb{q_{m,h}}_p^2+\sum_{T\in\mathcal{T}_h^m}\int_{\partial T}h_T^{-1}|\llbracket q_{m,0}-q_{m,b} \rrbracket|^2\mathrm{d}s\\
    	\leq\,& C\trb{q_{m,h}}_p^2+2\sum_{T\in\mathcal{T}_h^m}h_T^{-1}\| q_{m,0}-q_{m,b} \|^2_{\partial T}\\
    	\leq\,& C\trb{q_{m,h}}_p^2,
    \end{align*}
    The proof is complete.
\end{proof}
\begin{lemma}[\cite{MR3286455-2015-v0-vh,MR3508859-korn-inequality}]
	For $i=f,m$, \textcolor{black} {the following estimates hold:}
	\begin{align}
		\sum_{T\in\mathcal{T}_h^i}\|\mathbf{D}(\mathbf{v}_{i,0})\|_T^2&\leq C_3\trb{\mathbf{v}_{i,h}}_i^2,\quad\forall \,\mathbf{v}_{i,h}\in V_h^i,\label{nabla-weak-nabla2}\\
		\sum_{T\in\mathcal{T}_h^i}\|\mathbf{v}_{i,0}\|_T^2&\leq C_4\trb{\mathbf{v}_{i,h}}_i^2,\quad\forall\, \mathbf{v}_{i,h}\in V_h^{i,D}.\label{v0-vh2}
	\end{align}
\end{lemma}
\begin{lemma}
	For $i=f,m$ and each $\mathbf{v}_{i,h}\in V_h^{i,D}$, the interface component $\mathbf{v}_{i,b}$ satisfies
	\begin{align}
		\sum_{e\in\mathcal{E}_h^I}\|\mathbf{v}_{i,b}\|_e^2\leq C\trb{\mathbf{v}_{i,h}}_i^2.\label{vb-vh2}
	\end{align}
\end{lemma}
\begin{lemma}[\cite{MR3452926-WG-stokes-2016}]
	\textcolor{black}{ 
    The projection operators introduced earlier possess the following commutativity properties:}
	\begin{align}
		\nabla_w(\mathbf{Q}_h \mathbf{z}_i)&=\mathbb{Q}_h(\nabla\mathbf{z}_i),\quad\forall\,\mathbf{z}_i\in [H^1(\Omega_i)]^2,\quad i=f,m,\label{weak-gradient-property-vector}\\
		\nabla_w\cdot(\mathbf{Q}_h \mathbf{z}_i)&=\widehat{Q}_h(\nabla\cdot\mathbf{z}_i),\quad\forall\,\mathbf{z}_i\in H(div,\Omega_i),\quad i=f,m,\label{weak-divergence-property-vector}\\
		\widehat{\nabla}_w(Q_h y_m)&=\widehat{\mathbf{Q}}_h(\nabla y_m),\quad\forall\,y_m\in H^1(\Omega_m).\label{nabla-w-property-scale}
	\end{align}
\end{lemma}
According to the properties of the transpose operator together with \eqref{weak-gradient-property-vector}, it is easy to get
\begin{align}
	\mathbf{D}_w(\mathbf{Q}_h \mathbf{z}_i)=\mathbb{Q}_h\mathbf{D}(\mathbf{z}_i),\quad\forall\,\mathbf{z}_i\in [H^1(\Omega_i)]^2,\quad i=f,m.\label{weak-gradient-property-tensor}
\end{align}
\begin{lemma}[\cite{MR3452926-WG-stokes-2016,MR4799262-2025-WG-Biot}]
	\textcolor{black}{
    There are positive constants $\beta_1>0$ and $\beta_2>0$ such that}
	\begin{align}
		\sup_{\mathbf{v}_{f,h}\in V_h^{f,D}} \frac{b_f(\mathbf{v}_{f,h},q_{f,h})}{\trb{\mathbf{v}_{f,h}}_f}&\geq \beta_1\| q_{f,h}\|_{\Omega_f},\quad \forall \,q_{f,h}\in W_h^f,\label{inf-sup-Stokes}\\
		\sup_{\mathbf{v}_{m,h}\in V_h^{m,D}} \frac{b_m(\mathbf{v}_{m,h},\psi_{m,h})}{\trb{\mathbf{v}_{m,h}}_m}& \geq \beta_2\| \psi_{m,h}\|_{\Omega_m},\quad \forall \,\psi_{m,h}\in W_h^m.\label{inf-sup-Biot}
	\end{align}
\end{lemma}
\begin{theorem}
	The semi-discrete scheme \eqref{semi-scheme1}--\eqref{semi-scheme5} has a unique solution for any $t\in[0,T]$.
\end{theorem}
\begin{proof}
	\textcolor{black}{Let} $(\mathbf{u}_{f,h},\,p_{f,h},\,\mathbf{u}_{m,h},\,p_{m,h},\,\phi_{m,h})$ be the solution of the following homogeneous problem with $\mathbf{u}_{m,h}(0)=\mathbf{0},\,p_{m,h}(0)=0,\,\phi_{m,h}(0)=0$:
	\begin{align}
		&a_f(\mathbf{u}_{f,h},\mathbf{v}_{f,h})-b_f(\mathbf{v}_{f,h},p_{f,h})+a_{\Gamma}(\mathbf{u}_{f,h},\partial_t \mathbf{u}_{m,h};\mathbf{v}_{f,h},\mathbf{0})+b_{\Gamma}(\mathbf{0},p_{m,h};\mathbf{v}_{f,h},0)=0,\label{semi-unique1}\\
		&b_f(\mathbf{u}_{f,h},q_{f,h})=0,\label{semi-unique2}\\
		&a_m(\mathbf{u}_{m,h},\mathbf{v}_{m,h})-b_m(\mathbf{v}_{m,h},\phi_{m,h})+a_{\Gamma}(\mathbf{u}_{f,h},\partial_t \mathbf{u}_{m,h};\mathbf{0},\mathbf{v}_{m,h})+c_{\Gamma}(\mathbf{0},p_{m,h};\mathbf{v}_{m,h},0)=0,\label{semi-unique3}\\
		&\Bigl(c_0+\frac{\alpha^2}{\lambda}\Bigr)\bigl(\partial_t p_{m,0},q_{m,0}\bigr)_{\Omega_m}-\frac{\alpha}{\lambda}\bigl(\partial_t \phi_{m,h},q_{m,0}\bigr)_{\Omega_m}+a_p(p_{m,h},q_{m,h})+b_{\Gamma}(\mathbf{u}_{f,h},0;\mathbf{0},q_{m,h})\nonumber\\
		&\qquad+c_{\Gamma}(\partial_t \mathbf{u}_{m,h},0;\mathbf{0},q_{m,h})=0,\label{semi-unique4}\\
		&\frac{1}{\lambda}\bigl(\partial_t \phi_{m,h},\psi_{m,h}\bigr)_{\Omega_m}-\frac{\alpha}{\lambda}\bigl(\partial_t p_{m,0},\psi_{m,h}\bigr)_{\Omega_m}+b_m(\partial_t \mathbf{u}_{m,h},\psi_{m,h})=0,\label{semi-unique5}
	\end{align}
    for all $\mathbf{v}_{f,h}\in V_h^{f,D},\,q_{f,h}\in W_h^f,\,\mathbf{v}_{m,h}\in V_h^{m,D},\,q_{m,h}\in Z_h^D$ and $\psi_{m,h}\in W_h^m$. \\
    Then, in view of Theorem 2.4 in \cite{MR1946758-trivial-solution}, it suffices to prove that $(\mathbf{u}_{f,h},\,p_{f,h},\,\mathbf{u}_{m,h},\,p_{m,h},\,\phi_{m,h})\equiv\mathbf{0}.$ Setting $\mathbf{v}_{f,h}=\mathbf{u}_{f,h},\,q_{f,h}=p_{f,h},\,\mathbf{v}_{m,h}=\partial_t \mathbf{u}_{m,h},\,q_{m,h}=p_{m,h}$ and $\psi_{m,h}=\phi_{m,h}$ in \eqref{semi-unique1}--\eqref{semi-unique5}, and summing them together yields
    \begin{align}
    	&\dfrac{1}{2}\dfrac{\mathrm{d}}{\mathrm{d}t}\biggl(a_m(\mathbf{u}_{m,h},\mathbf{u}_{m,h})+c_0\|p_{m,0}\|^2_{\Omega_m}+\dfrac{1}{\lambda}\|\alpha p_{m,0}-\phi_{m,h}\|^2_{\Omega_m}\biggr)+a_f(\mathbf{u}_{f,h},\mathbf{u}_{f,h})+a_p(p_{m,h},p_{m,h})\nonumber\\
    	&\qquad+\beta_{\Gamma}\|(\mathbf{u}_{f,b}-\partial_t \mathbf{u}_{m,b})\cdot\bm{\tau}_{\Gamma}\|_{\Gamma}^2=0.\label{semi-unique6}
    \end{align}
    \textcolor{black}{Integrating the above equality over the interval $(0,t)$ and using the initial conditions gives}
    \begin{align*}
    	&2\mu_m\trb{\mathbf{u}_{m,h}(t)}_m^2+c_0\|p_{m,0}(t)\|^2_{\Omega_m}+\dfrac{1}{\lambda}\|\alpha p_{m,0}(t)-\phi_{m,h}(t)\|^2_{\Omega_m}\\
        &\quad+2\int_0^t\biggl(2\mu_f\trb{\mathbf{u}_{f,h}(r)}_f^2+\kappa\trb{p_{m,h}(r)}_p^2
    	+\sum_{e\in\mathcal{E}_h^I}\beta_{\Gamma}\|(\mathbf{u}_{f,b}(r)-\partial_t \mathbf{u}_{m,b}(r))\cdot\bm{\tau}_{\Gamma}\|_e^2\biggr)\,\mathrm{d}r=0.
    \end{align*}
    \textcolor{black}{The above identity implies} $\mathbf{u}_{m,h}(t)=\mathbf{0},\,\mathbf{u}_{f,h}(t)=\mathbf{0},\,p_{m,h}(t)=0$ and $\phi_{m,h}(t)=0$. Substituting these results, we obtain
    \begin{align*}
    	b_f(\mathbf{v}_{f,h},p_{f,h})=0.
    \end{align*}
    Combining the estimate \eqref{inf-sup-Stokes}, we conclude that $\|p_{f,h}\|^2=0$, and consequently, $p_{f,h}=0$. This completes the proof.
\end{proof}
\begin{theorem}
	Assume that $\mathbf f_f\in L^2(0,T;[L^2(\Omega_f)]^2),\,g\in L^2(0,T;L^2(\Omega_m)),\,\mathbf f_m\in H^1(0,T;[L^2(\Omega_m)]^2)$. Let $(\bu_{f,h}(t),\,p_{f,h}(t),\,\mathbf{u}_{m,h}(t),\,p_{m,h}(t),\,\phi_{m,h}(t))$ be the solution of the semi-discrete scheme \eqref{semi-scheme1}--\eqref{semi-scheme5}. Then, for any $0<t\leq T$, the following estimate holds
	\begin{align*}
		&2\mu_m\trb{\mathbf{u}_{m,h}(t)}_m^2+c_0\|p_{m,0}(t)\|^2_{\Omega_m}+\dfrac{1}{\lambda}\|\alpha p_{m,0}(t)-\phi_{m,h}(t)\|^2_{\Omega_m}+\int_0^t(2\mu_f\trb{\mathbf{u}_{f,h}(r)}_f^2+\kappa\trb{p_{m,h}(r)}_p^2\\
		&+\beta_{\Gamma}\|(\bu_{f,b}(r)-\partial_t \mathbf{u}_{m,b}(r))\cdot\bm{\tau}_{\Gamma}\|_{\Gamma}^2+\|p_{f,h}(r)\|^2_{\Omega_f}+\|\phi_{m,h}(r)\|^2_{\Omega_m})\,\mathrm{d}r\\
		\leq\ & C \biggl(2\mu_m\trb{\mathbf{u}_{m,h}(0)}_m^2+c_0\|p_{m,0}(0)\|^2_{\Omega_m}+\dfrac{1}{\lambda}\|\alpha p_{m,0}(0)-\phi_{m,h}(0)\|^2_{\Omega_m}+\int_0^t\bigl(\|\mathbf{f}_f(r)\|_{\Omega_f}^2+\|\partial_t \mathbf{f}_m(r)\|_{\Omega_m}^2\\
		&+\|g(r)\|_{\Omega_m}^2\bigr)\,\mathrm{d}r+\max\limits_{0\leq r\leq t}\|\mathbf{f}_m(r)\|_{\Omega_m}^2\biggr).
	\end{align*}
\end{theorem}
\begin{proof}
	Taking $\mathbf{v}_{f,h}=\mathbf{u}_{f,h},\,q_{f,h}=p_{f,h},\,\mathbf{v}_{m,h}=\partial_t \mathbf{u}_{m,h},\,q_{m,h}=p_{m,h} $ and $\psi_{m,h}=\phi_{m,h}$ in \eqref{semi-scheme1}--\eqref{semi-scheme5}, and summing them together, we have
	\begin{align}
		&\dfrac{1}{2}\dfrac{\mathrm{d}}{\mathrm{d}t}(a_m(\mathbf{u}_{m,h},\mathbf{u}_{m,h})+c_0\|p_{m,0}\|^2_{\Omega_m}+\dfrac{1}{\lambda}\|\alpha p_{m,0}-\phi_{m,h}\|^2_{\Omega_m})+a_f(\bu_{f,h},\bu_{f,h})+a_p(p_{m,h},p_{m,h})\nonumber\\
		&\qquad+\beta_{\Gamma}\|(\bu_{f,b}-\partial_t \mathbf{u}_{m,b})\cdot\bm{\tau}_{\Gamma}\|_{\Gamma}^2=(\mathbf{f}_f,\bu_{f,0})_{\Omega_f}+(\mathbf{f}_m,\partial_t \mathbf{u}_{m,0})_{\Omega_m}+(g,p_{m,0})_{\Omega_m}.\label{stability-proof1}
	\end{align}
    For brevity, we introduce the following notations
    \begin{align*}
    	&A(t)^2:=2\mu_m\trb{\mathbf{u}_{m,h}(t)}_m^2+c_0\|p_{m,0}(t)\|^2_{\Omega_m}+\dfrac{1}{\lambda}\|\alpha p_{m,0}(t)-\phi_{m,h}(t)\|^2_{\Omega_m},\\
    	&B(t)^2:=2\mu_f\trb{\mathbf{u}_{f,h}(t)}_f^2+\kappa\trb{p_{m,h}(t)}_p^2+\beta_{\Gamma}\|(\bu_{f,b}(t)-\partial_t \mathbf{u}_{m,b}(t))\cdot\bm{\tau}_{\Gamma}\|_{\Gamma}^2.
    \end{align*}
    Integrating \eqref{stability-proof1} \textcolor{black}{over the interval $(0,t)$} gives
    \begin{align}
    	\dfrac{1}{2}(A(t)^2-A(0)^2)+\int_0^t B(r)^2 \,\mathrm{d}r=&\int_0^t\bigl((\mathbf{f}_f(r),\bu_{f,0}(r))_{\Omega_f}+(\mathbf{f}_m(r),\partial_t \mathbf{u}_{m,0}(r))_{\Omega_m}\bigr)\,\mathrm{d}r\nonumber\\
    	&+\int_0^t(g(r),p_{m,0}(r))_{\Omega_m}\,\mathrm{d}r.\label{stability-proof2}
    \end{align}
    \textcolor{black}{
    The first and third terms on the right-hand side of \eqref{stability-proof2} are estimated with the help of \eqref{v0-vh1} and \eqref{v0-vh2}:}
    \begin{align*}
    	\int_0^t(\mathbf{f}_f(r),\bu_{f,0}(r))_{\Omega_f}\,\mathrm{d}r&\leq\int_0^t\|\mathbf{f}_f(r)\|_{\Omega_f}\|\bu_{f,0}(r)\|_{\Omega_f}\,\mathrm{d}r\leq\int_0^t\|\mathbf{f}_f(r)\|_{\Omega_f}\trb{\bu_{f,h}(r)}_f\,\mathrm{d}r,\\
    	\int_0^t(g(r),p_{m,0}(r))_{\Omega_m}\,\mathrm{d}r&\leq\int_0^t\|g(r)\|_{\Omega_m}\|p_{m,0}(r)\|_{\Omega_m}\,\mathrm{d}r\leq\int_0^t\|g(r)\|_{\Omega_m}\trb{p_{m,h}(r)}_p\,\mathrm{d}r.
    \end{align*}
    For the second source term, using integration by parts in time and the same estimates above,
    \begin{align*}
    	&\int_0^t(\mathbf{f}_m(r),\partial_t \mathbf{u}_{m,0}(r))_{\Omega_m}\,\mathrm{d}r\\
    	=\ &(\mathbf{f}_m(t),\mathbf{u}_{m,0}(t))_{\Omega_m}-(\mathbf{f}_m(0),\mathbf{u}_{m,0}(0))_{\Omega_m}-\int_0^t(\partial_t \mathbf{f}_m(r),\mathbf{u}_{m,0}(r))_{\Omega_m}\,\mathrm{d}r\\
    	\leq&\ \|\mathbf{f}_m(t)\|_{\Omega_m}\trb{\mathbf{u}_{m,h}(t)}_m+\|\mathbf{f}_m(0)\|_{\Omega_m}\trb{\mathbf{u}_{m,h}(0)}_m+\int_0^t\|\partial_t \mathbf{f}_m(r)\|_{\Omega_m}\trb{\mathbf{u}_{m,h}(r)}_m\,\mathrm{d}r.
    \end{align*}
    Then, we have
    \begin{align}
        \dfrac{1}{2}A(t)^2+\int_0^t B(r)^2 \,\mathrm{d}r\leq\ &\dfrac{1}{2}A(0)^2+C\int_0^t\bigl(\|\mathbf{f}_f(r)\|_{\Omega_f}B(r)+\|\partial_t \mathbf{f}_m(r)\|_{\Omega_m}A(r)+\|g(r)\|_{\Omega_m}B(r)\bigr)\,\mathrm{d}r\nonumber\\
        &+C\|\mathbf{f}_m(t)\|_{\Omega_m}A(t)+C\|\mathbf{f}_m(0)\|_{\Omega_m}A(0).\label{stability-proof3}
    \end{align}
    It follows from Young's inequality that
    \begin{align*}
    	\frac{1}{2}A(t)^2+\int_0^t B(r)^2 \,\mathrm{d}r\leq\ & \dfrac{1}{2}A(0)^2+C\int_0^t\bigl(\|\mathbf{f}_f(r)\|_{\Omega_f}^2+\|\partial_t \mathbf{f}_m(r)\|_{\Omega_m}^2+\|g(r)\|_{\Omega_m}^2\bigr)\,\mathrm{d}r+\int_0^tB(r)^2\,\mathrm{d}r\\
    	&+\frac{1}{4}\int_0^t A(r)^2\,\mathrm{d}r+C\bigl(\|\mathbf{f}_m(t)\|_{\Omega_m}^2+\|\mathbf{f}_m(0)\|_{\Omega_m}^2\bigr)+\frac{1}{4}A(t)^2+\frac{1}{4}A(0)^2,
    \end{align*}
    which implies
    \begin{align*}
    	\frac{1}{4}A(t)^2\leq&\  \dfrac{3}{4}A(0)^2+C\int_0^t\bigl(\|\mathbf{f}_f(r)\|_{\Omega_f}^2+\|\partial_t \mathbf{f}_m(r)\|_{\Omega_m}^2+\|g(r)\|_{\Omega_m}^2\bigr)\,\mathrm{d}r+C\bigl(\|\mathbf{f}_m(t)\|_{\Omega_m}^2+\|\mathbf{f}_m(0)\|_{\Omega_m}^2\bigr)\nonumber\\
    	&+\frac{1}{4}\int_0^tA(r)^2\,\mathrm{d}r.
    \end{align*}
    Since $0<t\le T$, the Gronwall inequality \cite{MR4799262-2025-WG-Biot} gives
    \begin{align}
    	A(t)^2\leq C\bigl(A(0)^2+K(t)\bigr),\label{stability-proof5}
    \end{align}
    where $K(t)=\int_0^t\bigl(\|\mathbf{f}_f(r)\|_{\Omega_f}^2+\|\partial_t \mathbf{f}_m(r)\|_{\Omega_m}^2+\|g(r)\|_{\Omega_m}^2\bigr)\,\mathrm{d}r+\max\limits_{0\leq r\leq t}\|\mathbf{f}_m(r)\|^2_{\Omega_m}$.
    
    Employing Young's inequality again with different coefficients in \eqref{stability-proof3}, we have
    \begin{align*}
    	\int_0^tB(r)^2\,\mathrm{d}r\leq&\  2A(0)^2+C\int_0^t\bigl(\|\mathbf{f}_f(r)\|_{\Omega_f}^2+\|\partial_t \mathbf{f}_m(r)\|_{\Omega_m}^2+\|g(r)\|_{\Omega_m}^2\bigr)\,\mathrm{d}r+\int_0^tA(r)^2\,\mathrm{d}r\nonumber\\
    	&+C\bigl(\|\mathbf{f}_m(t)\|_{\Omega_m}^2+\|\mathbf{f}_m(0)\|_{\Omega_m}^2\bigr).
    \end{align*}
    It follows from \eqref{stability-proof5} that
    \begin{align}
    	\int_0^tA(r)^2\,\mathrm{d}r\leq C\int_0^t\bigl(A(0)^2+K(r)\bigr)\,\mathrm{d}r\leq C\int_0^t\bigl(A(0)^2+K(t)\bigr)\,\mathrm{d}r\leq C \bigl(A(0)^2+K(t)\bigr).\label{stability-proof6}
    \end{align}
    Furthermore, we get
    \begin{align}
    	\int_0^tB(r)^2\,\mathrm{d}r\leq  C\bigl(A(0)^2+K(t)\bigr).\label{stability-proof7}
    \end{align}
   
    According to \eqref{semi-scheme1}, \eqref{vb-vh2} and \eqref{v0-vh2}, the following estimate holds
    \begin{align*}
    	&b_f(\mathbf{v}_{f,h},p_{f,h})\\
    	=\ &a_f(\mathbf{u}_{f,h},\mathbf{v}_{f,h})+a_{\Gamma}(\mathbf{u}_{f,h},\partial_t \mathbf{u}_{m,h};\mathbf{v}_{f,h},\mathbf{0})+b_{\Gamma}(\mathbf{0},p_{m,h};\mathbf{v}_{f,h},0)-(\mathbf{f}_f,\mathbf{v}_{f,0})_{\Omega_f}\\
    	\leq\ & C\bigl(\trb{\mathbf{u}_{f,h}}_f\trb{\mathbf{v}_{f,h}}_f+\|(\bu_{f,b}-\partial_t \mathbf{u}_{m,b})\cdot\bm{\tau}_{\Gamma}\|_{\Gamma}\|\mathbf{v}_{f,b}\|_{\Gamma}+\|p_{m,b}\|_{\Gamma}\|\mathbf{v}_{f,b}\|_{\Gamma}+\|\mathbf{f}_f\|_{\Omega_f}\|\mathbf{v}_{f,0}\|_{\Omega_f}\bigr)\\
    	\leq\ & C\bigl(\trb{\mathbf{u}_{f,h}}_f\trb{\mathbf{v}_{f,h}}_f+\|(\bu_{f,b}-\partial_t \mathbf{u}_{m,b})\cdot\bm{\tau}_{\Gamma}\|_{\Gamma}\trb{\mathbf{v}_{f,h}}_f+\trb{p_{m,h}}_p\trb{\mathbf{v}_{f,h}}_f+\|\mathbf{f}_f\|_{\Omega_f}\trb{\mathbf{v}_{f,h}}_f\bigr).
    \end{align*}
    \textcolor{black}{Invoking} the inf-sup condition \eqref{inf-sup-Stokes} \textcolor{black}{ gives}
    \begin{align}
    	\| p_{f,h}\|_{\Omega_f}^2\leq C\bigl(\trb{\mathbf{u}_{f,h}}_f^2+\|(\bu_{f,b}-\partial_t \mathbf{u}_{m,b})\cdot\bm{\tau}_{\Gamma}\|_{\Gamma}^2+\trb{p_{m,h}}_p^2+\|\mathbf{f}_f\|_{\Omega_f}^2\bigr).\label{stability-proof8}
    \end{align}
    \textcolor{black}{
    Integrating the above estimate over the interval $(0,t)$ and applying
\eqref{stability-proof7}, we know}
    \begin{align*}
    	\int_0^t\|p_{f,h}\|^2_{\Omega_f}\,\mathrm{d}r \leq C\int_0^t\bigl(B(r)^2+\|\mathbf{f}_f\|^2_{\Omega_f}\bigr)\,\mathrm{d}r\leq C\bigl(A(0)^2+K(t)\bigr).
    \end{align*}
   
    \textcolor{black}{
    Following the same procedure and utilizing the inf-sup condition \eqref{inf-sup-Biot}, estimates \eqref{stability-proof6} and \eqref{stability-proof7}, we further obtain}
    \begin{align*}
    	\int_0^t\|\phi_{m,h}\|^2_{\Omega_m}\,\mathrm{d}r &\leq C\int_0^t\bigl(\trb{\mathbf{u}_{m,h}}_m^2+\|(\bu_{f,b}-\partial_t \mathbf{u}_{m,b})\cdot\bm{\tau}_{\Gamma}\|_{\Gamma}^2+\trb{p_{m,h}}_p^2+\|\mathbf{f}_m\|_{\Omega_m}^2\bigr)\,\mathrm{d}r\\
    	&\leq C\int_0^t\bigl(A(r)^2+B(r)^2+\|\mathbf{f}_m\|^2_{\Omega_m}\bigr)\,\mathrm{d}r\\
    	&\leq C\bigl(A(0)^2+K(t)\bigr).
    \end{align*}
    This completes the proof.
\end{proof}

\begin{theorem}\label{existence-uniqueness}
	For each time level $t^{\bar{n}},\, 1\le \bar{n}\le N$, the fully discrete scheme \eqref{full-scheme1}--\eqref{full-scheme5} has a unique solution.
\end{theorem}
\begin{proof}
    \textcolor{black}{To prove uniqueness, we consider the corresponding homogeneous problem.
	Suppose that} $(\mathbf{u}_{f,h}^n,\,p_{f,h}^n,\,\mathbf{u}_{m,h}^n,\,p_{m,h}^n,\,\phi_{m,h}^n)$ where $n=1,\dots,\bar{n}$, satisfy the following homogeneous system:
	\begin{align}
		&a_f(\mathbf{u}_{f,h}^n,\mathbf{v}_{f,h})-b_f(\mathbf{v}_{f,h},p_{f,h}^n)+a_{\Gamma}(\mathbf{u}_{f,h}^n,\bar{\partial}_t \mathbf{u}_{m,h}^n;\mathbf{v}_{f,h},\mathbf{0})+b_{\Gamma}(\mathbf{0},p_{m,h}^n;\mathbf{v}_{f,h},0)=0,\label{full-unique1}\\
		&b_f(\mathbf{u}_{f,h}^n,q_{f,h})=0,\label{full-unique2}\\
		&a_m(\mathbf{u}_{m,h}^n,\mathbf{v}_{m,h})-b_m(\mathbf{v}_{m,h},\phi_{m,h}^n)+a_{\Gamma}(\mathbf{u}_{f,h}^n,\bar{\partial}_t \mathbf{u}_{m,h}^n;\mathbf{0},\mathbf{v}_{m,h})+c_{\Gamma}(\mathbf{0},p_{m,h}^n;\mathbf{v}_{m,h},0)=0,\label{full-unique3}\\
		&\Bigl(c_0+\frac{\alpha^2}{\lambda}\Bigr)(\bar{\partial}_t p_{m,0}^n,q_{m,0})_{\Omega_m}-\frac{\alpha}{\lambda}(\bar{\partial}_t \phi_{m,h}^n,q_{m,0})_{\Omega_m}+a_p(p_{m,h}^n,q_{m,h})+b_{\Gamma}(\mathbf{u}_{f,h}^n,0;\mathbf{0},q_{m,h})\nonumber\\
		&\qquad+c_{\Gamma}(\bar{\partial}_t \mathbf{u}_{m,h}^n,0;\mathbf{0},q_{m,h})=0,\label{full-unique4}\\
		&\frac{1}{\lambda}(\bar{\partial}_t \phi_{m,h}^n,\psi_{m,h})_{\Omega_m}-\frac{\alpha}{\lambda}(\bar{\partial}_t p_{m,0}^n,\psi_{m,h})_{\Omega_m}+b_m(\bar{\partial}_t \mathbf{u}_{m,h}^n,\psi_{m,h})=0,\label{full-unique5}
	\end{align}
	for all $\mathbf{v}_{f,h}\in V_h^{f,D},\,q_{f,h}\in W_h^f,\,\mathbf{v}_{m,h}\in V_h^{m,D},\,q_{m,h}\in Z_h^D$ and $\psi_{m,h}\in W_h^m$, with $\mathbf{u}_{m,h}^0=\mathbf{0},\,p_{m,h}^0=0,\,\phi_{m,h}^0=0$. \\
    Taking $\mathbf{v}_{f,h}=\mathbf{u}_{f,h}^n,\,q_{f,h}=p_{f,h}^n,\,\mathbf{v}_{m,h}=\bar{\partial}_t \mathbf{u}_{m,h}^n,\,q_{m,h}=p_{m,h}^n$ and $\psi_{m,h}=\phi_{m,h}^n$ in \eqref{full-unique1}--\eqref{full-unique5}, and summing the resulting \textcolor{black}{identities}, we obtain 
	 \begin{align*}
	 	&2\mu_f\trb{\mathbf{u}_{f,h}^n}_f^2+\kappa\trb{p_{m,h}^n}_p^2+\beta_{\Gamma}\|(\mathbf{u}_{f,b}^n-\bar{\partial}_t \mathbf{u}_{m,b}^n)\cdot\bm{\tau}_{\Gamma}\|^2_{\Gamma}+\frac{1}{\ell}a_m(\mathbf{u}_{m,h}^n-\mathbf{u}_{m,h}^{n-1},\mathbf{u}_{m,h}^n)\\
	 	&\qquad+\frac{c_0}{\ell}(p_{m,0}^n-p_{m,0}^{n-1},p_{m,0}^n)_{\Omega_m}+\frac{1}{\lambda\ell}((\alpha p_{m,0}^n-\phi_{m,h}^n)-(\alpha p_{m,0}^{n-1}-\phi_{m,h}^{n-1}),\alpha p_{m,0}^n-\phi_{m,h}^n)_{\Omega_m}=0.
	 \end{align*}
     Employing the inequality $\frac{1}{2}(\|a\|^2-\|b\|^2)\leq(a-b,a)$ to get
     \begin{align}
    	&2\mu_f\trb{\mathbf{u}_{f,h}^n}_f^2+\kappa\trb{p_{m,h}^n}_p^2+\beta_{\Gamma}\|(\mathbf{u}_{f,b}^n-\bar{\partial}_t \mathbf{u}_{m,b}^n)\cdot\bm{\tau}_{\Gamma}\|^2_{\Gamma}+\frac{1}{2\ell}\bigl(2\mu_m\trb{\mathbf{u}_{m,h}^n}_m^2-2\mu_m\trb{\mathbf{u}_{m,h}^{n-1}}_m^2\bigr)\nonumber\\
     	&\qquad+\frac{c_0}{2\ell}\bigl(\|p_{m,0}^n\|^2_{\Omega_m}-\|p_{m,0}^{n-1}\|^2_{\Omega_m}\bigr)+\frac{1}{2\lambda\ell}\bigl(\|\alpha p_{m,0}^n-\phi_{m,h}^n\|^2_{\Omega_m}-\|\alpha p_{m,0}^{n-1}-\phi_{m,h}^{n-1}\|^2_{\Omega_m}\bigr)\leq 0.\label{full-unique-proof1}
    \end{align}
    Setting $n=1,2,\dots,\bar{n}$ in \eqref{full-unique-proof1} and summing them together, we have
    \begin{align*}
    	&\sum_{n=1}^{\bar{n}}\bigl(2\mu_f\trb{\mathbf{u}_{f,h}^n}_f^2+\kappa\trb{p_{m,h}^n}_p^2+\beta_{\Gamma}\|(\mathbf{u}_{f,b}^n-\bar{\partial}_t \mathbf{u}_{m,b}^n)\cdot\bm{\tau}_{\Gamma}\|_{\Gamma}^2\bigr)+\frac{\mu_m}{\ell}\trb{\mathbf{u}_{m,h}^{\bar{n}}}_m^2+\frac{c_0}{2\ell}\|p_{m,0}^{\bar{n}}\|^2_{\Omega_m}\\
    	&\qquad+\frac{1}{2\lambda\ell}\|\alpha p_{m,0}^{\bar{n}}-\phi_{m,h}^{\bar{n}}\|^2_{\Omega_m}\leq 0.
    \end{align*}
    Since each term appearing on the left-hand side is nonnegative, we know 
    \begin{align}
    	\trb{\mathbf{u}_{f,h}^{\bar{n}}}_f=0,\quad\trb{p_{m,h}^{\bar{n}}}_p=0,\quad\trb{\mathbf{u}_{m,h}^{\bar{n}}}_m=0,\quad\|\alpha p_{m,0}^{\bar{n}}-\phi_{m,h}^{\bar{n}}\|_{\Omega_m}=0.\label{full-unique-proof2}
    \end{align}
    Thus,
    \[\mathbf u_{f,h}^{\bar n}=\mathbf 0,\quad p_{m,h}^{\bar n}=0,\quad \mathbf u_{m,h}^{\bar n}=\mathbf 0,\quad \phi_{m,h}^{\bar{n}}=0.\]
    Substituting the results \eqref{full-unique-proof2} back into \eqref{full-unique1} leads to $b_f(\mathbf{v}_{f,h},p_{f,h}^{\bar{n}})=0$. It follows from the inf-sup condition \eqref{inf-sup-Stokes} that $p_{f,h}^{\bar{n}}=0$. The proof is complete.
\end{proof}

\section{Error estimates}\label{section:error-estimate}
In this section, we are going to give error equations and derive error estimates for the semi-discrete and fully discrete schemes. Let $(\bu_f,\,p_f,\,\mathbf{u}_m,\,p_m,\,\phi_m)$ be the solution of the coupled problem \eqref{problem-eq1-s}--\eqref{problem-initial}, and $(\bu_{f,h},\,p_{f,h},\,\mathbf{u}_{m,h},\,p_{m,h},$ $\phi_{m,h})$ be the solution of the numerical scheme \eqref{semi-scheme1}--\eqref{semi-scheme5}. Let $\mathbf{e}_{f,h}=\mathbf{Q}_h \mathbf{u}_f-\mathbf{u}_{f,h},\,\varepsilon_{f,h}=\widehat{Q}_h p_f-p_{f,h},\,\mathbf{e}_{m,h}=\mathbf{Q}_h \mathbf{u}_m-\mathbf{u}_{m,h},\,\varepsilon_{m,h}=Q_h p_m-p_{m,h}$ and $\varepsilon_{\phi,h}=\widehat{Q}_h \phi_m-\phi_{m,h}$. 
\begin{lemma}
	For all $\mathbf{v}_{f,h}\in V_h^{f,D},\,q_{f,h}\in W_h^f,\,\mathbf{v}_{m,h}\in V_h^{m,D},\,q_{m,h}\in Z_h^D$ and $\psi_{m,h}\in W_h^m$, the error functions $\mathbf{e}_{f,h},\,\varepsilon_{f,h},\,\mathbf{e}_{m,h},\,\varepsilon_{m,h}$ and $\varepsilon_{\phi,h}$ satisfy 
	\begin{align}
		&a_f\left(\mathbf{e}_{f,h}, \mathbf{v}_{f,h}\right)-b_f(\mathbf{v}_{f,h}, \varepsilon_{f,h})+a_{\Gamma}(\mathbf{e}_{f,h},\partial_t \mathbf{e}_{m,h};\mathbf{v}_{f,h}, \mathbf{0})+b_{\Gamma}(\mathbf{0}, \varepsilon_{m,h};\mathbf{v}_{f,h},0)\nonumber\\
		&\qquad=s_f(\mathbf{Q}_h \mathbf{u}_f, \mathbf{v}_{f,h})-l_{\mathbf{u},f}(\mathbf{u}_f, \mathbf{v}_{f,h})+l_{ p,f}( p_f,\mathbf{v}_{f,h}),\label{semi-error-equation1}\\
    	&b_f(\mathbf{e}_{f,h},q_{f,h})=0,\label{semi-error-equation2}\\
		&a_m(\mathbf{e}_{m,h}, \mathbf{v}_{m,h})-b_m(\mathbf{v}_{m,h}, \varepsilon_{\phi,h})+a_{\Gamma}(\mathbf{e}_{f,h},\partial_t \mathbf{e}_{m,h}; \mathbf{0}, \mathbf{v}_{m,h})+c_{\Gamma}(\mathbf{0},\varepsilon_{m,h};\mathbf{v}_{m,h},0)\nonumber\\
		&\qquad=s_m(\mathbf{Q}_h \mathbf{u}_m, \mathbf{v}_{m,h})-l_{\mathbf{u},m}(\mathbf{u}_m, \mathbf{v}_{m,h})+l_{\phi,m}(\phi_m,\mathbf{v}_{m,h}),\label{semi-error-equation3}\\
		&\Bigl(c_0+\frac{\alpha^2}{\lambda}\Bigr)(\partial_t \varepsilon_{m,0},q_{m,0})_{\Omega_m}-\frac{\alpha}{\lambda}(\partial_t \varepsilon_{\phi,h},q_{m,0})_{\Omega_m}+a_p(\varepsilon_{m,h},q_{m,h})		+b_{\Gamma}(\mathbf{e}_{f,h},0; \mathbf{0},q_{m,h})\nonumber\\
		&\qquad+c_{\Gamma}(\partial_t \mathbf{e}_{m,h},0;\mathbf{0},q_{m,h})=s_p(Q_h p_m, q_{m,h})-l_{\phi_t,m}(\partial_t \phi_m, q_{m,h})-l_p(p_m,q_{m,h}),\label{semi-error-equation4}\\
		&\frac{1}{\lambda}(\varepsilon_{\phi,h}, \psi_{m,h})-\frac{\alpha}{\lambda}(\varepsilon_{m,0}, \psi_{m,h})+b_m(\mathbf{e}_{m,h}, \psi_{m,h}) =0,\label{semi-error-equation5}
	\end{align}
	where
	\begin{align*}
	   	& l_{\mathbf{u},i}(\mathbf{u}_i, \mathbf{v}_{i,h})=\sum_{T \in \mathcal{T}_h^i}2 \mu_i\langle(\mathbb{Q}_h \mathbf{D}(\mathbf{u}_i)-\mathbf{D}(\mathbf{u}_i)) \mathbf{n}, \mathbf{v}_{i,0}-\mathbf{v}_{i,b}\rangle_{\partial T},\quad i=f,m,\\
	   	& l_p(p_m, q_{m,h})=\sum_{T \in \mathcal{T}_h^m} \kappa\langle(\widehat{\mathbf{Q}}_h(\nabla p_m)-\nabla p_m) \cdot \mathbf{n}, q_{m,0}-q_{m,b}\rangle_{\partial T},\\
	   	& l_{p,f}(p_f,\mathbf{v}_{f,h})=\sum_{T \in \mathcal{T}_h^f}\langle(\widehat{Q}_h p_f-p_f)\mathbf{n}, \mathbf{v}_{f,0}-\mathbf{v}_{f,b} \rangle_{\partial T},\\
	   	& l_{\phi,m}(\phi_m,\mathbf{v}_{m,h})=\sum_{T \in \mathcal{T}_h^m}\langle(\widehat{Q}_h \phi_m-\phi_m)\mathbf{n}, \mathbf{v}_{m,0}-\mathbf{v}_{m,b} \rangle_{\partial T},\\
	   	& l_{\phi_t,m}(\partial_t \phi_m,q_{m,h})=\sum_{T \in \mathcal{T}_h^m} \frac{\alpha}{\lambda}(\widehat{Q}_h(\partial_t \phi_m)-\partial_t \phi_m, q_{m,0})_T.
	\end{align*}
\end{lemma}	
\begin{proof}
	Testing \eqref{problem-eq1-s} with $\mathbf{v}_{f,0}$, where $\mathbf{v}_{f,h}=\{\mathbf{v}_{f,0},\mathbf{v}_{f,b}\}\in V_h^{f,D}$, we have
	\begin{align*}
		(-\nabla\cdot\bm{\sigma}_f(\mathbf{u}_f,p_f),\mathbf{v}_{f,0})_{\Omega_f}=(\mathbf{f}_f,\mathbf{v}_{f,0})_{\Omega_f}.
	\end{align*}
    Using integration by parts to get
    \begin{align*}
    	(-\nabla\cdot\bm{\sigma}_f(\mathbf{u}_f,p_f),\mathbf{v}_{f,0})_{\Omega_f}=&\sum_{T \in \mathcal{T}_h^f}\bigl((\bm{\sigma}_f(\mathbf{u}_f,p_f),\nabla\mathbf{v}_{f,0})_T-\langle\bm{\sigma}_f(\mathbf{u}_f,p_f)\mathbf{n},\mathbf{v}_{f,0}\rangle_{\partial T}\bigr)\\
    	=&\sum_{T \in \mathcal{T}_h^f}\bigl((2\mu_f \mathbf{D}(\mathbf{u}_f),\mathbf{D}(\mathbf{v}_{f,0}))_T-(\nabla\cdot\mathbf{v}_{f,0},p_f)_T+\langle p_f\mathbf{n},\mathbf{v}_{f,0}-\mathbf{v}_{f,b}\rangle_{\partial T}\\
    	&-\langle2\mu_f \mathbf{D}(\mathbf{u}_f)\mathbf{n},\mathbf{v}_{f,0}-\mathbf{v}_{f,b}\rangle_{\partial T}\bigr)-\sum_{e\in\mathcal{E}_h^I} \langle\bm{\sigma}_f(\mathbf{u}_f,p_f)\mathbf{n}_{\Gamma},\mathbf{v}_{f,b}\rangle_e,
    \end{align*}
    \textcolor{black}{where we have used} $\sum\limits_{T \in \mathcal{T}_h^f}\langle2\mu_f \mathbf{D}(\mathbf{u}_f)\mathbf{n},\mathbf{v}_{f,b}\rangle_{\partial T\setminus(\partial T\cap\Gamma)}=0$ and $\sum\limits_{T \in \mathcal{T}_h^f}\langle p_f\mathbf{n},\mathbf{v}_{f,b}\rangle_{\partial T\setminus(\partial T\cap\Gamma)}=0$. \textcolor{black}According to the commutative property \eqref{weak-gradient-property-tensor}, \textcolor{black}{ we obtain}
    \begin{align*}
    	\sum_{T \in \mathcal{T}_h^f}(2\mu_f \mathbf{D}(\mathbf{u}_f),\mathbf{D}(\mathbf{v}_{f,0}))_T&=\sum_{T \in \mathcal{T}_h^f}(2\mu_f\mathbb{Q}_h \mathbf{D}(\mathbf{u}_f),\nabla\mathbf{v}_{f,0})_T\\
    	&=\sum_{T \in \mathcal{T}_h^f}2\mu_f\bigl(-(\nabla\cdot\mathbb{Q}_h\mathbf{D}(\mathbf{u}_f),\mathbf{v}_{f,0})_T+\langle\mathbb{Q}_h\mathbf{D}(\mathbf{u}_f)\mathbf{n},\mathbf{v}_{f,0}\rangle_{\partial T}\bigr)\\
    	&=\sum_{T \in \mathcal{T}_h^f}2\mu_f\bigl((\nabla_w\mathbf{v}_{f,h},\mathbb{Q}_h\mathbf{D}(\mathbf{u}_f))_T+\langle\mathbb{Q}_h\mathbf{D}(\mathbf{u}_f)\mathbf{n},\mathbf{v}_{f,0}-\mathbf{v}_{f,b}\rangle_{\partial T}\bigr)\\
    	&=\sum_{T \in \mathcal{T}_h^f}2\mu_f\bigl((\mathbf{D}_w(\mathbf{Q}_h \mathbf{u}_f),\mathbf{D}_w(\mathbf{v}_{f,h}))_T+\langle\mathbb{Q}_h\mathbf{D}(\mathbf{u}_f)\mathbf{n},\mathbf{v}_{f,0}-\mathbf{v}_{f,b}\rangle_{\partial T}\bigr).
    \end{align*}
    Similarly, we have
    \begin{align*}
    	-\sum_{T \in \mathcal{T}_h^f}(\nabla\cdot\mathbf{v}_{f,0},p_f)_T&=-\sum_{T \in \mathcal{T}_h^f}(\nabla\cdot\mathbf{v}_{f,0},\widehat{Q}_h p_f)_T\\
    	&=-\sum_{T \in \mathcal{T}_h^f}\bigl((\nabla_w\cdot\mathbf{v}_{f,h},\widehat{Q}_h p_f)_T+\langle\mathbf{v}_{f,0}-\mathbf{v}_{f,b},(\widehat{Q}_h p_f)\mathbf{n}\rangle_{\partial T}\bigr).
    \end{align*}
    For the interface terms, combining the interface conditions \eqref{problem-interface2} and \eqref{problem-interface4} together with the properties of the projection operators $\mathbf{Q}_h $ and $Q_h $ gives
    \begin{align*}
    	-\sum_{e\in\mathcal{E}_h^I} \langle\bm{\sigma}_f(\mathbf{u}_f,p_f)\mathbf{n}_{\Gamma},\mathbf{v}_{f,b}\rangle_e&=-\sum_{e\in\mathcal{E}_h^I} \langle(\bm{\sigma}_f(\mathbf{u}_f,p_f)\mathbf{n}_{\Gamma}\cdot\mathbf{n}_{\Gamma})\mathbf{n}_{\Gamma}+(\bm{\sigma}_f(\mathbf{u}_f,p_f)\mathbf{n}_{\Gamma}\cdot\bm{\tau}_{\Gamma})\bm{\tau}_{\Gamma},\mathbf{v}_{f,b}\rangle_e\\
    	&=\sum_{e\in\mathcal{E}_h^I}\bigl(\langle p_m\mathbf{n}_{\Gamma},\mathbf{v}_{f,b}\rangle_e+\beta_{\Gamma}\langle(\mathbf{u}_f-\partial_t \mathbf{u}_m)\cdot\bm{\tau}_{\Gamma},\mathbf{v}_{f,b}\cdot\bm{\tau}_{\Gamma}\rangle_e\bigr)\\
    	&=\sum_{e\in\mathcal{E}_h^I}\bigl(\langle (Q_b p_m )\mathbf{n}_{\Gamma},\mathbf{v}_{f,b}\rangle_e+\beta_{\Gamma}\langle(\mathbf{Q}_b \mathbf{u}_f-\mathbf{Q}_b(\partial_t \mathbf{u}_m))\cdot\bm{\tau}_{\Gamma},\mathbf{v}_{f,b}\cdot\bm{\tau}_{\Gamma}\rangle_e\bigr).
    \end{align*}
    Now, we arrive at
    \begin{align}
    	&a_f(\mathbf{Q}_h \mathbf{u}_f,\mathbf{v}_{f,h})-s_f(\mathbf{Q}_h \mathbf{u}_f,\mathbf{v}_{f,h})+\sum_{T \in \mathcal{T}_h^f}2\mu_f\langle(\mathbb{Q}_h\mathbf{D}(\mathbf{u}_f)-\mathbf{D}(\mathbf{u}_f))\mathbf{n},\mathbf{v}_{f,0}-\mathbf{v}_{f,b}\rangle_{\partial T}\nonumber\\
    	&\qquad-b_f(\mathbf{v}_{f,h},\widehat{Q}_h p_f)-\sum_{T \in \mathcal{T}_h^f}\langle(\widehat{Q}_h p_f-p_f)\mathbf{n},\mathbf{v}_{f,0}-\mathbf{v}_{f,b}\rangle_{\partial T}+a_{\Gamma}(\mathbf{Q}_h \mathbf{u}_f,\mathbf{Q}_h (\partial_t \mathbf{u}_m);\mathbf{v}_{f,h},\mathbf{0})\nonumber\\
    	&\qquad+b_{\Gamma}(\mathbf{0},Q_h p_m;\mathbf{v}_{f,h},0)=(\mathbf{f}_f,\mathbf{v}_{f,0})_{\Omega_f}.\label{error-equation-proof1}
    \end{align}
    Subtracting \eqref{semi-scheme1} from \eqref{error-equation-proof1}, we get \eqref{semi-error-equation1}.
    
    Testing \eqref{problem-eq2-s} with $q_{f,h}\in W_h^f$, we get
    \begin{align*}
    	(\nabla\cdot\bu_f,q_{f,h})_{\Omega_f}=0.
    \end{align*}
   Applying the property of the projection operator $\widehat{Q}_h$ and the commutative property \eqref{weak-divergence-property-vector} to get
    \begin{align*}
   	    (\nabla\cdot\bu_f,q_{f,h})_{\Omega_f}=\sum_{T \in \mathcal{T}_h^f}\bigl(\widehat{Q}_h(\nabla\cdot\bu_f),q_{f,h}\bigr)_T=\sum_{T \in \mathcal{T}_h^f}\bigl(\nabla_w\cdot(\mathbf{Q}_h \bu_f),q_{f,h}\bigr)_T,
       \end{align*}
   therefore,
   \begin{align}
   	   b_f(\mathbf{Q}_h \bu_f,q_{f,h})=0.\label{error-equation-proof2}
    \end{align}
    Subtracting \eqref{semi-scheme2} from \eqref{error-equation-proof2} yields  \eqref{semi-error-equation2}.
   
    Next, we test \eqref{problem-eq1-phi-b} with $\mathbf{v}_{m,0}$, where $\mathbf{v}_{m,h}=\{\mathbf{v}_{m,0},\mathbf{v}_{m,b}\}\in V_h^{m,D}$, to get
    \begin{align*}
   	    (-\nabla\cdot(2\mu_m \mathbf{D}(\mathbf{u}_m)-\phi_m \mathbf{I}),\mathbf{v}_{m,0})_{\Omega_m}=(\mathbf{f}_m,\mathbf{v}_{m,0})_{\Omega_m}.
    \end{align*}
    By an argument analogous to that leading to \eqref{error-equation-proof1}, we arrive at
    \begin{align}
   	    &a_m(\mathbf{Q}_h \mathbf{u}_m,\mathbf{v}_{m,h})-s_m(\mathbf{Q}_h \mathbf{u}_m,\mathbf{v}_{m,h})+\sum_{T \in \mathcal{T}_h^m}2\mu_m\langle(\mathbb{Q}_h\mathbf{D}(\mathbf{u}_m)-\mathbf{D}(\mathbf{u}_m))\mathbf{n},\mathbf{v}_{m,0}-\mathbf{v}_{m,b}\rangle_{\partial T}\nonumber\\
   	    &\qquad-b_m(\mathbf{v}_{m,h},\widehat{Q}_h \phi_m)-\sum_{T \in \mathcal{T}_h^m}\langle(\widehat{Q}_h \phi_m-\phi_m)\mathbf{n},\mathbf{v}_{m,0}-\mathbf{v}_{m,b}\rangle_{\partial T}+a_{\Gamma}(\mathbf{Q}_h \mathbf{u}_f,\mathbf{Q}_h (\partial_t \mathbf{u}_m);\mathbf{0},\mathbf{v}_{m,h})\nonumber\\
   	    &\qquad+c_{\Gamma}(\mathbf{0},Q_h p_m;\mathbf{v}_{m,h},0)=(\mathbf{f}_m,\mathbf{v}_{m,0})_{\Omega_m}.\label{error-equation-proof3}
    \end{align}
    Subtracting \eqref{semi-scheme3} from the above identity leads to \eqref{semi-error-equation3}.
   
    Testing \eqref{problem-eq2-phi-b} with $q_{m,0}$, where $q_{m,h}=\{q_{m,0},q_{m,b}\}\in Z_h^D$, yields
    \begin{align*}
   	    \biggl(\Bigl(c_0+\frac{\alpha^2}{\lambda}\Bigr)\partial_t  p_m-\frac{\alpha}{\lambda}\partial_t \phi_m-\nabla\cdot(\kappa\nabla p_m),q_{m,0}\biggr)_{\Omega_m}=(g,q_{m,0})_{\Omega_m}.
    \end{align*}
    Making use of the property of the projection operator $Q_0$, we derive
    \begin{align*}
    	\biggl(\Bigl(c_0+\frac{\alpha^2}{\lambda}\Bigr)\partial_t  p_m-\frac{\alpha}{\lambda}\partial_t \phi_m,q_{m,0}\biggr)_{\Omega_m}=\ &c_0(Q_0(\partial_t  p_m),q_{m,0})_{\Omega_m}+\frac{\alpha}{\lambda}(\alpha Q_0(\partial_t p_m)-\widehat{Q}_h(\partial_t \phi_m),q_{m,0})_{\Omega_m}\\
    	&+\frac{\alpha}{\lambda}(\widehat{Q}_h(\partial_t \phi_m)-\partial_t \phi_m,q_{m,0})_{\Omega_m}.
    \end{align*}
    It follows from $\sum\limits_{T \in \mathcal{T}_h^m}\langle(\kappa\nabla p_m)\cdot\mathbf{n},q_{m,b}\rangle_{\partial T\setminus(\partial T\cap\Gamma)}=0$, the interface condition \eqref{problem-interface1}, and the commutative property \eqref{nabla-w-property-scale} that
    \begin{align*}
    	(-\nabla\cdot(\kappa\nabla p_m),q_{m,0})_{\Omega_m}=&\sum_{T \in \mathcal{T}_h^m}\bigl((\widehat{\mathbf{Q}}_h(\kappa\nabla p_m),\nabla q_{m,0})_T-\langle(\kappa\nabla p_m)\cdot\mathbf{n},q_{m,0}-q_{m,b}\rangle_{\partial T}\bigr)\\
    	&+\sum_{e\in\mathcal{E}_h^I}\langle(\kappa\nabla p_m)\cdot\mathbf{n}_{\Gamma},q_{m,b}\rangle_e\\
    	=&\sum_{T \in \mathcal{T}_h^m}\bigl(-(\nabla\cdot\widehat{\mathbf{Q}}_h(\kappa\nabla p_m), q_{m,0})_T+\langle\widehat{\mathbf{Q}}_h(\kappa\nabla p_m)\cdot\mathbf{n},q_{m,0}\rangle_{\partial T}\\
    	&-\langle(\kappa\nabla p_m)\cdot\mathbf{n},q_{m,0}-q_{m,b}\rangle_{\partial T}\bigr)+\sum_{e\in\mathcal{E}_h^I}\langle(\kappa\nabla p_m)\cdot\mathbf{n}_{\Gamma},q_{m,b}\rangle_e\\
    	=&\sum_{T \in \mathcal{T}_h^m}\bigl((\widehat{\nabla}_wq_{m,h}, \widehat{\mathbf{Q}}_h(\kappa\nabla p_m))_T+\langle\widehat{\mathbf{Q}}_h(\kappa\nabla p_m)\cdot\mathbf{n},q_{m,0}-q_{m,b}\rangle_{\partial T}\\
    	&-\langle(\kappa\nabla p_m)\cdot\mathbf{n},q_{m,0}-q_{m,b}\rangle_{\partial T}\bigr)-\sum_{e\in\mathcal{E}_h^I}\langle(\mathbf{u}_f-\partial_t \mathbf{u}_m)\cdot\mathbf{n}_{\Gamma},q_{m,b}\rangle_e\\
    	=&\sum_{T \in \mathcal{T}_h^m}\bigl(\kappa(\widehat{\nabla}_w(Q_h p_m),\widehat{\nabla}_wq_{m,h})_T+\kappa\langle(\widehat{\mathbf{Q}}_h(\nabla p_m)-\nabla p_m)\cdot\mathbf{n},q_{m,0}-q_{m,b}\rangle_{\partial T}\bigr)\\
    	&-\sum_{e\in\mathcal{E}_h^I}\langle(\mathbf{Q}_b \mathbf{u}_f)\cdot\mathbf{n}_{\Gamma},q_{m,b}\rangle_e+\sum_{e\in\mathcal{E}_h^I}\langle\mathbf{Q}_b(\partial_t \mathbf{u}_m)\cdot\mathbf{n}_{\Gamma},q_{m,b}\rangle_e.
    \end{align*}
    Hence, 
    \begin{align}
    	&\Bigl(c_0+\frac{\alpha^2}{\lambda}\Bigr)(Q_0(\partial_t  p_m),q_{m,0})_{\Omega_m}-\frac{\alpha}{\lambda}\bigl(\widehat{Q}_h(\partial_t \phi_m),q_{m,0}\bigr)_{\Omega_m}+a_p(Q_h p_m,q_{m,h})\nonumber\\
    	&\qquad+\frac{\alpha}{\lambda}\bigl(\widehat{Q}_h(\partial_t \phi_m)-\partial_t \phi_m,q_{m,0}\bigr)_{\Omega_m}+\sum_{T \in \mathcal{T}_h^m}\kappa\langle(\widehat{\mathbf{Q}}_h(\nabla p_m)-\nabla p_m)\cdot\mathbf{n},q_{m,0}-q_{m,b}\rangle_{\partial T}\nonumber\\
    	&\qquad-s_p(Q_h p_m,q_{m,h})+b_{\Gamma}(\mathbf{Q}_h \mathbf{u}_f,0;\mathbf{0},q_{m,h})+c_{\Gamma}(\mathbf{Q}_h (\partial_t \mathbf{u}_m),0;\mathbf{0},q_{m,h})=(g,q_{m,0})_{\Omega_m}.\label{error-equation-proof4}
    \end{align}
    Taking the difference of \eqref{error-equation-proof4} and \eqref{semi-scheme4} leads to \eqref{semi-error-equation4}.
   
    Finally, testing \eqref{problem-eq3-phi-b} by $\psi_{m,h}\in W_h^m$ yields
    \begin{align*}
    	\biggl(\frac{1}{\lambda}\phi_m-\frac{\alpha}{\lambda} p_m+\nabla\cdot\mathbf{u}_m,\psi_{m,h}\biggr)_{\Omega_m}=0.
    \end{align*}
    Applying properties of $\widehat{Q}_h$, $Q_0$ and $\mathbf{Q}_h $ and the commutative property \eqref{weak-divergence-property-vector}, we have
    \begin{align*}
    	\biggl(\frac{1}{\lambda}\phi_m-\frac{\alpha}{\lambda} p_m+\nabla\cdot\mathbf{u}_m,\psi_{m,h}\biggr)_{\Omega_m}&=\frac{1}{\lambda}\bigl(\widehat{Q}_h \phi_m,\psi_{m,h}\bigr)_{\Omega_m}-\frac{\alpha}{\lambda}(Q_0 p_m,\psi_{m,h})_{\Omega_m}+\sum_{T \in \mathcal{T}_h^m}\bigl(\widehat{Q}_h(\nabla\cdot\mathbf{u}_m),\psi_{m,h}\bigr)_T\nonumber\\
    	&=\frac{1}{\lambda}\bigl(\widehat{Q}_h \phi_m,\psi_{m,h}\bigr)_{\Omega_m}-\frac{\alpha}{\lambda}(Q_0 p_m,\psi_{m,h})_{\Omega_m}+\sum_{T \in \mathcal{T}_h^m}(\nabla_w\cdot(\mathbf{Q}_h \mathbf{u}_m),\psi_{m,h})_T,
    \end{align*}
    which means
    \begin{align}
    	\frac{1}{\lambda}\bigl(\widehat{Q}_h \phi_m,\psi_{m,h}\bigr)_{\Omega_m}-\frac{\alpha}{\lambda}(Q_0 p_m,\psi_{m,h})_{\Omega_m}+b_m(\mathbf{Q}_h \mathbf{u}_m,\psi_{m,h})=0.\label{error-equation-proof5}
    \end{align}
    Subtracting \eqref{semi-scheme5} from \eqref{error-equation-proof5} gives \eqref{semi-error-equation5}.  This proves the result.
\end{proof}
\begin{lemma}
	Assume that $\bu_f\in[H^{k+1}(\Omega_f)]^2,\,p_f\in H^k(\Omega_f),\,\mathbf{u}_m\in[H^{k+1}(\Omega_m)]^2,\,p_m\in H^{k+1}(\Omega_m),\,\phi_m\in H^k(\Omega_m)$ and $\partial_t \phi_m\in H^k(\Omega_m)$. Then, for all $\mathbf{v}_{f,h}\in V_h^f,\,\mathbf{v}_{m,h}\in V_h^m$ and $q_{m,h}\in Z_h$, the following estimates hold
	\begin{align}
		|s_i(\mathbf{Q}_h \bu_i,\mathbf{v}_{i,h})|&\leq Ch^k\|\bu_i\|_{k+1}\trb{\mathbf{v}_{i,h}}_i,\quad i=f,m,\label{error-term1}\\
		|l_{\mathbf{u},i}(\bu_i,\mathbf{v}_{i,h})|&\leq Ch^k\|\bu_i\|_{k+1}\trb{\mathbf{v}_{i,h}}_i,\quad i=f,m,\label{error-term2}\\
		|l_{p,f}(p_f,\mathbf{v}_{f,h})|&\leq Ch^k\|p_f\|_k\trb{\mathbf{v}_{f,h}}_f,\label{error-term3}\\
		|l_{\phi,m}(\phi_m,\mathbf{v}_{m,h})|&\leq Ch^k\|\phi_m\|_k\trb{\mathbf{v}_{m,h}}_m,\label{error-term-phi}\\
		|s_p(Q_h p_m,q_{m,h})|&\leq Ch^k\|p_m\|_{k+1}\trb{q_{m,h}}_p,\label{error-term4}\\
		|l_p(p_m,q_{m,h})|&\leq Ch^k\|p_m\|_{k+1}\trb{q_{m,h}}_p,\label{error-term5}\\
		|l_{\phi_t,m}(\partial_t \phi_m,q_{m,h})|&\leq Ch^k\|\partial_t \phi_m\|_k\trb{q_{m,h}}_p.\label{error-term6}
	\end{align}
\end{lemma}
\begin{proof}
	It follows from the projection properties of $\mathbf{Q}_0$ and $\mathbf{Q}_b$, along with the trace inequality \cite{MR3452926-WG-stokes-2016} and \eqref{norm-vector-1} that
	\begin{align*}
		|s_i(\mathbf{Q}_h \bu_i,\mathbf{v}_{i,h})|&=\left|2\mu_i\sum_{T \in \mathcal{T}_h^i}h_T^{-1}\langle\mathbf{Q}_0 \mathbf{u}_i-\mathbf{Q}_b \bu_i,\mathbf{v}_{i,0}-\mathbf{v}_{i,b}\rangle_{\partial T}\right|\\
		&\leq C\left|\sum_{T \in \mathcal{T}_h^i}h_T^{-1}\langle\mathbf{Q}_0 \mathbf{u}_i-\bu_i,\mathbf{v}_{i,0}-\mathbf{v}_{i,b}\rangle_{\partial T}\right|+C\left|\sum_{T \in \mathcal{T}_h^i}h_T^{-1}\langle\bu_i-\mathbf{Q}_b \bu_i,\mathbf{v}_{i,0}-\mathbf{v}_{i,b}\rangle_{\partial T}\right|\\
		&\leq C\left(\sum_{T \in \mathcal{T}_h^i}h_T^{-1}\|\mathbf{Q}_0 \mathbf{u}_i-\bu_i\|^2_{\partial T}\right)^{\frac{1}{2}}\left(\sum_{T \in \mathcal{T}_h^i}h_T^{-1}\|\mathbf{v}_{i,0}-\mathbf{v}_{i,b}\|^2_{\partial T}\right)^{\frac{1}{2}}\\
		&\leq C\left(\sum_{T \in \mathcal{T}_h^i}(h_T^{-2}\|\mathbf{Q}_0 \mathbf{u}_i-\bu_i\|^2_T+\|\nabla(\mathbf{Q}_0 \mathbf{u}_i-\bu_i)\|^2_T)\right)^{\frac{1}{2}}\trb{\mathbf{v}_{i,h}}_i\\
		&\leq Ch^k\|\bu_i\|_{k+1}\trb{\mathbf{v}_{i,h}}_i.
	\end{align*}
    The estimates \eqref{error-term2}--\eqref{error-term5} follow from the same argument. It remains to estimate \(l_{\phi_t,m}\). Using the projection inequality and \eqref{v0-vh1}, we have
    \begin{align*}
    	|l_{\phi_t,m}(\partial_t \phi_m,q_{m,h})|&=\left|\sum_{T \in \mathcal{T}_h^m}\frac{\alpha}{\lambda}(\widehat{Q}_h(\partial_t \phi_m)-\partial_t \phi_m,q_{m,0})_T\right|\\
    	&\leq C\left(\sum_{T \in \mathcal{T}_h^m}\|\widehat{Q}_h(\partial_t \phi_m)-\partial_t \phi_m\|_T^2\right)^{\frac{1}{2}}\left(\sum_{T \in \mathcal{T}_h^m}\|q_{m,0}\|_T^2\right)^{\frac{1}{2}}\\
    	&\leq Ch^k\|\partial_t \phi_m\|_k\trb{q_{m,h}}_p.
    \end{align*}
    The proof is complete.
\end{proof}
\begin{theorem}
	Assume that $(\bu_f,p_f,\mathbf{u}_m,p_m,\phi_m)$ is the solution of the problem \eqref{problem-eq1-s}--\eqref{problem-initial}, and satisfies $\bu_f\in[H^{k+1}(\Omega_f)]^2\cap V^f,\,p_f\in H^k(\Omega_f),\,\mathbf{u}_m\in[H^{k+1}(\Omega_m)]^2\cap V^m,\,\partial_t \mathbf{u}_m\in[H^{k+1}(\Omega_m)]^2,\,p_m\in H^{k+1}(\Omega_m)\cap Z,\,\phi_m\in H^k(\Omega_m),\,\partial_t \phi_m\in H^k(\Omega_m)$. Let $(\bu_{f,h},\,p_{f,h};\,\mathbf{u}_{m,h},\,p_{m,h},\,\phi_{m,h})$ be the numerical solution of the semi-discrete scheme \eqref{semi-scheme1}--\eqref{semi-scheme5}. Then, for any $0\le t\le T$, there exists a constant $C>0$, such that
	\begin{align}
		&\int_0^t\left(\trb{\mathbf{e}_{f,h}(r)}_f^2+\trb{\varepsilon_{m,h}(r)}_p^2+\|\varepsilon_{f,h}(r)\|_{\Omega_f}^2+\|\varepsilon_{\phi,h}(r)\|_{\Omega_m}^2+\|(\mathbf{e}_{f,b}(r)-\partial_t \mathbf{e}_{m,b}(r))\cdot\bm{\tau}_{\Gamma}\|^2_{\Gamma}\right)\,\mathrm{d}r\nonumber\\
		&\qquad+\trb{\mathbf{e}_{m,h}(t)}_m^2\leq Ch^{2k}\mathcal{R}(t),
	\end{align}
	where
	\begin{align*}
		\mathcal{R}(t)=\ &\|\mathbf{u}_m(t)\|_{k+1}^2+\|\phi_m(t)\|_k^2+\int_0^t\bigl(\|\mathbf{u}_f(r)\|_{k+1}^2+\|p_f(r)\|_k^2+\|\mathbf{u}_m(r)\|^2_{k+1}+\|\partial_t \mathbf{u}_m(r)\|_{k+1}^2\\
		&+\|p_m(r)\|_{k+1}^2+\|\phi_m(r)\|_k^2+\|\partial_t \phi_m(r)\|_k^2\bigr)\,\mathrm{d}r.
	\end{align*}
\end{theorem}
\begin{proof}
	Differentiating \eqref{semi-error-equation5} with respect to time yields the following identity
	\begin{align}
		\frac{1}{\lambda}(\partial_t \varepsilon_{\phi,h}, \psi_{m,h})-\frac{\alpha}{\lambda}(\partial_t \varepsilon_{m,0}, \psi_{m,h})+b_m(\partial_t \mathbf{e}_{m,h},\psi_{m,h})=0.\label{semi-error-equation6}
	\end{align}
	We take $\mathbf{v}_{f,h}=\mathbf{e}_{f,h},\,q_{f,h}=\varepsilon_{f,h},\,\mathbf{v}_{m,h}=\partial_t \mathbf{e}_{m,h},\,q_{m,h}=\varepsilon_{m,h}$ in the error equations \eqref{semi-error-equation1}--\eqref{semi-error-equation4}, and take $\psi_{m,h}=\varepsilon_{\phi,h}$ in \eqref{semi-error-equation6}, then add the resulting identities to get
	\begin{align}
		&2\mu_f\trb{\mathbf{e}_{f,h}}_f^2+\mu_m\frac{\mathrm{d}}{\mathrm{d}t}\trb{\mathbf{e}_{m,h}}_m^2+\frac{c_0}{2}\frac{\mathrm{d}}{\mathrm{d}t}\|\varepsilon_{m,0}\|^2_{\Omega_m}+\frac{1}{2\lambda}\frac{\mathrm{d}}{\mathrm{d}t}\|\alpha\varepsilon_{m,0}-\varepsilon_{\phi,h}\|^2_{\Omega_m}+\beta_{\Gamma}\|(\mathbf{e}_{f,b}-\partial_t \mathbf{e}_{m,b})\cdot\bm{\tau}_{\Gamma}\|^2_{\Gamma}\nonumber\\
		&\qquad+\kappa\trb{\varepsilon_{m,h}}_p^2=s_f(\mathbf{Q}_h \mathbf{u}_f, \mathbf{e}_{f,h})-l_{\mathbf{u},f}(\mathbf{u}_f, \mathbf{e}_{f,h})+l_{ p,f}(p_f,\mathbf{e}_{f,h})+s_m(\mathbf{Q}_h \mathbf{u}_m, \partial_t \mathbf{e}_{m,h})-l_{\mathbf{u},m}(\mathbf{u}_m, \partial_t \mathbf{e}_{m,h})\nonumber\\
		&\qquad+l_{\phi,m}(\phi_m,\partial_t \mathbf{e}_{m,h})+s_p(Q_h p_m, \varepsilon_{m,h})-l_{\phi_t,m}(\partial_t \phi_m, \varepsilon_{m,h})-l_p(p_m,\varepsilon_{m,h}).\label{semi-error-estimate1}
	\end{align}
	\textcolor{black}{The following identities hold:}
	\begin{align*}
		s_m(\mathbf{Q}_h \mathbf{u}_m, \partial_t \mathbf{e}_{m,h})&=\frac{\mathrm{d}}{\mathrm{d}t}s_m(\mathbf{Q}_h \mathbf{u}_m, \mathbf{e}_{m,h})-s_m( \mathbf{Q}_h(\partial_t \mathbf{u}_m), \mathbf{e}_{m,h}),\\
		l_{\mathbf{u},m}(\mathbf{u}_m, \partial_t \mathbf{e}_{m,h})&=\frac{\mathrm{d}}{\mathrm{d}t}l_{\mathbf{u},m}(\mathbf{u}_m,\mathbf{e}_{m,h})-l_{\mathbf{u},m}(\partial_t \mathbf{u}_m,\mathbf{e}_{m,h}),\\
		l_{\phi,m}(\phi_m,\partial_t \mathbf{e}_{m,h})&=\frac{\mathrm{d}}{\mathrm{d}t}l_{\phi,m}(\phi_m,\mathbf{e}_{m,h})-l_{\phi,m}(\partial_t\phi_m,\mathbf{e}_{m,h}).
	\end{align*}
	\textcolor{black}{Integrating \eqref{semi-error-estimate1} over the time interval (0,t), we obtain}
	\begin{align}
		&\int_0^t\bigl(2\mu_f\trb{\mathbf{e}_{f,h}(r)}_f^2+\kappa\trb{\varepsilon_{m,h}(r)}_p^2+\beta_{\Gamma}\|(\mathbf{e}_{f,b}(r)-\partial_t \mathbf{e}_{m,b}(r))\cdot\bm{\tau}_{\Gamma}\|^2_{\Gamma}\bigr)\,\mathrm{d}r+\mu_m\trb{\mathbf{e}_{m,h}(t)}_m^2\nonumber\\
		&\qquad+\frac{c_0}{2}\|\varepsilon_{m,0}(t)\|^2_{\Omega_m}+\frac{1}{2\lambda}\|\alpha\varepsilon_{m,0}(t)-\varepsilon_{\phi,h}(t)\|^2_{\Omega_m}=s_m(\mathbf{Q}_h \mathbf{u}_m(t), \mathbf{e}_{m,h}(t))-l_{\mathbf{u},m}(\mathbf{u}_m(t),\mathbf{e}_{m,h}(t))\nonumber\\
		&\qquad+l_{\phi,m}(\phi_m(t),\mathbf{e}_{m,h}(t))+\int_0^t(s_f(\mathbf{Q}_h \mathbf{u}_f(r),\mathbf{e}_{f,h}(r))-l_{\mathbf{u},f}(\mathbf{u}_f(r), \mathbf{e}_{f,h}(r))+l_{ p,f}(p_f(r),\mathbf{e}_{f,h}(r))\nonumber\\
		&\qquad-s_m( \mathbf{Q}_h(\partial_t \mathbf{u}_m(r)), \mathbf{e}_{m,h}(r))+l_{\mathbf{u},m}(\partial_t \mathbf{u}_m(r), \mathbf{e}_{m,h}(r))-l_{\phi,m}(\partial_t \phi_m(r),\mathbf{e}_{m,h}(r))\nonumber\\
		&\qquad+s_p(Q_h p_m(r), \varepsilon_{m,h}(r))-l_{\phi_t,m}(\partial_t \phi_m(r), \varepsilon_{m,h}(r))-l_p(p_m(r),\varepsilon_{m,h}(r)))\,\mathrm{d}r,\label{semi-error-estimate2}
	\end{align}
	where we have used initial conditions $\mathbf{e}_{m,h}(0)=\mathbf{0},\,\varepsilon_{m,h}(0)=0$ and $\varepsilon_{\phi,h}(0)=0$. Denote the right-hand side of identity \eqref{semi-error-estimate2} by $I_{\rm rhs}(t)$. Define
	\begin{align*}
		\mathcal{A}(t)=\ &\int_0^t(\|\mathbf{u}_f(r)\|_{k+1}^2+\|p_f(r)\|_k^2+\|\partial_t \mathbf{u}_m(r)\|_{k+1}^2+\|p_m(r)\|_{k+1}^2+\|\partial_t \phi_m(r)\|_k^2)\,\mathrm{d}r,\\
		\mathcal{B}(t)=\ &\|\mathbf{u}_m(t)\|_{k+1}^2+\|\phi_m(t)\|_k^2+\mathcal{A}(t).
	\end{align*}
	Then, using the estimates \eqref{error-term1}--\eqref{error-term6}, and Young's inequality, we get
	\begin{align*}
		I_{\rm rhs}(t)\leq\ & Ch^k\biggl((\|\mathbf{u}_m(t)\|_{k+1}+\|\phi_m(t)\|_k)\trb{\mathbf{e}_{m,h}(t)}_m+\int_0^t\bigl((\|\mathbf{u}_f(r)\|_{k+1}+\|p_f(r)\|_k)\trb{\mathbf{e}_{f,h}(r)}_f\\
		&+(\|\partial_t \mathbf{u}_m(r)\|_{k+1}+\|\partial_t \phi_m(r)\|_k)\trb{\mathbf{e}_{m,h}(r)}_m+(\|p_m(r)\|_{k+1}+\|\partial_t \phi_m(r)\|_k)\trb{\varepsilon_{m,h}(r)}_p\bigr)\,\mathrm{d}r\biggr)\\
		\leq\ & Ch^{2k}\mathcal{B}(t)+\frac{\mu_m}{2}\trb{\mathbf{e}_{m,h}(t)}_m^2+\int_0^t\biggl(\frac{\mu_m}{2}\trb{\mathbf{e}_{m,h}(r)}_m^2+2\mu_f\trb{\mathbf{e}_{f,h}(r)}_f^2+\kappa\trb{\varepsilon_{m,h}(r)}_p^2\biggr)\,\mathrm{d}r,
	\end{align*}
	which implies
	\begin{align*}
		\frac{\mu_m}{2}\trb{\mathbf{e}_{m,h}(t)}_m^2\leq\ & Ch^{2k}\mathcal{B}(t)+\int_0^t\frac{\mu_m}{2}\trb{\mathbf{e}_{m,h}(r)}_m^2\,\mathrm{d}r.
	\end{align*}
	Using the Gronwall inequality yields
	\begin{align}
		\trb{\mathbf{e}_{m,h}(t)}_m^2\leq\ Ch^{2k}\mathcal{B}(t).\label{semi-error-estimate3}
	\end{align}
   
	Using Young's inequality in \eqref{semi-error-estimate2} with a different choice of coefficients gives
	\begin{align*}
		I_{\rm rhs}(t)\leq\ & Ch^{2k}\mathcal{B}(t)+\mu_m\trb{\mathbf{e}_{m,h}(t)}_m^2+\int_0^t\biggl(\mu_m\trb{\mathbf{e}_{m,h}(r)}_m^2+\mu_f\trb{\mathbf{e}_{f,h}(r)}_f^2+\frac{\kappa}{2}\trb{\varepsilon_{m,h}(r)}_p^2\biggr)\,\mathrm{d}r.
	\end{align*}
    Then, we get
    \begin{align*}
    	\int_0^t\biggl(\mu_f\trb{\mathbf{e}_{f,h}(r)}_f^2+\frac{\kappa}{2}\trb{\varepsilon_{m,h}(r)}_p^2+\beta_{\Gamma}\|(\mathbf{e}_{f,b}(r)-\partial_t \mathbf{e}_{m,b}(r))\cdot\bm{\tau}_{\Gamma}\|^2_{\Gamma}\biggr)\,\mathrm{d}r \leq Ch^{2k}\mathcal{B}(t)+\int_0^t\mu_m\trb{\mathbf{e}_{m,h}(r)}_m^2\,\mathrm{d}r.
    \end{align*}
    According to \eqref{semi-error-estimate3}, we know
    \begin{align}
    	\trb{\mathbf{e}_{m,h}(r)}_m^2\leq Ch^{2k}\mathcal{B}(r)\leq Ch^{2k}\bigl(\|\mathbf{u}_m(r)\|_{k+1}^2+\|\phi_m(r)\|_k^2+\mathcal{A}(t)\bigr),\label{semi-error-estimate4}
    \end{align}
    where $0<r\leq t$. Integrating \eqref{semi-error-estimate4} \textcolor{black}{over the interval $(0,t)$} and using $0\le t\le T$, we have
    \begin{align}
    	\int_0^t\trb{\mathbf{e}_{m,h}(r)}_m^2\,\mathrm{d}r&\leq Ch^{2k}\biggl(\int_0^t\bigl(\|\mathbf{u}_m(r)\|_{k+1}^2+\|\phi_m(r)\|_k^2+\mathcal{A}(t)\bigr)\,\mathrm{d}r\biggr)\nonumber\\
    	&\leq Ch^{2k}\biggl(\int_0^t\bigl(\|\mathbf{u}_m(r)\|_{k+1}^2+\|\phi_m(r)\|_k^2\bigr)\,\mathrm{d}r+\mathcal{A}(t)\biggr).\label{semi-error-estimate5}
    \end{align}
    Therefore,
    \begin{align}
    	&\int_0^t\bigl(\trb{\mathbf{e}_{f,h}(r)}_f^2+\trb{\varepsilon_{m,h}(r)}_p^2+\|(\mathbf{e}_{f,b}(r)-\partial_t \mathbf{e}_{m,b}(r))\cdot\bm{\tau}_{\Gamma}\|^2_{\Gamma}\bigr)\,\mathrm{d}r\leq Ch^{2k}\mathcal{R}(t).\label{semi-error-estimate6}
    \end{align}
   
    From the first error equation \eqref{semi-error-equation1}, we get
    \begin{align*}
    	b_f(\mathbf{v}_{f,h}, \varepsilon_{f,h})=\ &a_f\left(\mathbf{e}_{f,h}, \mathbf{v}_{f,h}\right)+a_{\Gamma}(\mathbf{e}_{f,h},\partial_t \mathbf{e}_{m,h};\mathbf{v}_{f,h}, \mathbf{0})+b_{\Gamma}(\mathbf{0}, \varepsilon_{m,h};\mathbf{v}_{f,h},0)-s_f(\mathbf{Q}_h \mathbf{u}_f, \mathbf{v}_{f,h})\\
    	&+l_{\mathbf{u},f}(\mathbf{u}_f, \mathbf{v}_{f,h})-l_{ p,f}(p_f,\mathbf{v}_{f,h}).
    \end{align*}
    \textcolor{black}{By applying the estimates
\eqref{vb-vh1}, \eqref{vb-vh2}, and the error bounds
\eqref{error-term1}--\eqref{error-term3},} we have
    \begin{align*}
    	b_f(\mathbf{v}_{f,h}, \varepsilon_{f,h})\leq\ & C\bigl(\trb{\mathbf{e}_{f,h}}_f\trb{\mathbf{v}_{f,h}}_f+\|(\mathbf{e}_{f,b}-\partial_t \mathbf{e}_{m,b})\cdot\bm{\tau}_{\Gamma}\|_{\Gamma}\|\mathbf{v}_{f,b}\|_{\Gamma}+\|\mathbf{v}_{f,b}\|_{\Gamma}\|\varepsilon_{m,b}\|_{\Gamma}\bigr)\\
    	&+Ch^k(\|\mathbf{u}_f\|_{k+1}+\|p_f\|_k)\trb{\mathbf{v}_{f,h}}_f\\
    	\leq\ & C\bigl(\trb{\mathbf{e}_{f,h}}_f\trb{\mathbf{v}_{f,h}}_f+\|(\mathbf{e}_{f,b}-\partial_t \mathbf{e}_{m,b})\cdot\bm{\tau}_{\Gamma}\|_{\Gamma}\trb{\mathbf{v}_{f,h}}_f+\trb{\mathbf{v}_{f,h}}_f\trb{\varepsilon_{m,h}}_p\bigr)\\
    	&+Ch^k(\|\mathbf{u}_f\|_{k+1}+\|p_f\|_k)\trb{\mathbf{v}_{f,h}}_f.
    \end{align*}
    \textcolor{black}{Employing} the inf-sup condition \eqref{inf-sup-Stokes} \textcolor{black}{gives}
    \begin{align*}
    	\|\varepsilon_{f,h}\|_{\Omega_f}^2&\leq C\bigl(\trb{\mathbf{e}_{f,h}}_f^2+\|(\mathbf{e}_{f,b}-\partial_t \mathbf{e}_{m,b})\cdot\bm{\tau}_{\Gamma}\|^2_{\Gamma}+\trb{\varepsilon_{m,h}}_p^2\bigr)+Ch^{2k}(\|\mathbf{u}_f\|^2_{k+1}+\|p_f\|^2_k).
    \end{align*}
    \textcolor{black}{ Integrating the above estimate over the interval $(0,t)$ and using \eqref{semi-error-estimate6}, we further obtain}
    \begin{align}
    	\int_0^t\|\varepsilon_{f,h}(r)\|_{\Omega_f}^2\,\mathrm{d}r\leq\ & C\int_0^t\bigl(\trb{\mathbf{e}_{f,h}(r)}_f^2+\|(\mathbf{e}_{f,b}(r)-\partial_t \mathbf{e}_{m,b}(r))\cdot\bm{\tau}_{\Gamma}\|^2_{\Gamma}+\trb{\varepsilon_{m,h}(r)}_p^2\bigr)\,\mathrm{d}r\nonumber\\
    	&+Ch^{2k}\int_0^t\bigl(\|\mathbf{u}_f(r)\|^2_{k+1}+\|p_f(r)\|^2_k\bigr)\,\mathrm{d}r\nonumber\\
    	\leq\ & Ch^{2k}\mathcal{R}(t).\label{semi-error-estimate7}
    \end{align}
    \textcolor{black}{Following the same procedure and combining the inf-sup condition \eqref{inf-sup-Biot} with \eqref{error-term4}--\eqref{error-term6} and \eqref{semi-error-estimate5}, we have}
    \begin{align*}
    	\int_0^t\|\varepsilon_{\phi,h}(r)\|_{\Omega_m}^2\,\mathrm{d}r\leq Ch^{2k}\mathcal{R}(t).
    \end{align*}
    The proof is complete.
\end{proof}
\textcolor{black}{Let $(\mathbf{u}_{f,h}^n,p_{f,h}^n;\mathbf{u}_{m,h}^n,p_{m,h}^n,\phi_{m,h}^n)$ denote the solution of the discrete WG scheme
\eqref{full-scheme1}--\eqref{full-scheme5} at the time level $t^n$. 
Define the corresponding error functions by} $\mathbf{e}_{f,h}^n=\mathbf{Q}_h \mathbf{u}_f^n-\mathbf{u}_{f,h}^n,\,\varepsilon_{f,h}^n=\widehat{Q}_h p_f^n-p_{f,h}^n,\,\mathbf{e}_{m,h}^n=\mathbf{Q}_h \mathbf{u}_m^n-\mathbf{u}_{m,h}^n,\,\varepsilon_{m,h}^n=Q_h p_m^n-p_{m,h}^n$ and $\varepsilon_{\phi,h}^n=\widehat{Q}_h \phi_m^n-\phi_{m,h}^n$. Then the following fully discrete error equations hold.
\begin{lemma}
	For all $\mathbf{v}_{f,h}\in V_h^{f,D},\,q_{f,h}\in W_h^f,\,\mathbf{v}_{m,h}\in V_h^{m,D},\,q_{m,h}\in Z_h^D$ and $\psi_{m,h}\in W_h^m$, the error functions $\mathbf{e}_{f,h}^n,\,\varepsilon_{f,h}^n,\,\mathbf{e}_{m,h}^n,\,\varepsilon_{m,h}^n$ and $\varepsilon_{\phi,h}^n$ satisfy 
	\begin{align}
		&a_f\left(\mathbf{e}_{f,h}^n, \mathbf{v}_{f,h}\right)-b_f(\mathbf{v}_{f,h}, \varepsilon_{f,h}^n)+a_{\Gamma}(\mathbf{e}_{f,h}^n,\bar{\partial}_t \mathbf{e}_{m,h}^n;\mathbf{v}_{f,h}, \mathbf{0})+b_{\Gamma}(\mathbf{0}, \varepsilon_{m,h}^n;\mathbf{v}_{f,h},0)\nonumber\\
		&\qquad=s_f(\mathbf{Q}_h \mathbf{u}_f^n, \mathbf{v}_{f,h})-l_{\mathbf{u},f}(\mathbf{u}_f^n, \mathbf{v}_{f,h})+l_{p,f}( p_f^n,\mathbf{v}_{f,h})+l_{\Gamma}^1(\mathbf{u}_m^n,\mathbf{v}_{f,h}),\label{full-error-equation1}\\
		&b_f(\mathbf{e}_{f,h}^n,q_{f,h})=0,\label{full-error-equation2}\\
		&a_m(\mathbf{e}_{m,h}^n, \mathbf{v}_{m,h})-b_m(\mathbf{v}_{m,h}, \varepsilon_{\phi,h}^n)+a_{\Gamma}(\mathbf{e}_{f,h}^n,\bar{\partial}_t \mathbf{e}_{m,h}^n ; \mathbf{0}, \mathbf{v}_{m,h})+c_{\Gamma}(\mathbf{0},\varepsilon_{m,h}^n;\mathbf{v}_{m,h},0)\nonumber\\
		&\qquad=s_m(\mathbf{Q}_h \mathbf{u}_m^n, \mathbf{v}_{m,h})-l_{\mathbf{u},m}(\mathbf{u}_m^n, \mathbf{v}_{m,h})+l_{\phi,m}(\phi_m^n,\mathbf{v}_{m,h})-l_{\Gamma}^2(\mathbf{u}_m^n,\mathbf{v}_{m,h}),\label{full-error-equation3}\\
		&\Bigl(c_0+\frac{\alpha^2}{\lambda}\Bigr)(\bar{\partial}_t \varepsilon_{m,0}^n,q_{m,0})_{\Omega_m}-\frac{\alpha}{\lambda}(\bar{\partial}_t \varepsilon_{\phi,h}^n,q_{m,0})_{\Omega_m}+a_p(\varepsilon_{m,h}^n,q_{m,h})+b_{\Gamma}(\mathbf{e}_{f,h}^n,0; \mathbf{0},q_{m,h})\nonumber\\
		&\qquad+c_{\Gamma}(\bar{\partial}_t \mathbf{e}_{m,h}^n,0;\mathbf{0},q_{m,h})=s_p(Q_h p_m^n, q_{m,h})-l_{\phi_t,m}((\partial_t \phi_m)^n, q_{m,h})-l_p(p_m^n,q_{m,h})\nonumber\\
		&\qquad-l_0^1(p_m^n,q_{m,0})+l_0^2(\phi_m^n,q_{m,0})-l_{\Gamma}^3(\mathbf{u}_m^n,q_{m,h}),\label{full-error-equation4}\\
		&\frac{1}{\lambda}(\varepsilon_{\phi,h}^n, \psi_{m,h})-\frac{\alpha}{\lambda}(\varepsilon_{m,0}^n, \psi_{m,h})+b_m(\mathbf{e}_{m,h}^n, \psi_{m,h})=0,\label{full-error-equation5}
	\end{align}
	where
	\begin{align*}
		&l_{\Gamma}^1(\mathbf{u}_m^n, \mathbf{v}_{f,h})=\sum_{e\in\mathcal{E}_h^I}\beta_{\Gamma}\langle((\partial_t \mathbf{u}_m)^n-\bar{\partial}_t \mathbf{u}_m^n)\cdot\bm{\tau}_{\Gamma},\mathbf{v}_{f,b}\cdot\bm{\tau}_{\Gamma}\rangle_e,\\
		&l_{\Gamma}^2(\mathbf{u}_m^n, \mathbf{v}_{m,h})=\sum_{e\in\mathcal{E}_h^I}\beta_{\Gamma}\langle((\partial_t \mathbf{u}_m)^n-\bar{\partial}_t \mathbf{u}_m^n)\cdot\bm{\tau}_{\Gamma},\mathbf{v}_{m,b}\cdot\bm{\tau}_{\Gamma}\rangle_e,\\
		&l_{\Gamma}^3(\mathbf{u}_m^n, q_{m,h})=\sum_{e\in\mathcal{E}_h^I}\langle((\partial_t \mathbf{u}_m)^n-\bar{\partial}_t \mathbf{u}_m^n)\cdot\mathbf{n}_{\Gamma},q_{m,b}\rangle_e,\\
		&l_0^1(p_m^n,q_{m,0})=\Bigl(c_0+\frac{\alpha^2}{\lambda}\Bigr)((\partial_t p_m)^n-\bar{\partial}_t p_m^n,q_{m,0})_{\Omega_m}, \\
		&l_0^2(\phi_m^n,q_{m,0})=\frac{\alpha}{\lambda}(\widehat{Q}_h((\partial_t \phi_m)^n-\bar{\partial}_t \phi_m^n),q_{m,0})_{\Omega_m}.
	\end{align*}
\end{lemma}	
\begin{proof}
	We take $t=t^n$ in \eqref{error-equation-proof1} to get
	\begin{align}
		(\mathbf{f}_f^n,\mathbf{v}_{f,0})_{\Omega_f}=&a_f(\mathbf{Q}_h \mathbf{u}_f^n,\mathbf{v}_{f,h})-s_f(\mathbf{Q}_h \mathbf{u}_f^n,\mathbf{v}_{f,h})+l_{\mathbf{u},f}(\mathbf{u}_f^n,\mathbf{v}_{f,h})-b_f(\mathbf{v}_{f,h},\widehat{Q}_h p_f^n)\nonumber\\
		&-l_{p,f}(p_f^n,\mathbf{v}_{f,h})+a_{\Gamma}(\mathbf{Q}_h \mathbf{u}_f^n,\mathbf{Q}_h (\partial_t \mathbf{u}_m)^n;\mathbf{v}_{f,h},\mathbf{0})+b_{\Gamma}(\mathbf{0},Q_h p_m^n;\mathbf{v}_{f,h},0).\label{full-error-equation-proof1}
	\end{align}
    Subtracting \eqref{full-scheme1} from \eqref{full-error-equation-proof1} leads to
    \begin{align}
    	&a_f(\mathbf{e}_{f,h}^n,\mathbf{v}_{f,h})-b_f(\mathbf{v}_{f,h},\varepsilon_{f,h}^n)+a_{\Gamma}(\mathbf{e}_{f,h}^n,\mathbf{Q}_h (\partial_t \mathbf{u}_m)^n-\bar{\partial}_t \mathbf{u}_{m,h}^n;\mathbf{v}_{f,h},\mathbf{0})+b_{\Gamma}(\mathbf{0},\varepsilon_{m,h}^n;\mathbf{v}_{f,h},0)\nonumber\\
    	=\ &s_f(\mathbf{Q}_h \mathbf{u}_f^n,\mathbf{v}_{f,h})-l_{\mathbf{u},f}(\mathbf{u}_f^n,\mathbf{v}_{f,h})+l_{p,f}(p_f^n,\mathbf{v}_{f,h}).\label{full-error-equation-proof2}
    \end{align}
    \textcolor{black}{By the projection property of $\mathbf{Q}_b$, we obtain}
    \begin{align}
    	&a_{\Gamma}(\mathbf{e}_{f,h}^n,\mathbf{Q}_h (\partial_t \mathbf{u}_m)^n-\bar{\partial}_t \mathbf{u}_{m,h}^n;\mathbf{v}_{f,h},\mathbf{0})\nonumber\\
    	=\ &\beta_{\Gamma}\left\langle\left(\mathbf{e}_{f,b}^n-\mathbf{Q}_b(\partial_t \mathbf{u}_m)^n+\bar{\partial}_t (\mathbf{Q}_b \mathbf{u}_m^n)-\bar{\partial}_t (\mathbf{Q}_b \mathbf{u}_m^n-\mathbf{u}_{m,b}^n)\right)\cdot\bm{\tau}_{\Gamma},\mathbf{v}_{f,b}\cdot\bm{\tau}_{\Gamma}\right\rangle_{\Gamma}\nonumber\\
    	=\ &a_{\Gamma}(\mathbf{e}_{f,h}^n,\bar{\partial}_t \mathbf{e}_{m,h}^n;\mathbf{v}_{f,h},\mathbf{0})-\beta_{\Gamma}\left\langle\left((\partial_t \mathbf{u}_m)^n-\bar{\partial}_t \mathbf{u}_m^n\right)\cdot\bm{\tau}_{\Gamma},\mathbf{v}_{f,b}\cdot\bm{\tau}_{\Gamma}\right\rangle_{\Gamma}.\label{full-error-equation-proof3}
    \end{align}
    Substituting \eqref{full-error-equation-proof3} into \eqref{full-error-equation-proof2}, we get the error equation \eqref{full-error-equation1}. Similarly, we obtain the remaining error equations \eqref{full-error-equation2}--\eqref{full-error-equation5}.
\end{proof}
\begin{theorem}
	Assume that the time step satisfies $0<\ell<1$ and that some
	regularity assumptions hold:
	\begin{align*}
		&\mathbf{u}_f\in L^{\infty}(0,T;[H^{k+1}(\Omega_f)]^2\cap V^f),
		\qquad
		p_f\in L^{\infty}(0,T;H^k(\Omega_f)),\\
		&\mathbf{u}_m\in L^{\infty}(0,T;[H^{k+1}(\Omega_m)]^2\cap V^m),
		\qquad
		\partial_t\mathbf{u}_m\in L^2(0,T;[H^{k+1}(\Omega_m)]^2),\\
		&\partial_{tt}\mathbf{u}_m\in L^2(0,T;[L^2(\Gamma)]^2),\\
		&p_m\in L^{\infty}(0,T;H^{k+1}(\Omega_m)\cap Z),
		\qquad
		\partial_{tt}p_m\in L^2(0,T;L^2(\Omega_m)),\\
		&\phi_m\in L^{\infty}(0,T;H^k(\Omega_m)),
		\qquad
		\partial_t \phi_m\in L^{\infty}(0,T;H^k(\Omega_m)),
		\qquad
		\partial_{tt}\phi_m\in L^2(0,T;L^2(\Omega_m)).
	\end{align*}
    Then, for each $1\le \bar n\le N$, there exists a constant $C>0$, independent of $h$, $\ell$, $\lambda$ and
    $c_0$, such that
    \begin{align}
    	&\ell\sum_{n=1}^{\bar{n}}\bigl(\trb{\mathbf{e}_{f,h}^n}_f^2+\trb{\varepsilon_{m,h}^n}_p^2+\|\varepsilon_{f,h}^n\|_{\Omega_f}^2+\|\varepsilon_{\phi,h}^n\|_{\Omega_m}^2+\|(\mathbf{e}_{f,b}^n-\bar{\partial}_t \mathbf{e}_{m,b}^n)\cdot\bm{\tau}_{\Gamma}\|^2_{\Gamma}\bigr)+\trb{\mathbf{e}_{m,h}^{\bar n}}_m^2\nonumber\\
    	\leq\ &Ch^{2k}\mathcal{R}_h(\bar{n})+C\ell^2\mathcal{R}_{\ell}(\bar{n}),
    \end{align}
    where
    \begin{align*}
    	\mathcal{R}_h(\bar{n}):=\ &\max_{1\le n\le \bar n}\|\mathbf{u}_f^n\|^2_{k+1}+\max_{1\le n\le \bar n}\|p_f^n\|^2_k+\max_{1\le n\le \bar n}\|p_m^n\|^2_{k+1}+\max_{1\le n\le \bar n}\|\partial_t \phi_m^n\|^2_k+\max_{1\le n\le \bar n}\|\mathbf{u}_m^n\|^2_{k+1}\nonumber\\
    	&+\max_{1\le n\le \bar n}\|\phi_m^n\|^2_k+\int_0^{t^{\bar{n}}}(\|\partial_t \phi_m(r)\|^2_k+\|\partial_t \mathbf{u}_m(r)\|^2_{k+1})\,\mathrm{d}r,\\
    	\mathcal{R}_{\ell}(\bar{n}):=\ &\int_0^{t^{\bar{n}}}(\|\partial_{tt}\mathbf{u}_m(r)\|^2_{\Gamma}+\|\partial_{tt}p_m(r)\|^2_{\Omega_m}+\|\partial_{tt}\phi_m(r)\|^2_{\Omega_m})\,\mathrm{d}r.
    \end{align*}
\end{theorem}
\begin{proof}
    Applying the backward difference operator to \eqref{full-error-equation5}, we obtain
    \begin{align}
    	\frac{1}{\lambda}(\bar{\partial}_t \varepsilon_{\phi,h}^n, \psi_{m,h})-\frac{\alpha}{\lambda}(\bar{\partial}_t \varepsilon_{m,0}^n, \psi_{m,h})+b_m(\bar{\partial}_t \mathbf{e}_{m,h}^n, \psi_{m,h}) =0.\label{full-error-equation6}
    \end{align}
    Taking $
    \mathbf{v}_{f,h}=\mathbf{e}_{f,h}^n,\,q_{f,h}=\varepsilon_{f,h}^n,\,\mathbf{v}_{m,h}=\bar{\partial}_t\mathbf{e}_{m,h}^n,\,q_{m,h}=\varepsilon_{m,h}^n$ in \eqref{full-error-equation1}--\eqref{full-error-equation4}, and taking $\psi_{m,h}=\varepsilon_{\phi,h}^n$ in \eqref{full-error-equation6}. Adding the resulting identities gives
    \begin{align}
    	&2\mu_f\trb{\mathbf{e}_{f,h}^n}_f^2+\frac{1}{\ell}a_m(\mathbf{e}_{m,h}^n,\mathbf{e}_{m,h}^n-\mathbf{e}_{m,h}^{n-1})+\frac{1}{\lambda\ell}((\alpha\varepsilon_{m,0}^n-\varepsilon_{\phi,h}^n)-(\alpha\varepsilon_{m,0}^{n-1}-\varepsilon_{\phi,h}^{n-1}),\alpha\varepsilon_{m,0}^n-\varepsilon_{\phi,h}^n)\nonumber\\
    	&\qquad+\kappa\trb{\varepsilon_{m,h}^n}_p^2+\frac{c_0}{\ell}(\varepsilon_{m,0}^n-\varepsilon_{m,0}^{n-1},\varepsilon_{m,0}^n)+\beta_{\Gamma}\|(\mathbf{e}_{f,b}^n-\bar{\partial}_t \mathbf{e}_{m,b}^n)\cdot\bm{\tau}_{\Gamma}\|^2_{\Gamma}=s_f(\mathbf{Q}_h \bu_f^n,\mathbf{e}_{f,h}^n)\nonumber\\
    	&\qquad-l_{\mathbf{u},f}(\bu_f^n,\mathbf{e}_{f,h}^n)+l_{p,f}(p_f^n,\mathbf{e}_{f,h}^n)+l_{\Gamma}^1(\mathbf{u}_m^n,\mathbf{e}_{f,h}^n)+s_m(\mathbf{Q}_h \mathbf{u}_m^n,\bar{\partial}_t \mathbf{e}_{m,h}^n)-l_{\mathbf{u},m}(\mathbf{u}_m^n,\bar{\partial}_t \mathbf{e}_{m,h}^n)\nonumber\\
    	&\qquad+l_{\phi,m}(\phi_m^n,\bar{\partial}_t \mathbf{e}_{m,h}^n)-l^2_{\Gamma}(\mathbf{u}_m^n,\bar{\partial}_t \mathbf{e}_{m,h}^n)+s_p(Q_h p_m^n,\varepsilon_{m,h}^n)-l_{\phi_t,m}((\partial_t \phi_m)^n,\varepsilon_{m,h}^n)-l_p(p_m^n,\varepsilon_{m,h}^n)\nonumber\\
    	&\qquad-l_0^1(p_m^n,\varepsilon_{m,0}^n)+l_0^2(\phi_m^n,\varepsilon_{m,0}^n)-l_{\Gamma}^3(\mathbf{u}_m^n,\varepsilon_{m,h}^n).\label{full-error-estimate-proof1}
    \end{align}
    \textcolor{black}{Let $I_{\rm lhs}$ represent the left-hand side of \eqref{full-error-estimate-proof1}. With the help of the inequality $\frac{1}{2}(\|a\|^2-\|b\|^2)\leq(a-b,a)$}, we know that 
    \begin{align*}
    	I_{\rm lhs}\geq\ & 2\mu_f\trb{\mathbf{e}_{f,h}^n}_f^2+\frac{\mu_m}{\ell}\bigl(\trb{\mathbf{e}_{m,h}^n}^2_m-\trb{\mathbf{e}_{m,h}^{n-1}}_m^2\bigr)+\frac{1}{2\lambda\ell}\bigl(\|\alpha\varepsilon_{m,0}^n-\varepsilon_{\phi,h}^n\|^2_{\Omega_m}-\|\alpha\varepsilon_{m,0}^{n-1}-\varepsilon_{\phi,h}^{n-1}\|^2_{\Omega_m}\bigr)\\
    	&+\kappa\trb{\varepsilon_{m,h}^n}_p^2+\frac{c_0}{2\ell}\bigl(\|\varepsilon_{m,0}^n\|^2_{\Omega_m}-\|\varepsilon_{m,0}^{n-1}\|^2_{\Omega_m}\bigr)+\beta_{\Gamma}\|(\mathbf{e}_{f,b}^n-\bar{\partial}_t \mathbf{e}_{m,b}^n)\cdot\bm{\tau}_{\Gamma}\|^2_{\Gamma}.
    \end{align*}
    \textcolor{black}{Substituting the above estimate into \eqref{full-error-estimate-proof1}. Then we sum the resulting inequality over $n=1,\ldots,\bar n$, and multiply by $\ell$ to get}
    \begin{align}
    	&\ell\sum_{n=1}^{\bar{n}}\bigl(2\mu_f\trb{\mathbf{e}_{f,h}^n}_f^2+\kappa\trb{\varepsilon_{m,h}^n}_p^2+\beta_{\Gamma}\|(\mathbf{e}_{f,b}^n-\bar{\partial}_t \mathbf{e}_{m,b}^n)\cdot\bm{\tau}_{\Gamma}\|^2_{\Gamma}\bigr)+\mu_m\trb{\mathbf{e}_{m,h}^{\bar{n}}}_m^2+\frac{1}{2\lambda }\|\alpha\varepsilon_{m,0}^{\bar{n}}-\varepsilon_{\phi,h}^{\bar{n}}\|_{\Omega_m}^2\nonumber\\
    	&\qquad+\frac{c_0}{2}\|\varepsilon_{m,0}^{\bar{n}}\|_{\Omega_m}^2\leq \ell\sum_{n=1}^{\bar{n}}\biggl(s_f(\mathbf{Q}_h \bu_f^n,\mathbf{e}_{f,h}^n)-l_{\mathbf{u},f}(\bu_f^n,\mathbf{e}_{f,h}^n)+l_{p,f}(p_f^n,\mathbf{e}_{f,h}^n)+l_{\Gamma}^1(\mathbf{u}_m^n,\mathbf{e}_{f,h}^n)\nonumber\\
    	&\qquad+s_p(Q_h p_m^n,\varepsilon_{m,h}^n)-l_{\phi_t,m}((\partial_t \phi_m)^n,\varepsilon_{m,h}^n)-l_p(p_m^n,\varepsilon_{m,h}^n)-l_0^1(p_m^n,\varepsilon_{m,0}^n)+l_0^2(\phi_m^n,\varepsilon_{m,0}^n)\nonumber\\
    	&\qquad-l_{\Gamma}^3(\mathbf{u}_m^n,\varepsilon_{m,h}^n)+\frac{1}{\ell}(s_m(\mathbf{Q}_h \mathbf{u}_m^n,\mathbf{e}_{m,h}^n-\mathbf{e}_{m,h}^{n-1})-l_{\mathbf{u},m}(\mathbf{u}_m^n,\mathbf{e}_{m,h}^n-\mathbf{e}_{m,h}^{n-1})+l_{\phi,m}(\phi_m^n,\mathbf{e}_{m,h}^n-\mathbf{e}_{m,h}^{n-1})\nonumber\\
    	&\qquad-l^2_{\Gamma}(\mathbf{u}_m^n,\mathbf{e}_{m,h}^n-\mathbf{e}_{m,h}^{n-1}))\biggr),\label{full-error-estimate-proof2}
    \end{align}
    where we have used $\mathbf{e}_{m,h}^0=\mathbf{0},\,\varepsilon_{m,h}^0=0$ and $\varepsilon_{\phi,h}^0=0$.

    \textcolor{black}{To bound the terms of the right-hand side of \eqref{full-error-estimate-proof2}, we use \eqref{error-term1} and Young's inequality,}
    \begin{align}
    	\ell\sum_{n=1}^{\bar{n}}s_f(\mathbf{Q}_h \bu_f^n,\mathbf{e}_{f,h}^n)&\leq \ell\sum_{n=1}^{\bar{n}}(Ch^k\|\mathbf{u}_f^n\|_{k+1}\trb{\mathbf{e}_{f,h}^n}_f)\nonumber\\
    	&\leq \ell\sum_{n=1}^{\bar{n}}(Ch^{2k}\|\mathbf{u}_f^n\|^2_{k+1}+\frac{2\mu_f}{3}\trb{\mathbf{e}_{f,h}^n}_f^2)\nonumber\\
    	&\leq Ch^{2k}\max_{1\le n\le \bar n}\|\mathbf{u}_f^n\|^2_{k+1}+\frac{\ell}{3}\sum_{n=1}^{\bar{n}}2\mu_f\trb{\mathbf{e}_{f,h}^n}_f^2.\label{full-error-estimate-proof-s_s}
    \end{align}
    Repeating the same argument as \eqref{full-error-estimate-proof-s_s} and using \eqref{error-term2}--\eqref{error-term6} gives
    \begin{align}
   	    &\ell\sum_{n=1}^{\bar{n}}l_{\mathbf{u},f}(\bu_f^n,\mathbf{e}_{f,h}^n)\leq Ch^{2k}\max_{1\le n\le \bar n}\|\mathbf{u}_f^n\|^2_{k+1}+\frac{\ell}{3}\sum_{n=1}^{\bar{n}}2\mu_f\trb{\mathbf{e}_{f,h}^n}_f^2,\label{full-error-estimate-proof-l_s}\\
   	    &\ell\sum_{n=1}^{\bar{n}}l_{p,f}(p_f^n,\mathbf{e}_{f,h}^n)\leq Ch^{2k}\max_{1\le n\le \bar n}\|p_f^n\|^2_k+\frac{\ell}{3}\sum_{n=1}^{\bar{n}}2\mu_f\trb{\mathbf{e}_{f,h}^n}_f^2,\label{full-error-estimate-proof-l_sp}\\
   	    &\ell\sum_{n=1}^{\bar{n}}s_p(Q_h p_m^n,\varepsilon_{m,h}^n)\leq Ch^{2k}\max_{1\le n\le \bar n}\|p_m^n\|^2_{k+1}+\frac{\ell}{6}\sum_{n=1}^{\bar{n}}\kappa\trb{\varepsilon_{m,h}^n}_p^2,\label{full-error-estimate-proof-s_p}\\
   	    &\ell\sum_{n=1}^{\bar{n}}l_{\phi_t,m}(\partial_t \phi_m^n,\varepsilon_{m,h}^n)\leq Ch^{2k}\max_{1\le n\le \bar n}\|\partial_t \phi_m^n\|_k^2+\frac{\ell}{6}\sum_{n=1}^{\bar{n}}\kappa\trb{\varepsilon_{m,h}^n}_p^2,\label{full-error-estimate-proof-l_phi}\\
   	    &\ell\sum_{n=1}^{\bar{n}}l_p(p_m^n,\varepsilon_{m,h}^n)\leq Ch^{2k}\max_{1\le n\le \bar n}\|p_m^n\|^2_{k+1}+\frac{\ell}{6}\sum_{n=1}^{\bar{n}}\kappa\trb{\varepsilon_{m,h}^n}_p^2,\label{full-error-estimate-proof-l_w}
    \end{align}
    Next, employing the initial condition $\mathbf{e}_{m,h}^0=\mathbf{0}$, together with \eqref{error-term1} to get
    \begin{align}
    	&\sum_{n=1}^{\bar{n}}s_m(\mathbf{Q}_h \mathbf{u}_m^n,\mathbf{e}_{m,h}^n-\mathbf{e}_{m,h}^{n-1})\nonumber\\
    	=\ &\sum_{n=1}^{\bar{n}}\bigl(s_m(\mathbf{Q}_h \mathbf{u}_m^n,\mathbf{e}_{m,h}^n)-s_m(\mathbf{Q}_h \mathbf{u}_m^{n-1},\mathbf{e}_{m,h}^{n-1})-s_m(\mathbf{Q}_h(\mathbf{u}_m^n-\mathbf{u}_m^{n-1}),\mathbf{e}_{m,h}^{n-1})\bigr)\nonumber\\
    	=\ &s_m(\mathbf{Q}_h \mathbf{u}_m^{\bar{n}},\mathbf{e}_{m,h}^{\bar{n}})-\sum_{n=1}^{\bar{n}}\int_{t^{n-1}}^{t^n}s_m(\mathbf{Q}_h (\partial_t \mathbf{u}_m(r)),\mathbf{e}_{m,h}^{n-1})\,\mathrm{d}r\nonumber\\
    	\leq\ & Ch^k\|\mathbf{u}_m^{\bar{n}}\|_{k+1}\trb{\mathbf{e}_{m,h}^{\bar{n}}}_m+\sum_{n=1}^{\bar{n}}\int_{t^{n-1}}^{t^n}Ch^k\|\partial_t \mathbf{u}_m(r)\|_{k+1}\trb{\mathbf{e}_{m,h}^{n-1}}_m\,\mathrm{d}r\nonumber\\
    	\leq\ & Ch^k\|\mathbf{u}_m^{\bar{n}}\|_{k+1}\trb{\mathbf{e}_{m,h}^{\bar{n}}}_m+\sum_{n=1}^{\bar{n}}\left(\int_{t^{n-1}}^{t^n}(Ch^k\|\partial_t \mathbf{u}_m(r)\|_{k+1})^2\,\mathrm{d}r\right)^{\frac{1}{2}}\left(\int_{t^{n-1}}^{t^n}\trb{\mathbf{e}_{m,h}^{n-1}}^2_m\,\mathrm{d}r\right)^{\frac{1}{2}}\nonumber\\
    	\leq\ & Ch^{2k}\|\mathbf{u}_m^{\bar{n}}\|^2_{k+1}+\frac{\mu_m}{6}\trb{\mathbf{e}_{m,h}^{\bar{n}}}_m^2+Ch^{2k}\int_0^{t^{\bar{n}}}\|\partial_t \mathbf{u}_m(r)\|^2_{k+1}\,\mathrm{d}r+\frac{\ell}{6}\sum_{n=1}^{\bar{n}}\mu_m\trb{\mathbf{e}_{m,h}^n}^2_m.\label{full-error-estimate-proof-s_m}
    \end{align}
    Arguing as in \eqref{full-error-estimate-proof-s_m}, and using \eqref{error-term2} and \eqref{error-term-phi}, we have
    \begin{align}
    	&\sum_{n=1}^{\bar{n}}l_{\mathbf{u},m}(\mathbf{u}_m^n,\mathbf{e}_{m,h}^n-\mathbf{e}_{m,h}^{n-1})\nonumber\\
    	\leq\ & Ch^{2k}\|\mathbf{u}_m^{\bar{n}}\|^2_{k+1}+\frac{\mu_m}{6}\trb{\mathbf{e}_{m,h}^{\bar{n}}}_m^2+Ch^{2k}\int_0^{t^{\bar{n}}}\|\partial_t \mathbf{u}_m(r)\|^2_{k+1}\,\mathrm{d}r+\frac{\ell}{6}\sum_{n=1}^{\bar{n}}\mu_m\trb{\mathbf{e}_{m,h}^n}_m^2,\label{full-error-estimate-proof-l_b}
    \end{align}
    and
    \begin{align}
    	&\sum_{n=1}^{\bar{n}}l_{\phi,m}(\phi_m^n,\mathbf{e}_{m,h}^n-\mathbf{e}_{m,h}^{n-1})\nonumber\\
    	\leq\ & Ch^{2k}\|\phi_m^{\bar{n}}\|^2_k+\frac{\mu_m}{6}\trb{\mathbf{e}_{m,h}^{\bar{n}}}_m^2+Ch^{2k}\int_0^{t^{\bar{n}}}\|\partial_t \phi_m(r)\|^2_k\,\mathrm{d}r+\frac{\ell}{6}\sum_{n=1}^{\bar{n}}\mu_m\trb{\mathbf{e}_{m,h}^n}_m^2.\label{full-error-estimate-proof-l_bphi}	
    \end{align}
    
    For any sufficiently smooth function $\phi$, we know 
    \begin{align}
    	\|(\partial_t\phi)^n-\bar{\partial}_t\phi^n\|_D^2
    	\le
    	C\ell\int_{t^{n-1}}^{t^n}
    	\|\partial_{tt}\phi(r)\|_D^2\,\mathrm{d}r,
    	\label{full-error-estimate-proof-truncation}
    \end{align}
    where $D$ denotes either $\Omega_m$ or $\Gamma$. Using Young's inequality, \eqref{full-error-estimate-proof-truncation}, \eqref{v0-vh1}, and the \textcolor{black}{parameter} assumptions $\lambda\in[\lambda_0,\infty)$ and $0\le c_0\le c_1$, we bound the term $l_0^1$ as:
    \begin{align}
    	\ell\sum_{n=1}^{\bar{n}}l_0^1(p_m^n,\varepsilon_{m,0}^n)&=\ell\sum_{n=1}^{\bar{n}}\Bigl(c_0+\frac{\alpha^2}{\lambda}\Bigr)((\partial_t p_m)^n-\bar{\partial}_t p_m^n,\varepsilon_{m,0}^n)_{\Omega_m}\nonumber\\
    	&\leq \ell\sum_{n=1}^{\bar{n}}\Bigl(c_0+\frac{\alpha^2}{\lambda}\Bigr)\|(\partial_t p_m)^n-\bar{\partial}_t p_m^n\|_{\Omega_m}\|\varepsilon_{m,0}^n\|_{\Omega_m}\nonumber\\
    	&\leq C\ell\sum_{n=1}^{\bar{n}}\Bigl(c_0+\frac{\alpha^2}{\lambda}\Bigr)\left(\ell\int_{t^{n-1}}^{t^n}\|\partial_{tt}p_m\|^2_{\Omega_m}\,\mathrm{d}r\right)^{\frac{1}{2}}\|\varepsilon_{m,0}^n\|_{\Omega_m}\nonumber\\
    	&\leq C\Bigl(c_0+\frac{\alpha^2}{\lambda}\Bigr)^2\ell^2\int_{0}^{t^{\bar{n}}}\|\partial_{tt}p_m\|^2_{\Omega_m}\,\mathrm{d}r+C\ell\sum_{n=1}^{\bar{n}}\kappa\|\varepsilon_{m,0}^n\|^2_{\Omega_m}\nonumber\\
    	&\leq C\ell^2\int_{0}^{t^{\bar{n}}}\|\partial_{tt}p_m\|^2_{\Omega_m}\,\mathrm{d}r+\frac{\ell}{6}\sum_{n=1}^{\bar{n}}\kappa\trb{\varepsilon_{m,h}^n}_p^2.\label{full-error-estimate-proof-l_01}
    \end{align}
    Applying the same argument as in \eqref{full-error-estimate-proof-l_01} together with \eqref{vb-vh1} yields
    \begin{align}
    	\ell\sum_{n=1}^{\bar{n}}l_0^2(\phi_m^n,\varepsilon_{m,0}^n)\leq C\ell^2\int_{0}^{t^{\bar{n}}}\|\partial_{tt}\phi_m\|^2_{\Omega_m}\,\mathrm{d}r+\frac{\ell}{6}\sum_{n=1}^{\bar{n}}\kappa\trb{\varepsilon_{m,h}^n}_p^2,\label{full-error-estimate-proof-l_02}
    \end{align}
    and
    \begin{align}
    	\ell\sum_{n=1}^{\bar{n}}l_{\Gamma}^3(\mathbf{u}_m^n, \varepsilon_{m,h}^n)\leq C\ell^2\int_0^{t^{\bar{n}}}\|\partial_{tt}\mathbf{u}_m(r)\|^2_{\Gamma}\,\mathrm{d}r+\frac{\ell}{6}\sum_{n=1}^{\bar{n}}\kappa\trb{\varepsilon_{m,h}^n}_p^2,\label{full-error-estimate-proof-l_gamma3}
    \end{align}
    The two tangential interface remainder terms are estimated together as follows:
    \begin{align}
    	&\ell\sum_{n=1}^{\bar{n}}l_{\Gamma}^1(\mathbf{u}_m^n, \mathbf{e}_{f,h}^n)-\ell\sum_{n=1}^{\bar{n}}l_{\Gamma}^2(\mathbf{u}_m^n, \bar{\partial}_t \mathbf{e}_{m,h}^n)\nonumber\\
    	=\ &\ell\sum_{n=1}^{\bar{n}}\beta_{\Gamma}\langle((\partial_t \mathbf{u}_m)^n-\bar{\partial}_t \mathbf{u}_m^n)\cdot\bm{\tau}_{\Gamma},(\mathbf{e}_{f,b}^n-\bar{\partial}_t \mathbf{e}_{m,b}^n)\cdot\bm{\tau}_{\Gamma}\rangle_{\Gamma},\nonumber\\
    	\leq\ & C\ell^2\int_0^{t^{\bar{n}}}\|\partial_{tt}\mathbf{u}_m(r)\|^2_{\Gamma}\,\mathrm{d}r+\ell\sum_{n=1}^{\bar{n}}\beta_{\Gamma}\|(\mathbf{e}_{f,b}^n-\bar{\partial}_t \mathbf{e}_{m,b}^n)\cdot\bm{\tau}_{\Gamma}\|^2_{\Gamma},\label{full-error-estimate-proof-l_gamma12}
    \end{align}
    We combine the above estimates \eqref{full-error-estimate-proof-s_s}--\eqref{full-error-estimate-proof-l_gamma12} with \eqref{full-error-estimate-proof2} to obtain
    \begin{align*}
    	\trb{\mathbf{e}_{m,h}^{\bar{n}}}_m^2\leq C\mathcal{R}_h(\bar{n})h^{2k}+C\mathcal{R}_{\ell}(\bar{n})\ell^2+\ell\sum_{n=1}^{\bar{n}}\trb{\mathbf{e}_{m,h}^n}^2_m.
    \end{align*}
    An application of the discrete Gronwall inequality \cite{MR4799262-2025-WG-Biot} then yields
    \begin{align}
    	\trb{\mathbf{e}_{m,h}^{\bar{n}}}_m^2\leq C\mathcal{R}_h(\bar{n})h^{2k}+C\mathcal{R}_{\ell}(\bar{n})\ell^2.\label{full-error-estimate-proof15}
    \end{align}
    
    Employing Young's inequality with different coefficients in \eqref{full-error-estimate-proof-s_s}--\eqref{full-error-estimate-proof-l_gamma12}, we have
    \begin{align*}
    	\ell\sum_{n=1}^{\bar{n}}\bigl(\trb{\mathbf{e}_{f,h}^n}_f^2+\trb{\varepsilon_{m,h}^n}_p^2+\|(\mathbf{e}_{f,b}^n-\bar{\partial}_t \mathbf{e}_{m,b}^n)\cdot\bm{\tau}_{\Gamma}\|^2_{\Gamma}\bigr)\leq C\mathcal{R}_h(\bar{n})h^{2k}+C\mathcal{R}_{\ell}(\bar{n})\ell^2+C\max\limits_{1\le j\le \bar n}\trb{\mathbf{e}_{m,h}^j}_m^2,
    \end{align*}
    Since \eqref{full-error-estimate-proof15} holds for every $1\le j\le \bar n$, we know
    \begin{align}
    	\max_{1\le j\le \bar n}\trb{\mathbf e_{m,h}^j}_m^2 \le C\mathcal{R}_h(\bar n)h^{2k}+C\mathcal{R}_{\ell}(\bar n)\ell^2,\label{full-error-estimate-proof16}
    \end{align}
    and  consequently,  
    \begin{align}
    	\ell\sum_{n=1}^{\bar{n}}\bigl(\trb{\mathbf{e}_{f,h}^n}_f^2+\trb{\varepsilon_{m,h}^n}_p^2+\|(\mathbf{e}_{f,b}^n-\bar{\partial}_t \mathbf{e}_{m,b}^n)\cdot\bm{\tau}_{\Gamma}\|^2_{\Gamma}\bigr)\leq C\mathcal{R}_h(\bar{n})h^{2k}+C\mathcal{R}_{\ell}(\bar{n})\ell^2.\label{full-error-estimate-proof17}
    \end{align}
   
    It remains to estimate errors of the pressure and the total pressure. From \eqref{full-error-equation1}, the Cauchy--Schwarz inequality, \eqref{vb-vh2}, \eqref{error-term1}--\eqref{error-term3}, and \eqref{full-error-estimate-proof-truncation}, we obtain
    \begin{align*}
    	b_f(\mathbf{v}_{f,h}, \varepsilon_{f,h}^n)=\ &a_f\left(\mathbf{e}_{f,h}^n, \mathbf{v}_{f,h}\right)+a_{\Gamma}(\mathbf{e}_{f,h}^n,\bar{\partial}_t \mathbf{e}_{m,h}^n;\mathbf{v}_{f,h}, \mathbf{0})+b_{\Gamma}(\mathbf{0}, \varepsilon_{m,h}^n;\mathbf{v}_{f,h},0)-s_f(\mathbf{Q}_h \mathbf{u}_f^n,\mathbf{v}_{f,h})\\
    	&+l_{\mathbf{u},f}(\mathbf{u}_f^n, \mathbf{v}_{f,h})-l_{s, p}(\mathbf{v}_{f,h}, p_f^n)-l_{\Gamma}^1(\mathbf{u}_m^n,\mathbf{v}_{f,h})\\
    	\leq\ &C\trb{\mathbf{e}_{f,h}^n}_f\trb{\mathbf{v}_{f,h}}_f+C\|(\mathbf{e}_{f,b}^n-\bar{\partial}_t \mathbf{e}_{m,b}^n)\cdot\bm{\tau}_{\Gamma}\|_{\Gamma}\trb{\mathbf{v}_{f,h}}_f+C\trb{\mathbf{v}_{f,h}}_f\trb{\varepsilon_{m,h}^n}_p\\
    	&+Ch^k\|\mathbf{u}_f^n\|_{k+1}\trb{\mathbf{v}_{f,h}}_f+Ch^k\|p_f^n\|_k\trb{\mathbf{v}_{f,h}}_f+C\ell^{\frac{1}{2}}\left(\int_{t^{n-1}}^{t^n}\|\partial_{tt}\mathbf{u}_m(r)\|^2_{\Gamma}\,\mathrm{d}r\right)^{\frac{1}{2}}\trb{\mathbf{v}_{f,h}}_f.
    \end{align*}
    Using the inf-sup condition \eqref{inf-sup-Stokes} to get
    \begin{align*}
    	\|\varepsilon_{f,h}^n\|_{\Omega_f}^2\leq\ & C\biggl(\trb{\mathbf{e}_{f,h}^n}_f^2+\|(\mathbf{e}_{f,b}^n-\bar{\partial}_t \mathbf{e}_{m,b}^n)\cdot\bm{\tau}_{\Gamma}\|_{\Gamma}^2+\trb{\varepsilon_{m,h}^n}_p^2+h^{2k}(\|\mathbf{u}_f^n\|_{k+1}^2+\|p_f^n\|_k^2)\\
    	&+\ell\int_{t^{n-1}}^{t^n}\|\partial_{tt}\mathbf{u}_m(r)\|^2_{\Gamma}\,\mathrm{d}r\biggr).
    \end{align*}
    \textcolor{black}{Multiplying the above estimate by $\ell$, summing it over $n=1,\ldots,\bar n$, and applying \eqref{full-error-estimate-proof17}, we have}
    \begin{align}
    	\ell\sum_{n=1}^{\bar{n}}\|\varepsilon_{f,h}^n\|_{\Omega_f}^2\leq\ & C\ell\sum_{n=1}^{\bar{n}}(\trb{\mathbf{e}_{f,h}^n}_f^2+\trb{\varepsilon_{m,h}^n}_p^2+\|(\mathbf{e}_{f,b}^n-\bar{\partial}_t \mathbf{e}_{m,b}^n)\cdot\bm{\tau}_{\Gamma}\|^2_{\Gamma})+Ch^{2k}\max_{1\le n\le \bar n}\|\mathbf{u}_f^n\|_{k+1}^2\nonumber\\
    	&+Ch^{2k}\max_{1\le n\le \bar n}\|p_f^n\|_k^2+C\ell^2\int_{0}^{t^{\bar{n}}}\|\partial_{tt}\mathbf{u}_m(r)\|^2_{\Gamma}\,\mathrm{d}r\nonumber\\
    	\leq\ &C\mathcal{R}_h(\bar{n})h^{2k}+C\mathcal{R}_{\ell}(\bar{n})\ell^2.\label{full-error-estimate-proof18}
    \end{align}
    Finally, following a similar procedure as in \eqref{full-error-estimate-proof18}, and using \eqref{full-error-estimate-proof16},
    \begin{align*}
    	\ell\sum_{n=1}^{\bar{n}}\|\varepsilon_{\phi,h}^n\|_{\Omega_m}^2\leq C\mathcal{R}_h(\bar{n})h^{2k}+C\mathcal{R}_{\ell}(\bar{n})\ell^2.
    \end{align*}
    \textcolor{black}{Therefore, the desired error estimates are achieved.}
\end{proof}

\section{Numerical results}\label{Section:example}
\begin{figure}[!htbp]
	\centering
	\includegraphics[width=0.37\linewidth]{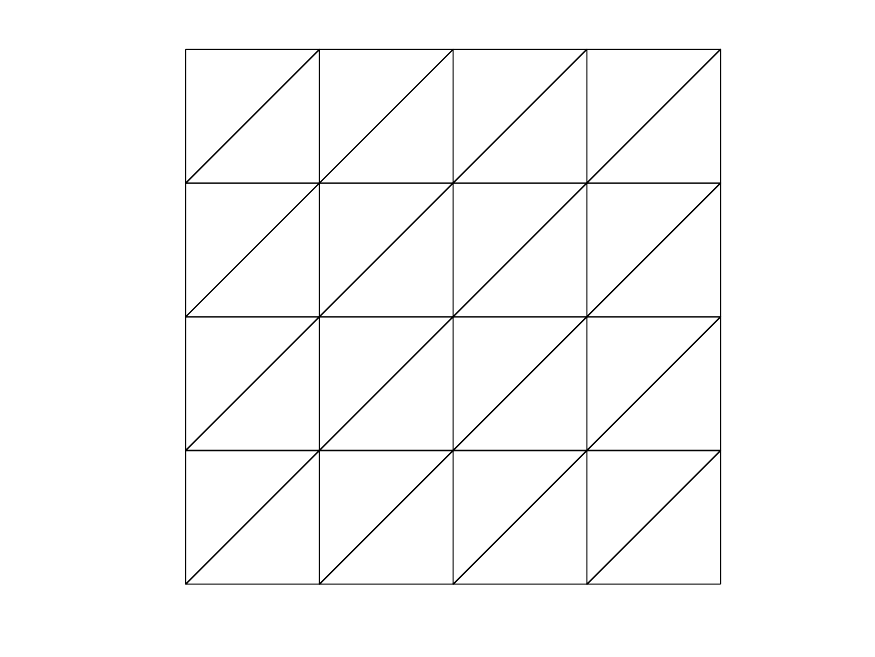}
	\hspace{-0.08\linewidth}
	\includegraphics[width=0.37\linewidth]{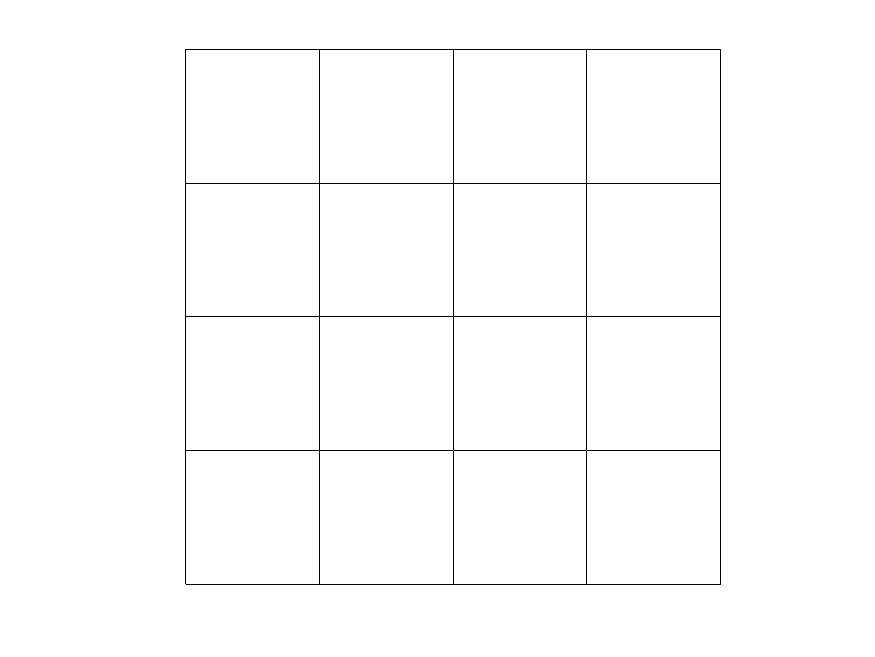}
	\hspace{-0.08\linewidth}
	\includegraphics[width=0.37\linewidth]{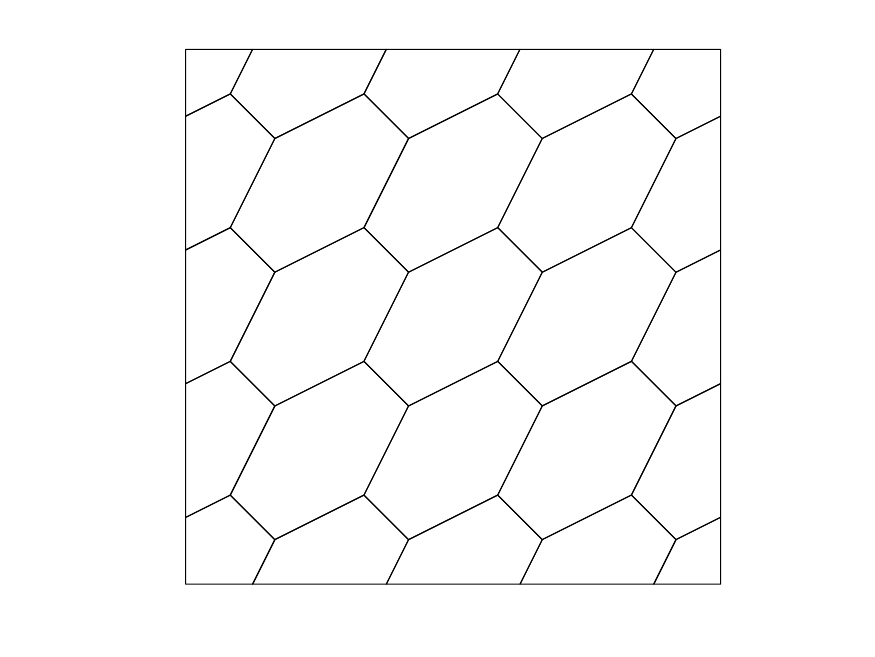}
	\caption{The first-level meshes of the three types of meshes: triangular (left), quadrilateral (middle), and hexagonal (right).}\label{mesh}
\end{figure}

\begin{figure}[!htbp]
	\centering
	\includegraphics[width=0.48\linewidth]{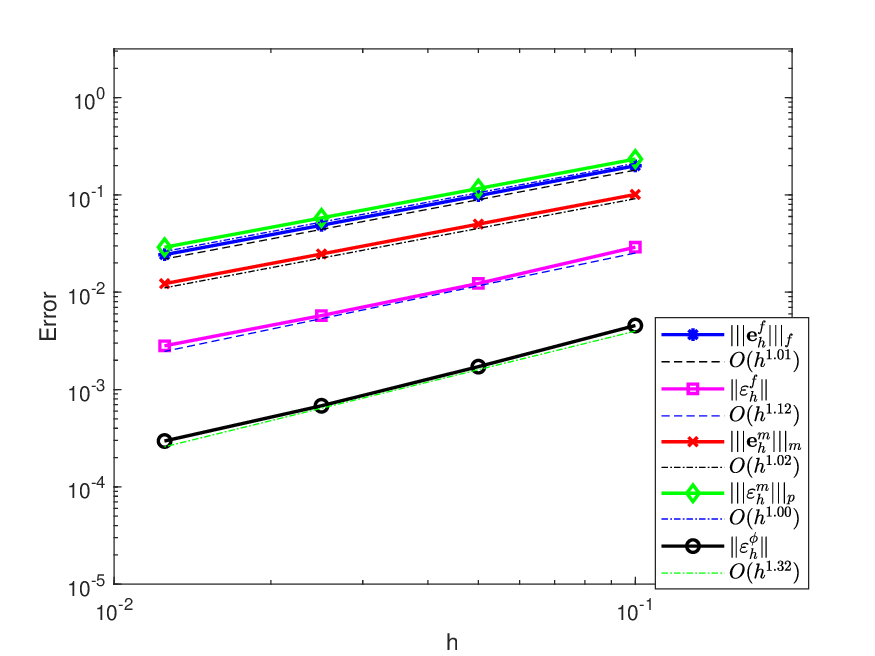}
	\includegraphics[width=0.48\linewidth]{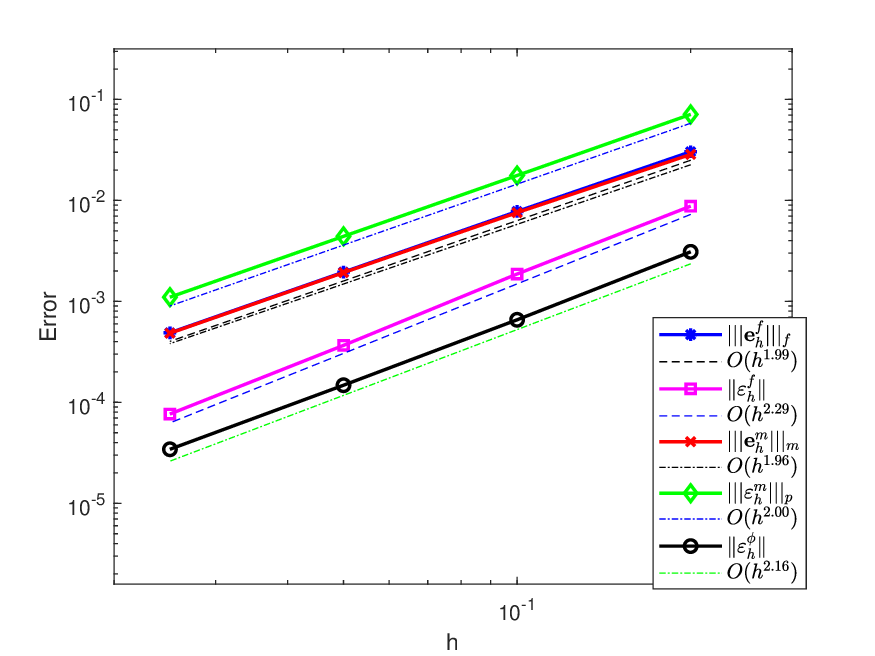}
	\caption{The numerical results of Example~\ref{numexample1} on the triangular meshes with $k = 1$ (left) and with $k = 2$ (right).}\label{triangle}
\end{figure}

\begin{figure}[!htbp]
	\centering
	\includegraphics[width=0.48\linewidth]{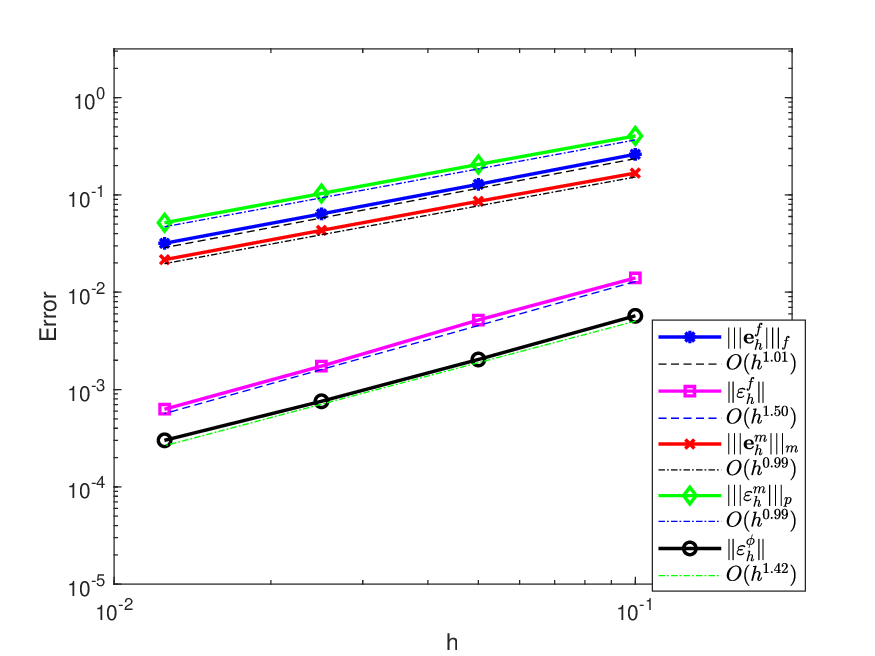}
	\includegraphics[width=0.48\linewidth]{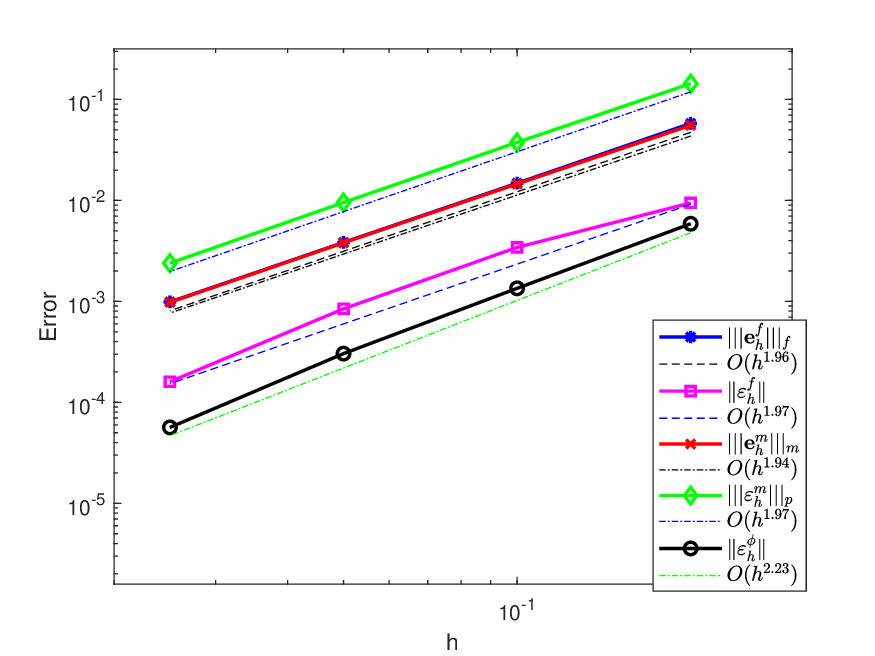}
	\caption{The numerical results of Example~\ref{numexample1} on the rectangular meshes with $k = 1$ (left) and with $k = 2$ (right).}\label{rectangle}
\end{figure}

\begin{figure}[!htbp]
	\centering
	\includegraphics[width=0.48\linewidth]{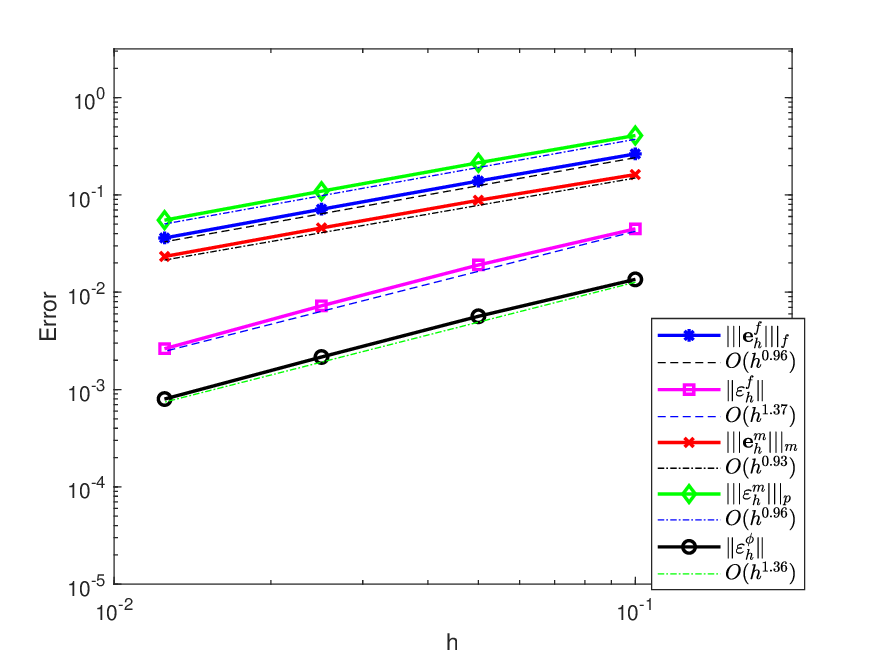}
	\includegraphics[width=0.48\linewidth]{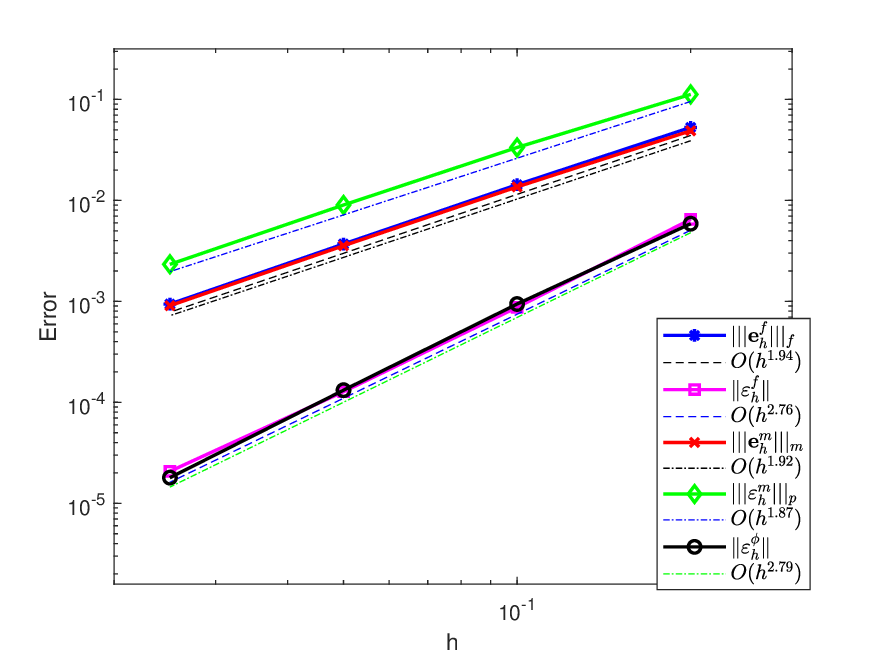}
	\caption{The numerical results of Example~\ref{numexample1} on the hexagonal meshes with $k = 1$ (left) and with $k = 2$ (right).}\label{polygon}
\end{figure}

In this section, three examples are presented to demonstrate the accuracy and robustness of the WG algorithm. 
\begin{example}\label{numexample1}\cite{MR4362198-2022-Esaim-mixed}
{\color{black}In this example, we test the convergence of the algorithm. The computational domain is taken as} $\Omega=(0,1)\times(-1,1)$ with $\Omega_f=(0,1)\times(0,1),\ \Omega_m=(0,1)\times(-1,0)$ and the interface $\Gamma=(0,1)\times\{0\}$. Let $\Gamma_f^{N}=(0,1)\times\{1\}$, $\Gamma_f^D=\Gamma_f\setminus\Gamma_f^{N}$, $\Gamma_m^{N}=\Gamma_m^{N_1}=(0,1)\times\{-1\}$, $\Gamma_m^{D}=\Gamma_m^{D_1}=\Gamma_m\setminus\Gamma_m^{N}$. We consider the time interval $J=[0,0.1]$. The exact solutions are
	\begin{align*}
		\mathbf{u}_f&=\Bigg(\begin{array}{c}
			\pi\cos(\pi t)(-3x+\cos y)\\
			\pi\cos(\pi t)(y+1)
		\end{array}\Bigg),& &p_f=e^t\sin(\pi x)\cos\Bigl(\frac{1}{2}\pi y\Bigr)+2\pi\cos(\pi t),\\
    	\mathbf{u}_m&=\Bigg(\begin{array}{c}
    		\sin(\pi t)(-3x+\cos y)\\
    		\sin(\pi t)(y+1)
    	\end{array}\Bigg),& &p_m=e^t\sin(\pi x)\cos\Bigl(\frac{1}{2}\pi y\Bigr).
    \end{align*}
    The total pressure is determined by
    $\phi_m=\alpha p_m-\lambda\nabla\cdot\mathbf u_m$.
\end{example}
We {\color{black}choose} the parameters as $\mu_f=1,\,\mu_m=1,\,c_0=1,\,\alpha=1,\,\lambda=1,\,\kappa=1$, and $\alpha_{\rm BJS}=1$. {\color{black}The numerical algorithm is implemented on three types of meshes (See Figure~\ref{mesh})}. The corresponding results with $k=1$ and $k=2$ are presented in Figures~\ref{triangle}--Figures~\ref{polygon}. {\color{black}As we can see, all variables achieve} the optimal convergence rates in the corresponding norms, in agreement with the theoretical error estimates in Section~\ref{section:error-estimate}.

\begin{example}\label{numexample2}\cite{MR4799262-2025-WG-Biot,MR4126527-Stokes-Darcy-WG-FEM}
	In this example, we check the locking-free property of {\color{black}the proposed method}. Let $\Omega=(0,1)\times(0,2)$ with $\Omega_f=(0,1)\times(1,2),\ \Omega_m=(0,1)\times(0,1)$ and the interface $\Gamma=(0,1)\times\{1\}$. The boundaries are {\color{black}divided} as $\Gamma_f^N = (0,1)\times\{2\},\,\Gamma_f^D = \Gamma_f \setminus \Gamma_f^N,\,\Gamma_m^{N} = \Gamma_m^{N_1} = (0,1)\times\{0\},\,\Gamma_m^{D} = \Gamma_m^{D_1} = \Gamma_m \setminus \Gamma_m^{N}$. Let $T=1$ be the fine time. The exact solutions are
\begin{align*}
	\mathbf{u}_f&=\Bigg(\begin{array}{c}
		-t\cos(\pi x)\sin(\pi y)\\
		t\sin(\pi x)\cos(\pi y)
	\end{array}\Bigg),& &p_f=t\sin(\pi x),\\
	\mathbf{u}_m&=e^{-t}\Bigg(\begin{array}{c}
		\sin(2\pi y)(-1+\cos (2\pi x))+\dfrac{1}{\mu_m+\lambda}\sin(\pi x)\sin(\pi y)\\
		\sin(2\pi x)(1-\cos (2\pi y))+\dfrac{1}{\mu_m+\lambda}\sin(\pi x)\sin(\pi y)
	\end{array}\Bigg),& &p_m=e^{-t}\sin(\pi x)\sin(\pi y).
\end{align*}
The total pressure is determined by $\phi_m=\alpha p_m-\lambda\nabla\cdot\mathbf u_m$.
\end{example}

{\color{black}To observe the variation of numerical results with the Lam\'e parameter $\lambda$, we first fix the other physical parameters as follows:  $\mu_f = 1$, $\mu_m = 0.5$, $c_0 = 0$, $\alpha = 1$, $\kappa = 1$, and $\alpha_{\rm BJS} = 1$. Subsequently, we increase the value of $\lambda$ with $\lambda=$ 1, 100, 10000.}
Table~\ref{example2-k=1}--Table~\ref{example2-k=2} present the errors and the corresponding convergence rates employing triangular meshes with $k = 1$ and $k = 2$, respectively, {\color{black}which demonstrate that our method is locking-free.} 
\begin{table}[!htbp]
	\centering
	\caption{The errors and convergence rates in Example~\ref{numexample2} on triangular meshes with $k=1,\,\ell=h$.}\label{example2-k=1}
	\setlength{\belowcaptionskip}{-0.1cm}
	\footnotesize
	\begin{tabular}{ccccccccccc}
		\hline level & $\trb{\mathbf{e}_{f,h}}_f$ & Rate & $\|\varepsilon_{f,h}\|$ & Rate & $\trb{\mathbf{e}_{m,h}}_m$ & Rate & $\|\varepsilon_{m,h}\|$ & Rate & $\|\varepsilon_{\phi,h}\|$ & Rate\\
		\hline \multicolumn{11}{c}{$\lambda=1$}\\ 
		\hline
		2  & 6.4485E-01 &   & 2.6631E-01 &   &5.7845E-01  &  &1.2450E-01  & &1.0836E-01  & \\
		
		3  & 3.2226E-01 & 1.00  & 1.3129E-01 & 1.02  & 2.8658E-01 & 1.01 &6.0346E-02  &1.04 &5.3248E-02  &1.03  \\
		
		4  & 1.6096E-01 & 1.00  & 6.5300E-02 & 1.01  & 1.4262E-01 & 1.01 &2.9783E-02  & 1.02&2.6450E-02  & 1.01  \\
		
		5 & 8.0418E-02 & 1.00  & 3.2589E-02 & 1.00  & 7.1142E-02 & 1.00 &1.4810E-02  &1.01 &1.3211E-02  &1.00   \\
		
		\hline \multicolumn{11}{c}{$\lambda=100$}\\ 
		\hline
		2  & 6.4481E-01 &   & 2.6621E-01 &   &5.8054E-01  &  &1.2316E-01  & &1.8803E-01  & \\
		
		3  & 3.2226E-01 & 1.00  & 1.3126E-01 & 1.02  & 2.8818E-01 & 1.01 &6.0107E-02  &1.03 &8.8150E-02 &1.09  \\
		
		4  & 1.6096E-01 & 1.00  & 6.5294E-02 & 1.01  & 1.4348E-01 & 1.01 &2.9723E-02  & 1.02&4.3000E-02  & 1.04  \\
		
		5 & 8.0418E-02 & 1.00  & 3.2588E-02 & 1.00  & 7.1577E-02 & 1.00 &1.4786E-02  &1.01 &2.1365E-02  &1.01   \\
		\hline \multicolumn{11}{c}{$\lambda=10000$}\\ 
		\hline
		2  & 6.4481E-01 &   & 2.6621E-01 &   &5.8080E-01  &  &1.2315E-01  & &1.8958E-01  & \\
		
		3  & 3.2226E-01 & 1.00  & 1.3126E-01 & 1.02  & 2.8832E-01 & 1.01 &6.0106E-02  &1.03 &8.8765E-02 &1.09  \\
		
		4  & 1.6096E-01 & 1.00  & 6.5294E-02 & 1.01  & 1.4355E-01 & 1.01 &2.9723E-02  & 1.02&4.3275E-02  & 1.04  \\
		
		5 & 8.0418E-02 & 1.00  & 3.2588E-02 & 1.00  & 7.1611E-02 & 1.00 &1.4786E-02  &1.01 &2.1498E-02  &1.01   \\
		\hline
	\end{tabular}
\end{table}
\begin{table}[!htbp]
	\centering
	\caption{The errors and convergence rates in Example~\ref{numexample2} on triangular meshes with $k=2,\,\ell=h^2$.}\label{example2-k=2}
	\setlength{\belowcaptionskip}{-0.1cm}
	\footnotesize
	\begin{tabular}{ccccccccccc}
		\hline level & $\trb{\mathbf{e}_{f,h}}_f$ & Rate & $\|\varepsilon_{f,h}\|$ & Rate & $\trb{\mathbf{e}_{m,h}}_m$ & Rate & $\|\varepsilon_{m,h}\|$ & Rate & $\|\varepsilon_{\phi,h}\|$ & Rate\\
		\hline \multicolumn{11}{c}{$\lambda=1$}\\ 
		\hline
		1  & 1.9887E-01 &   & 9.0009E-02 &   &3.0570E-01  &  &4.8750E-02  & &6.8174E-02  & \\
		
		2  & 5.1193E-02 & 1.96  & 2.2056E-02 & 2.03  & 8.3495E-02 & 1.87 &1.1849E-02  &2.04 &1.9142E-02  &1.83  \\
		
		3  & 1.2962E-02 & 1.98  & 5.3862E-03 & 2.03  & 2.1646E-02 & 1.95 &2.9105E-03  & 2.03 &4.8828E-03  & 1.97  \\
		
		4 & 3.2575E-03 & 1.99  & 1.3283E-03 & 2.02  & 5.4883E-03 & 1.98 &7.2263E-04  &2.01 &1.2202E-03  &2.00   \\
		
		\hline \multicolumn{11}{c}{$\lambda=100$}\\ 
		\hline
		1  & 1.9884E-01 &   & 9.0020E-02 &   &3.0152E-01  &  &4.7864E-02  & &1.1391E-01  & \\
		
		2  & 5.1191E-02 & 1.96  & 2.2050E-02 & 2.03  & 8.2699E-02 & 1.87 &1.1744E-02  &2.03 &3.0918E-02 &1.88  \\
		
		3  & 1.2962E-02 & 1.98  & 5.3856E-03 & 2.03  & 2.1507E-02 & 1.94 &2.8983E-03  & 2.02 &7.7159E-03  & 2.00  \\
		
		4 & 3.2575E-03 & 1.99  & 1.3282E-03 & 2.02  & 5.4611E-03 & 1.98 &7.2054E-04  &2.01 &1.9064E-03  &2.02   \\
		\hline \multicolumn{11}{c}{$\lambda=10000$}\\ 
		\hline
		1  & 1.9884E-01 &   & 9.0020E-02 &   &3.0156E-01  &  &4.7855E-02  & &1.1471E-01  & \\
	
		2  & 5.1191E-02 & 1.96  & 2.2050E-02 & 2.03  & 8.2716E-02 & 1.87 &1.1743E-02  &2.03 &3.1116E-02 &1.88  \\
		
		3  & 1.2962E-02 & 1.98  & 5.3855E-03 & 2.03  & 2.1512E-02 & 1.94 &2.8982E-03  & 2.02 &7.7620E-03  & 2.00  \\
		
		4 & 3.2575E-03 & 1.99  & 1.3282E-03 & 2.02  & 5.4625E-03 & 1.98 &7.2054E-04  &2.01 &1.9174E-03  &2.02   \\
		\hline
	\end{tabular}
\end{table}
\begin{example}\label{numexample3}
    In this example, we consider the example without exact solutions constructed in \cite{MR4362198-2022-Esaim-mixed}.  Let $\Omega = (0,2)\times(-1,1)$, with $\Omega_f=(0,2)\times(0,1)$, $\Omega_m=(0,2)\times(-1,0)$ and the interface $\Gamma=(0,2)\times\{0\}$.   
    The body force and source terms are zero, and homogeneous initial conditions are prescribed. Let $\Gamma_f^D=\Gamma_f$, $\Gamma_m^{N}=\Gamma_m^{D_1}=(0,2)\times\{-1\}$, $\Gamma_m^{D}=\Gamma_m^{N_1}=\Gamma_m\setminus\Gamma_m^{N}$. Time interval is $[0, 3]$. The boundary conditions are
	\begin{align*}
		&\mathbf{u}_f=\Bigg(\begin{array}{c}
			-20y(y-1)(2-x)\\
			0
		\end{array}\Bigg),& & \text{on}~\Gamma_f^D, 
		\\
		&\mathbf{u}_m=\mathbf{0},\quad\ (-\kappa\nabla p_m)\cdot\mathbf{n}=0,& & \text{on}~\Gamma_m^{D}=\Gamma_m^{N_1}
		\\
		&p_m=0,\quad\ (2\mu_m\mathbf{D}(\mathbf{u}_m)-\phi_m\mathbf I)\mathbf{n}=0,& &\text{on}~\Gamma_m^{D_1}=\Gamma_m^{N}.
	\end{align*}
 The parameters $\mu_f$, $\alpha$ and $\alpha_{\rm BJS}$ are fixed to be unity. Different values of $\kappa$, $c_0$, $\mu_m$ and $\lambda$ are chosen in Table \ref{sen_params}.
\begin{table}[htbp]
  \centering
  \caption{Parameters for the sensitivity analysis.}
  \label{sen_params}
  \begin{tabular}{ccccc}
    \hline
    Case & $\kappa$ & $c_0$ & $\lambda$ & $\mu_m$ \\
    \hline
    I    & $1$      & $1$   & $1$       & $1$      \\
    II   & $10^{-4}$& $10^{-4}$ & $10^{6}$ & $1$    \\
    III  & $10^{-4}$& $10^{-4}$ & $10^{6}$ & $10^{6}$\\
    \hline
  \end{tabular}
\end{table}
\end{example}
 This example {\color{black}is computed on uniform triangular meshes with $k=1$}.
 The time step is taken as $\Delta \ell = 0.06$. Since there is no exact solution to this example, we evaluate the numerical results mainly through their dynamical behaviors. In Figure \ref{velocity-fig} - Figure \ref{pressure-fig}, we present the numerical results at $T=3$ with $h=1/64$ for three different cases. Here are some  observations.
\begin{enumerate}
    \item The numerical solution is free of Poisson locking or nonphysical pressure oscillations under extreme physical parameters (nearly incompressible, low permeability and storage coefficient).
    \item In the upper Stokes domain, the overall velocity patterns are similar among three cases. For Case II and Case III, the vertical velocity component near the coupling interface is evidently affected by pressure and solid deformation in the Biot subdomain, and the normal velocity ($\mathbf{n}_m = [0,~1]^T$, $\mathbf{u}\cdot\mathbf{n}_m = \mathbf{u}_2$) across the interface decreases.
    \item With simultaneous growth of Lam\'e constants $\mu_m$ and $\lambda$, the poroelastic pressure and stress increase significantly, while the solid displacement reduces greatly, compared with Case II where only $\lambda$ is large and $\mu_m$ remains small.
\end{enumerate}

\noindent
\begin{minipage}{\linewidth}
\captionsetup{type=figure}
	\centering
	\begin{subfigure}[b]{0.36\linewidth}
		\centering
		\includegraphics[width=\linewidth,trim=0 100 0 0,clip]{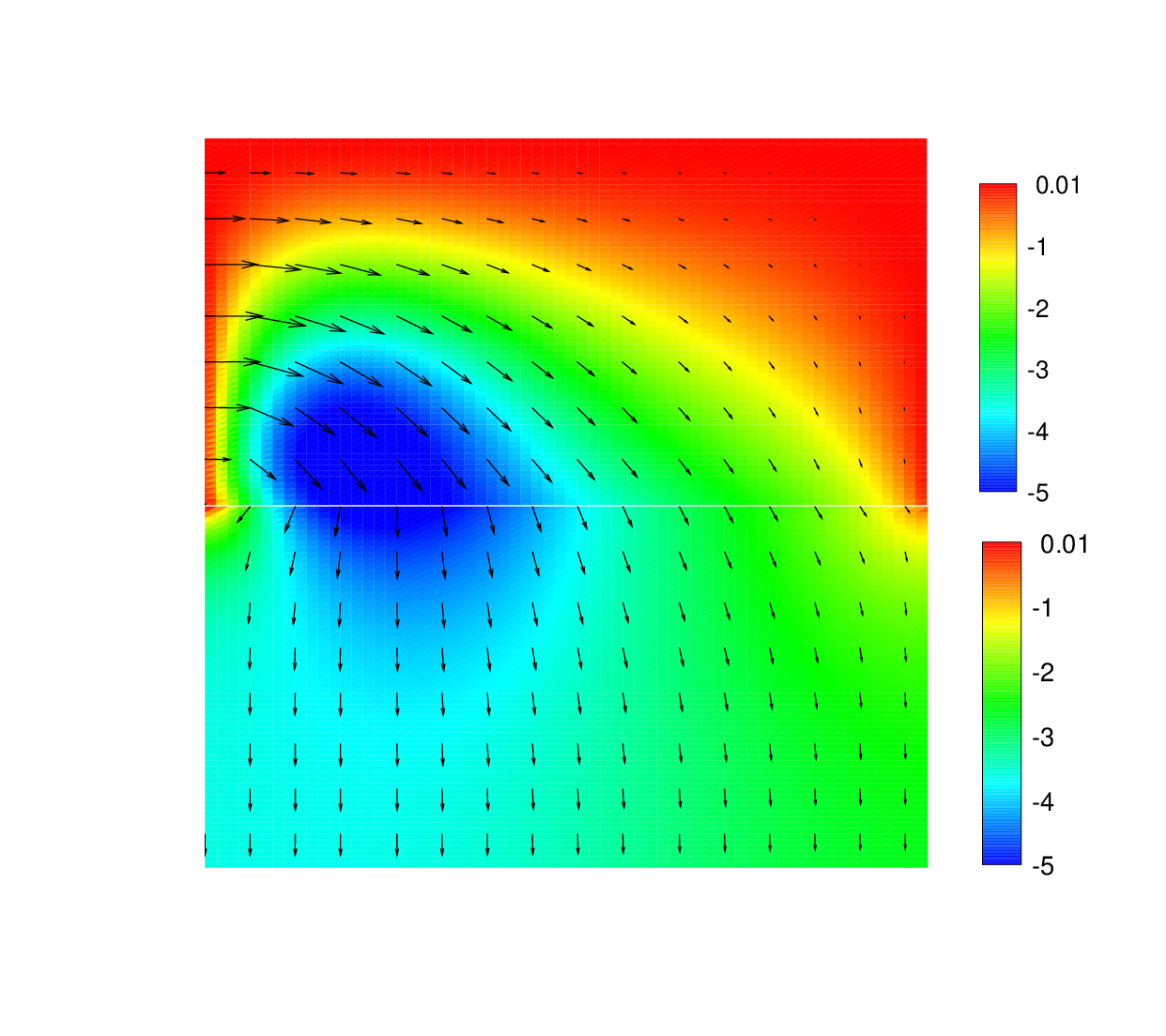}
		\caption*{ }
	\end{subfigure}
	\hspace{-0.055\linewidth}
	\begin{subfigure}[b]{0.36\linewidth}
		\centering
		\includegraphics[width=\linewidth,trim=0 100 0 0,clip]{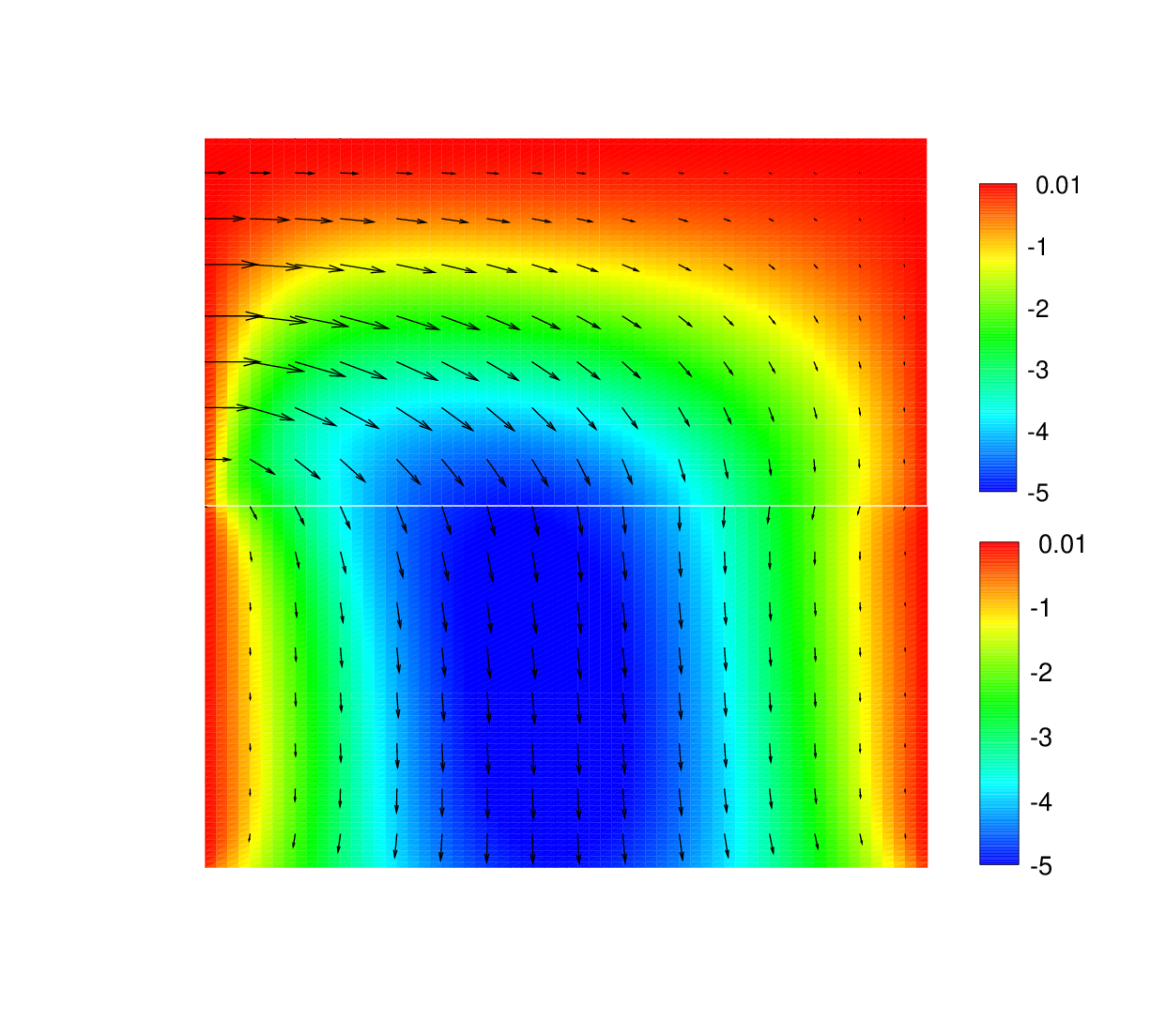}
		\caption*{ }
	\end{subfigure}
	\hspace{-0.055\linewidth}
	\begin{subfigure}[b]{0.36\linewidth}
		\centering
		\includegraphics[width=\linewidth,trim=0 100 0 0,clip]{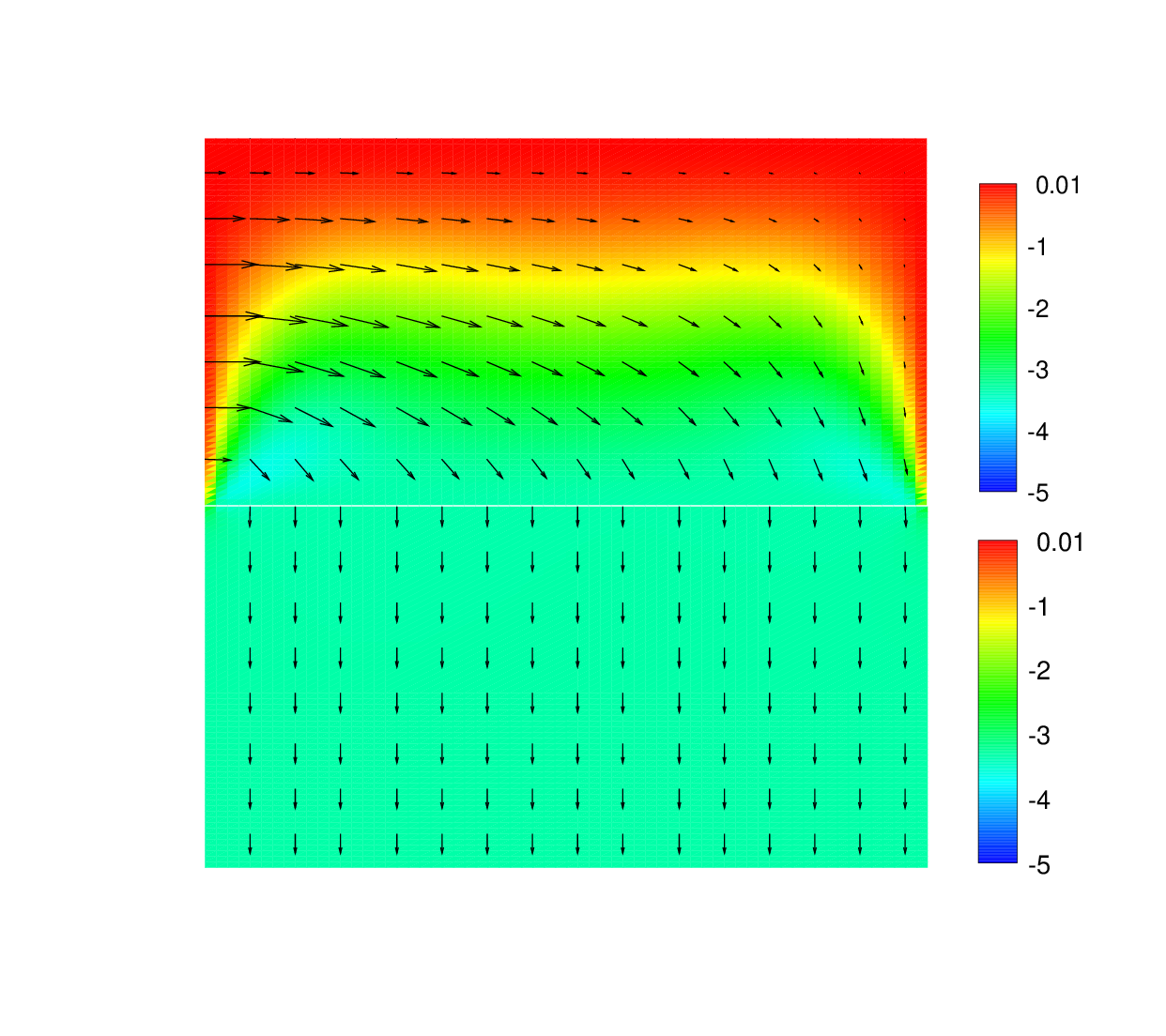}
		\caption*{ }
	\end{subfigure}
    \caption{Velocities  $\mathbf{u}_{f}$ and $\partial_t \mathbf{u}_m-\kappa\nabla p_m$ (arrow), $\mathbf{u}_{f,2}$ and $(\partial_t \mathbf{u}_m-\kappa\nabla p_m)_2$ (contour);  Case I (left), Case II (middle), Case III (right).}\label{velocity-fig}
\end{minipage}

\noindent
\begin{minipage}{\linewidth}
\captionsetup{type=figure}
	\centering
	\begin{subfigure}[b]{0.35\linewidth}
		\centering
		\includegraphics[width=\linewidth,trim=0 200 0 0,clip]{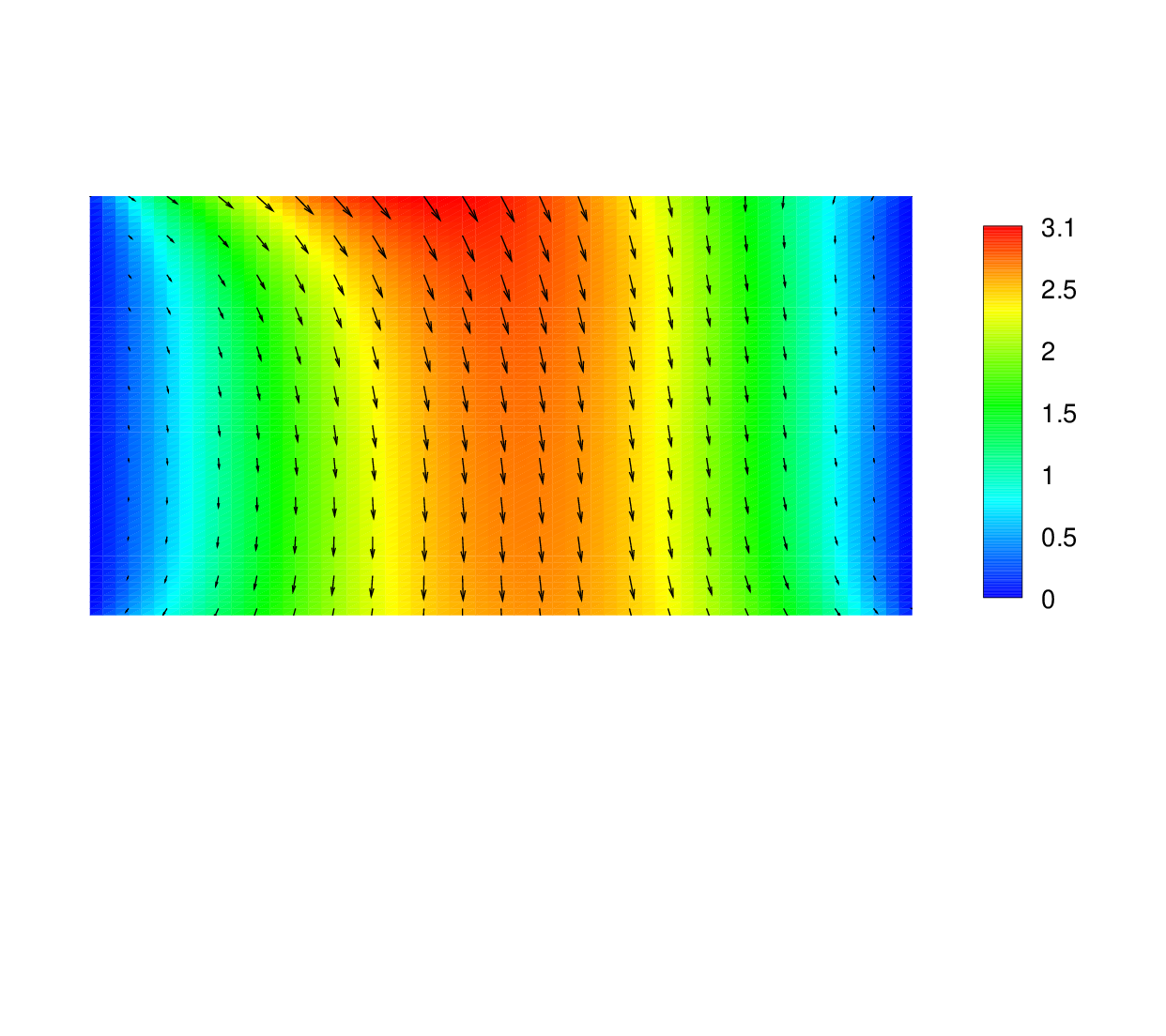}
		\caption*{ }
	\end{subfigure}
	\hspace{-0.055\linewidth}
	\begin{subfigure}[b]{0.35\linewidth}
		\centering
		\includegraphics[width=\linewidth,trim=0 200 0 0,clip]{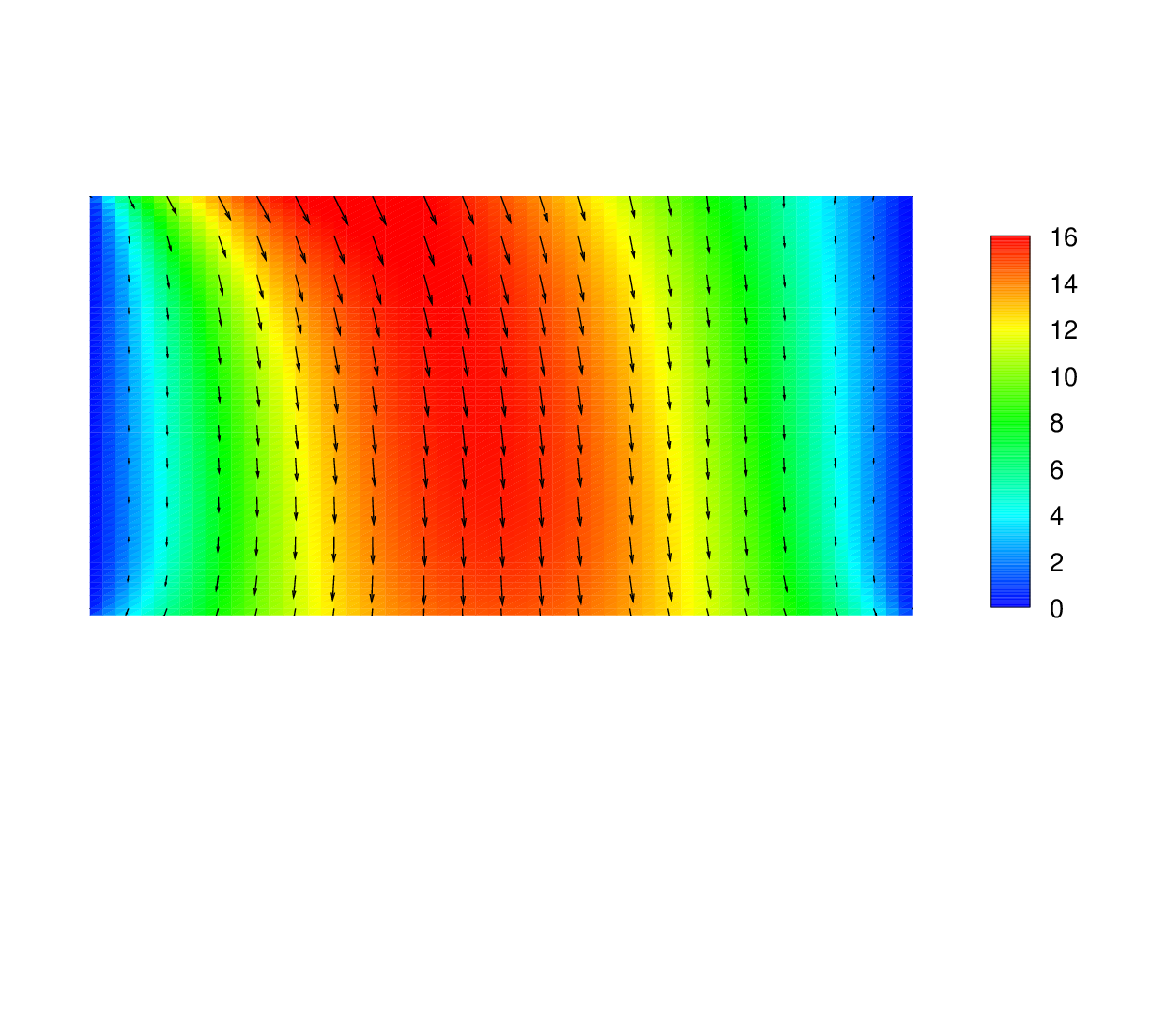}
		\caption*{ }
	\end{subfigure}
	\hspace{-0.055\linewidth}
	\begin{subfigure}[b]{0.35\linewidth}
		\centering
		\includegraphics[width=\linewidth,trim=0 200 0 0,clip]{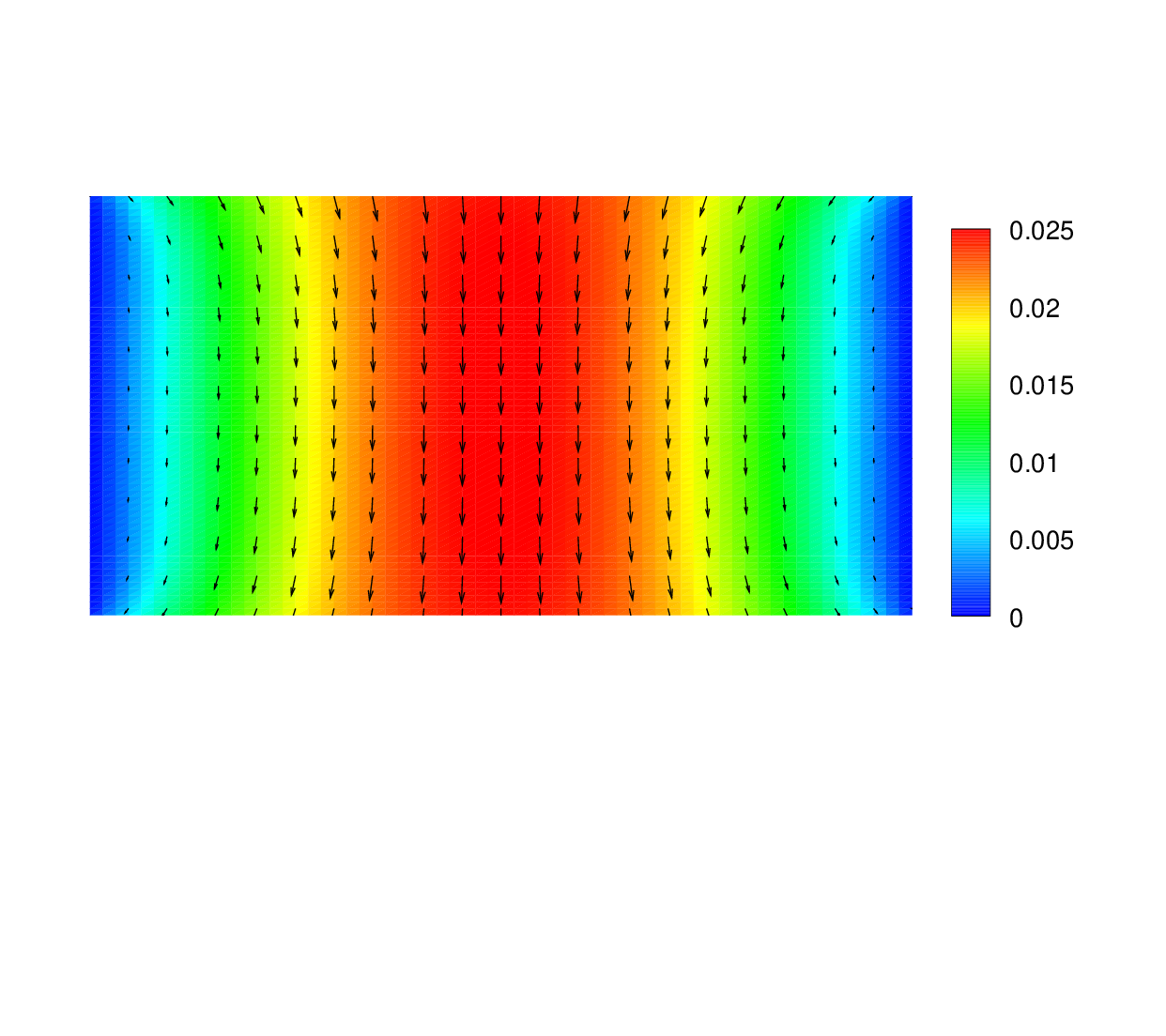}
		\caption*{ }
	\end{subfigure}
    \caption{Displacement $\mathbf{u}_m$ (arrow), $|\mathbf{u}_m|$ (contour); Case I (left), Case II (middle), Case III (right).}\label{displacement-fig}
\end{minipage}

\noindent
\begin{minipage}{\linewidth}
\captionsetup{type=figure}
	\centering
	\begin{subfigure}[b]{0.35\linewidth}
		\centering
		\includegraphics[width=\linewidth,trim=0 200 0 0,clip]{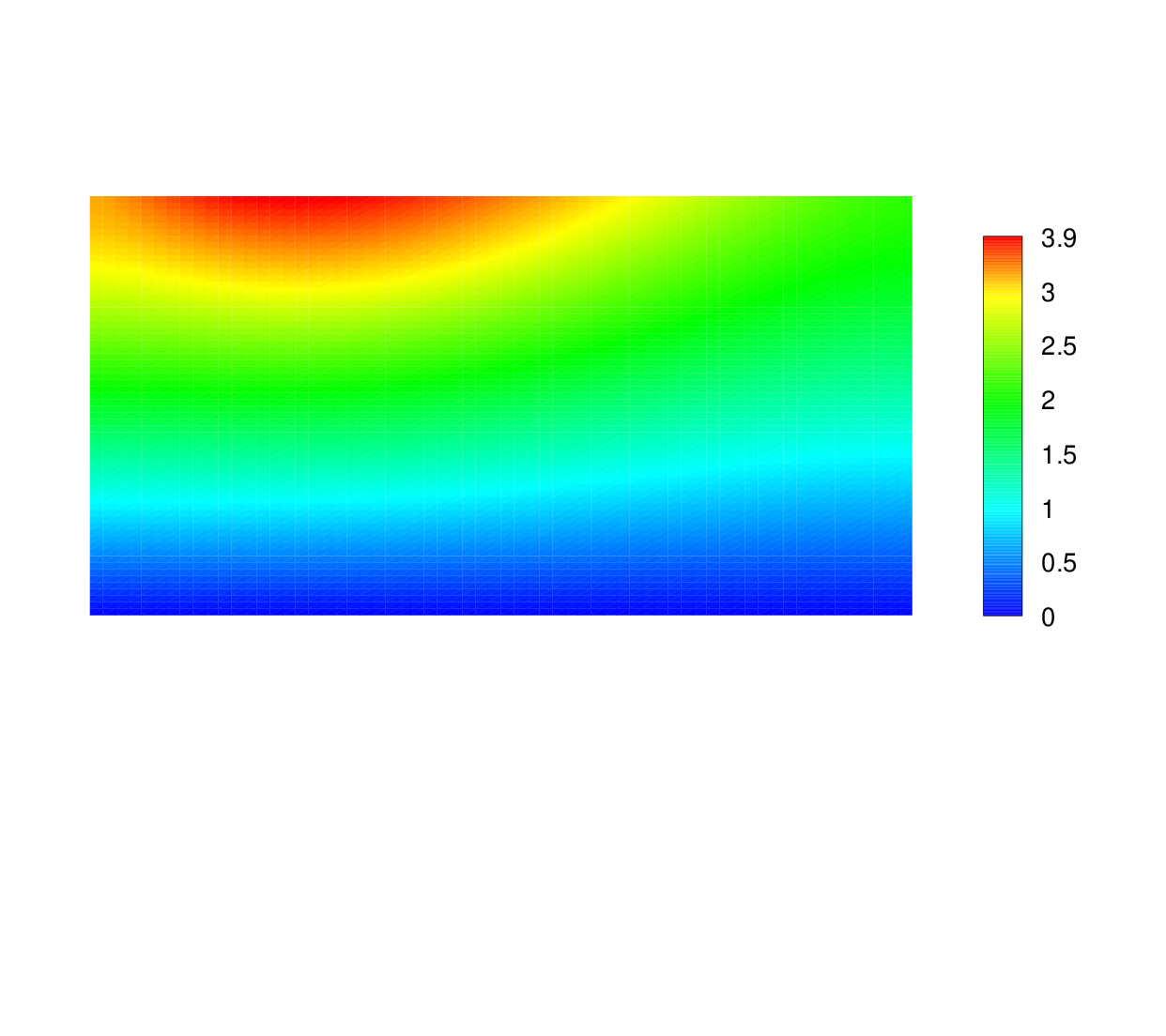}
		\caption*{ }
		\label{fig6:sub7}
	\end{subfigure}
	\hspace{-0.055\linewidth}
	\begin{subfigure}[b]{0.35\linewidth}
		\centering
		\includegraphics[width=\linewidth,trim=0 200 0 0,clip]{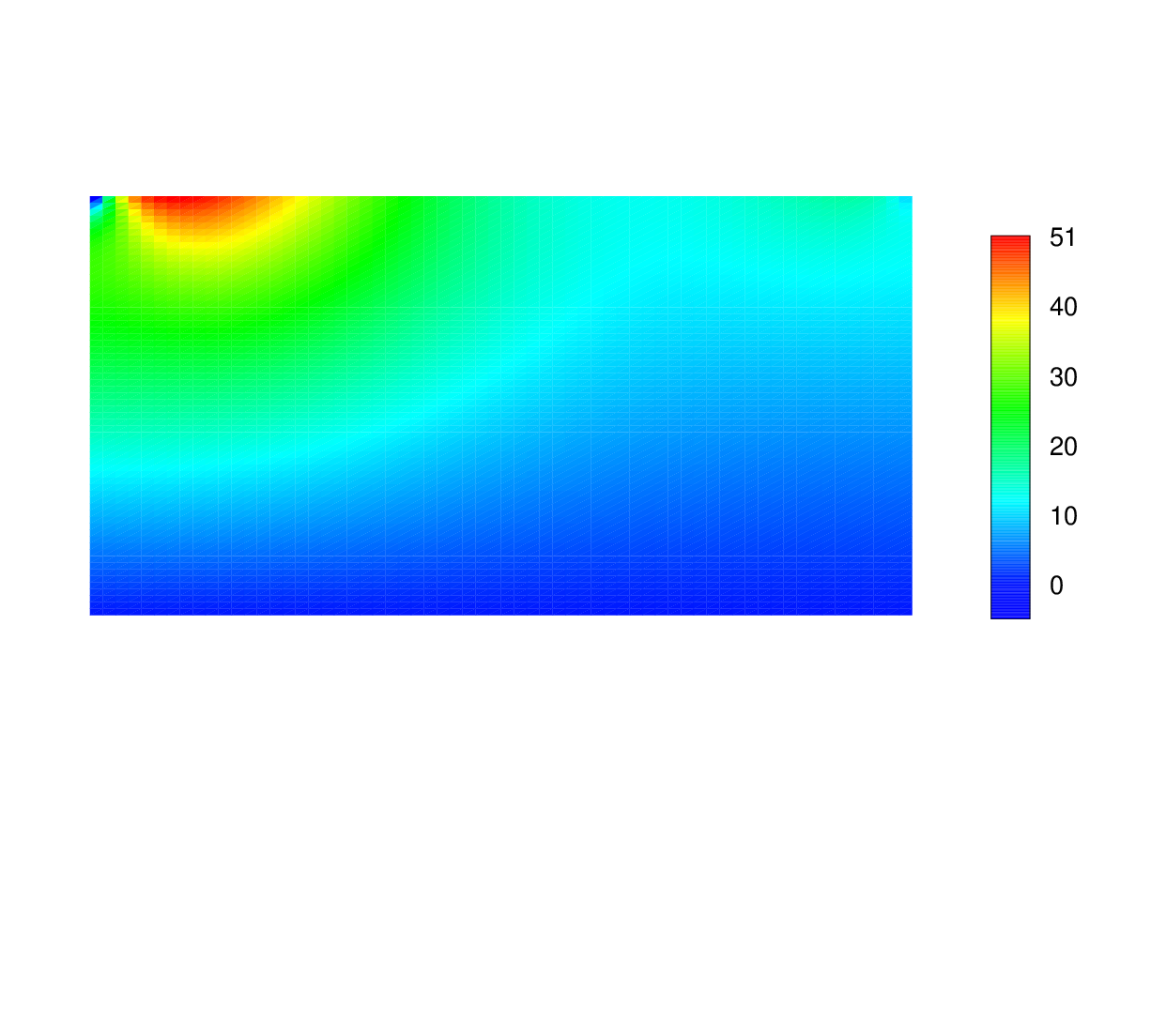}
		\caption*{ }
		\label{fig6:sub8}
	\end{subfigure}
	\hspace{-0.055\linewidth}
	\begin{subfigure}[b]{0.35\linewidth}
		\centering
		\includegraphics[width=\linewidth,trim=0 200 0 0,clip]{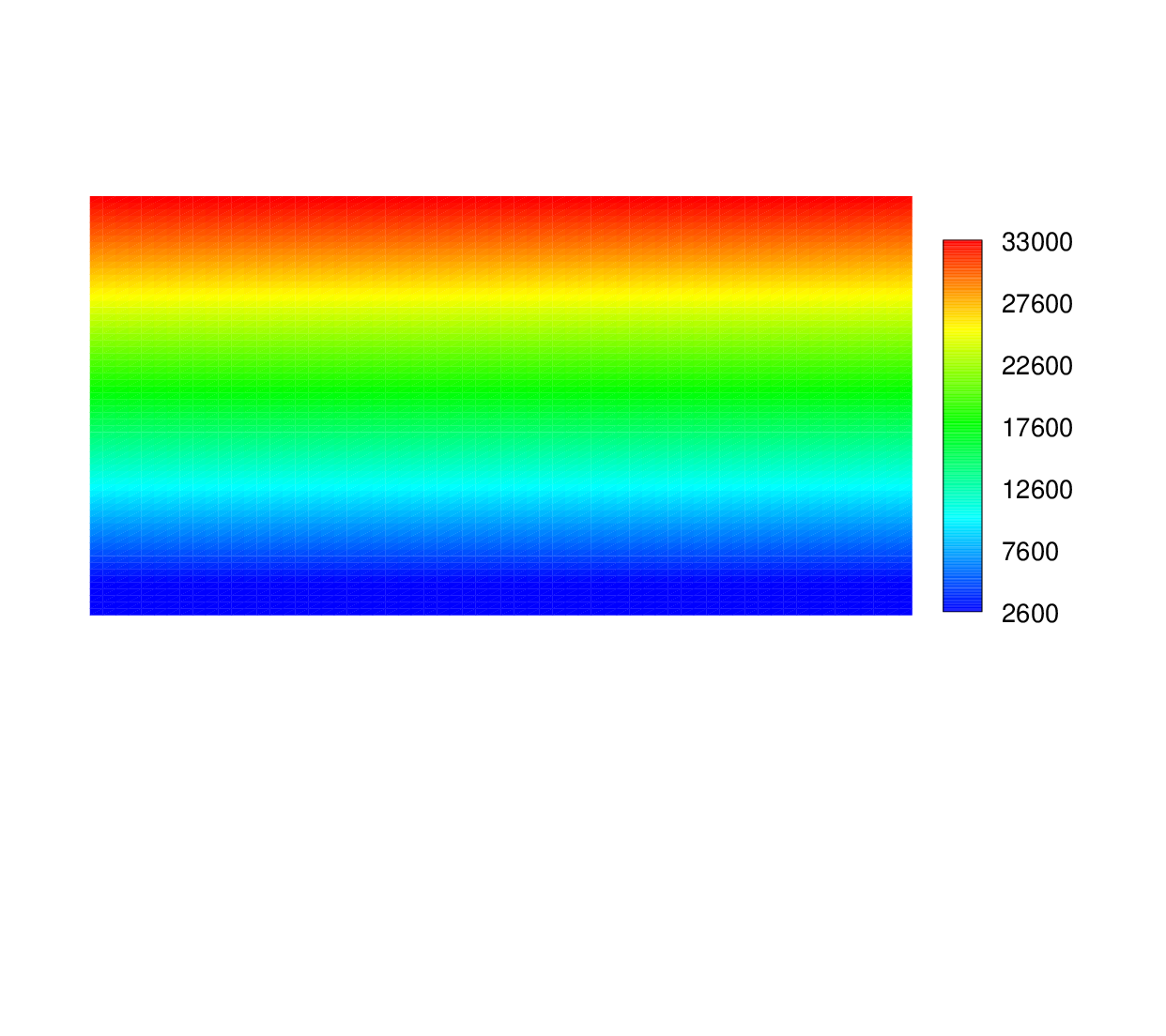}
		\caption*{ }
		\label{fig6:sub9}
	\end{subfigure}
    \caption{Poroelastic pressure ${p}_m$ (contour); Case I (left), Case II (middle), Case III (right).}\label{pressure-fig}
\end{minipage}

\section{Conclusion}\label{sec:7}
In this work, a stable numerical method is provided for Stokes--Biot model. Compatible weak Galerkin spaces are selected for spatial discretization and backward Euler method is employed to discrete time. The rigorously analysis of well-posedness and convergence are established for semi-discrete and fully discrete schemes. The proposed method supports polygonal mesh partitions  and exhibits stability respect to some physical parameters. 
{\color{black}However, the WG discretization leads to a large‑scale algebraic system. In future, we aim to design decoupled algorithms to reduce computational cost of the coupled Stokes--Biot model.}

%


\end{document}